\documentclass{amsart}

\usepackage{amsfonts,latexsym,amsmath,graphicx,psfrag}
\usepackage{amssymb}
\usepackage{hyperref}
\usepackage[left=3.7cm,right=3.7cm]{geometry}

\usepackage[usenames]{color}

\usepackage{graphicx}
\usepackage{graphics} 

\newtheorem{definition}{Definition}[section] 
\newtheorem{theorem}[definition]{Theorem} 
\newtheorem{lemma}[definition]{Lemma} 
 
\newtheorem{proposition}[definition]{Proposition} 
\newtheorem{corollary}[definition]{Corollary} 
\newtheorem{remark}[definition]{Remark}

\def\btext{} 

\numberwithin{equation}{section}

\newtheorem{predf}{Definition}[section]

\newtheorem{prex}[section]{Example}

\def\hatdN{D^a_Y}
\def\hatd{D^a_Y}
 \def\tilded{d}
 
 \def\YY{Y}

\def \Cmainone {C_4}
\def \Cmaintwo {C_5}
\def \Cmainthree {C_0}
\def \Cmainfour {C_1}
\def \Cmainfive {C_2}
\def \Cmainfivadded {C_3}

\def \CexpLip {C_{A}}
\def \CfirstToponogov {C_{6}}
\def \CsecondToponogov {C_7}
\def \CthirdToponogov {C_8}
\def \Calphafirst {C_9}
\def \Calphabeta {C_{10}}
\def \Caux {C'_{10}}
\def \Calphasecond {C_{11}}
\def \Cbetadisfirst {C_{12}}
\def \betadissecond{C_{13}}
\def \Cclosedirection {C_{14}}
\def \Cclosedistance {C_{15}}
\def \Cangleseparation {C_{16}}
\def \Cnetseparation {C_{17}}
\def \Cs {C_{18}}
\def \Cxiv {C_{19}}
\def \Csl {C'_{19}}
\def \CE {C_{20}}
\def \CEv {C_{21}}
\def \Cmainfivelast{C_{22}}
\def \Clast{C_{23}}

\def\crneighbor{\rho_0}
\def\cvdetlower{c_1}
\def\ccalphafirst{c_2}
\def\cvdetneighbor{c_3}
\def\cvdetwidth{c_4}

\def\epszero{\e_0}
\def\epsone{\e_1}
\def\hatepsone{{\hat\e}_1}
\def\epsX{\e_2}

\newcommand{\notmid}{\mid\kern-0.5em\not\kern0.5em}

\newcommand{\R}{\mathbb{R}}

\def\Z{{\mathbb Z}}

\newcommand{\ignore}[1]{}

\newcommand{\barr}{\begin{array}}
\newcommand{\earr}{\end{array}}

\def\bfo{\begin{eqnarray*}}
\def\efo{\end{eqnarray*}}
\def\ba{\begin{eqnarray*}}
\def\ea{\end{eqnarray*}}
\def\beq{\begin{eqnarray}}
\def\eeq{\end{eqnarray}}

\def\hat{\widehat}
\def\tilde{\widetilde}

\renewcommand{\H}{{\mathbb H}}

\def\dim{\hbox{dim}\,}

\def\dist{\hbox{dist}}

\def\det{\hbox{det}}

\def\bra{\langle}
\def\cet{\rangle}

\def\e{\varepsilon}
\def\p{\partial}
\def\a{\alpha}

\def\tilde{\widetilde}

\def\picture #1 by #2 (#3){
   \vsquare to #2{
     \hrule width #1 height 0pt depth 0pt
     \hfill
     \special{picture #3}}}

\def\scaledpicture #1 by #2 (#3 scaled #4){{
   \dimen0=#1 \dimen1=#2
   \divide\dimen0 by 1000 \multiply\dimen0 by #4
   \divide\dimen1 by 1000 \multiply\dimen1 by #4
   \picture \dimen0 by \dimen1 (#3 scaled #4)}}

\def \Z {{\Bbb {Z}}}

\def \R {{\Bbb {R}}}

\def \H2s {H^{s+1}_0(\partial M\times [0,T/2])}

\def \dist {\hbox{dist}}

\def \det {\hbox{det}}
\def\bra{\langle}
\def\cet{\rangle}
\def \e {\varepsilon}

\def \a {\alpha}

\def \pa0 {\partial _0}

\def \p {\partial}

\def\e{\varepsilon}

\def\tilde{\widetilde}
\def \bar{\overline }

\def \mbeq {\begin {eqnarray}}
\def \meeq {\end {eqnarray}}

\newcommand{\pd}{\partial}

\title[Reconstruction and interpolation of manifolds II]{Reconstruction and interpolation of manifolds II: \\
Inverse problems with partial data for distances observations and for the heat kernel}

\author[Fefferman, Ivanov, Lassas,  Lu,  Narayanan]
{Charles Fefferman, Sergei Ivanov,
Matti Lassas, \\ Jinpeng Lu, Hariharan Narayanan\\ \
\\
{Dedicated to the memory of Yaroslav V. Kurylev}}

\address{\hspace{0.1cm} \linebreak
Charles Fefferman, Princeton University, Mathematics Department,
Fine Hall, Washington Road,
Princeton NJ, 08544-1000, USA.\hspace{12cm} \vspace{3mm}
\linebreak
 Sergei Ivanov,
St.~Petersburg Department of Steklov Institute of Mathematics, Russian Academy of Sciences, 
27 Fontanka, 191023 St.~Petersburg, Russia.
\hspace{6cm} 
\vspace{3mm}
\linebreak
Matti Lassas,
University of Helsinki,
Department of Mathematics and Statistics, P.O. Box 68,
00014, Helsinki, Finland.
\hspace{6cm} 
\vspace{3mm}
\linebreak
Jinpeng Lu,
University of Helsinki,
Department of Mathematics and Statistics, P.O. Box 68,
00014, Helsinki, Finland.
\hspace{6cm} 
\vspace{3mm}
 \linebreak
Hariharan Narayanan, 
School of Technology and Computer Science, 
Tata Institute for Fundamental Research, 
Mumbai 400005, India.
 }

\begin{document}


\maketitle

\vspace{-8mm}



\begin{abstract}

We consider how a closed Riemannian manifold $M$ and its metric tensor $g$
can be approximately reconstructed from local distance measurements. Moreover, we consider
an inverse problem of determining $(M,g)$ from limited knowledge on the heat kernel. 
In the part 1 of the paper, we considered the approximate construction of a smooth
manifold in the case when one is given the noisy distances $\tilde d(x,y)=d(x,y)+\varepsilon_{x,y}$ for all points $x,y\in X$,  where $X$ is a $\delta$-dense subset of $M$ and 
$|\varepsilon_{x,y}|<\delta$.
 In this part 2 of the paper, we consider a similar problem with partial data, that is, 
the approximate construction of the manifold $(M,g)$ when we are given
$\tilde d(x,y)$ for $x\in X$ and $y \in U\cap X$, where $U$ is an open subset of $M$.
In addition, 
we consider the inverse problem of determining the manifold $(M,g)$ with non-negative
Ricci curvature from noisy observations of the heat kernel $G(y,z,t)$. We show
that a manifold approximating $(M,g)$ can be determined in a stable way, when for some unknown source points
$z_j$ in $X\setminus U$, we are given the values of the heat kernel
$G(y,z_k,t)$  for $y\in X\cap U$ and $t\in (0,1)$ with a multiplicative noise. We also give a uniqueness result for the inverse problem in the case when the data does not contain noise and consider applications in manifold learning. A novel feature of the inverse problem for the heat kernel is that the set $M\setminus U$ 
containing the sources and the observation set $U$  are disjoint.
\end{abstract}

\section{Introduction}

Let $(M,g)$ be a closed connected Riemannian manifold of dimension $n \geq 2$. We consider $(M,g)$ in the following class of Riemannian manifolds with bounded geometry given by
\begin{eqnarray}\label{basic 1}
\quad \hbox{diam}(M) \leq \Lambda,
\quad \hbox{inj}(M) \geq  \Lambda^{-1},
\quad |\hbox{Sec}_M| \leq \Lambda^2,
\end {eqnarray}
where $\Lambda\geq 1$, $\hbox{diam}(M)$ denotes the diameter of $M$, $\hbox{inj}(M)$ denotes the injectivity radius of $M$, and $\hbox{Sec}_M$ denotes the sectional curvature of $M$.


Let $U\subset M$ be an open subset, and assume that $U$ contains an open ball $B(x_0,R)$ of radius $R>\Lambda^{-1}$ centered at some $x_0\in U$. 
We call the subset $U$ the measurement domain.
We say that $Y\subset U$ is an
 $\varepsilon$-net in $U$ if the  $\varepsilon$ neighborhood of $Y$ in $M$ contains the set $U$.
We also say that a set $Y$ is $\varepsilon$-dense in $U$ when it is an $\varepsilon$-net in $U$.

The goal of this paper is to show that 
when  $Y$ is  an $\e$-net in $U$ and
$X$ is  an $\e$-net in $M$, then the approximate distances $d(x,y)$ between the points
$x\in X$  and $y\in Y$ determine an approximation of the whole manifold $M$.
The part 1 of the paper, \cite{FIKLN}, considers the case when the observation domain $U$ is the whole manifold $M$.

\begin{figure}[htbp]


\begin{center}
\includegraphics[width=7cm]{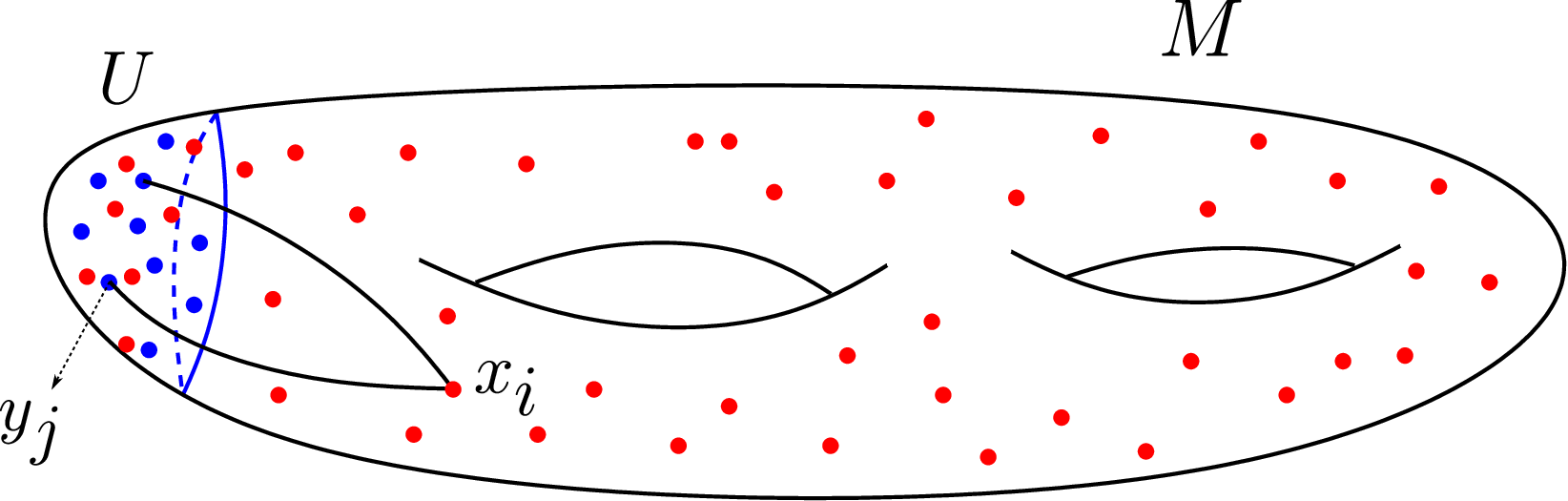} 
\label{pict1aaD}
\end{center}

\caption{\it  
On a closed manifold $(M,g)$, we consider the distances $d(x_i,y_j)$ from 
the blue  points $y_j\in Y$ in
the open subset  $U\subset M$ to  the red points $x_i$ filling the manifold $M$. These distances 
with measurement errors  $\e_{i,j}$ define the noisy distance vectors $\vec R_{i}=[\hat R_{i,j}]$,
where 
$\hat R_{i,j}=d(x_i,y_j)+\e_{i,j}$ and $|\e_{i,j}|<\e_1$.
The inverse problem is to construct  an approximation of the manifold $(M,g)$ from these data.
}
\end{figure}


\subsubsection{Formulation of data for the inverse problem: distance vector data}

Let $\epszero>0$ be a small parameter and let $Y$ be a finite $\epszero$-net in $U$,
\begin{equation}\label{set Y}
Y=\{y_j\in U: j=0,1,\dots,J\}.
\end{equation}
Assume  $y_0\in Y$ is such that $d(x_0,y_0)<\epszero$, where $d(x,y)=\dist_{M}(x,y)$ is  the distance 
between points $x,y\in M$.

The {\it distance vector data} consist of a finite set of vectors $\vec R_i\in \R^{J+1}$, $i=1,2,\dots,I,$ given by
\beq
\vec R_{i}=[\hat R_{i,j}]_{j\in \{0,\dots,J\}}\in \R^{J+1}.
\eeq
We assume that the vectors $\vec R_{i}$ satisfy the following conditions for some small parameter $\epsone>0$:
%
%
\medskip

\noindent
(a1) For any $i=1,2,\dots,I$, there exists a point $x\in M$ such that
\ba\label{r jp correspondence with matrix}
\big|\hat R_{i,j}-d(x,y_j) \big|<\epsone, \quad \hbox{for all }\ j=0,1,\dots,J.
\ea

\noindent
(a2) For any $x\in M$, there is $i\in \{1,2,\dots,I\}$   such that 
\ba\label{second r jp correspondence with matrix}
\big|\hat R_{i,j}-d(x,y_j) \big|<\epsone, \quad \hbox{for all }j=0,1,\dots,J.
\ea


\subsubsection{An alterernative formulation of data: approximate interior distance function data.}

For $x\in M$, we consider the distance function $r^U_x:U\to \R$ defined by
\beq
r^U_x(y)=d(x,y),\quad y\in U.
\eeq
Let $\epszero>0$ and the set $Y$ given in \eqref{set Y} be a finite $\epszero$-net in $U$.
Let
\beq\label{def RY}
{\mathcal R}_Y(M):=\{ r_x^U \big|_Y:\ x\in M\}\subset \ell^\infty(Y)=\R^{J+1}
\eeq
be the set of the restrictions of distance functions $r_x^U$ onto the finite subset $Y$.
The motivation of these functions is that they are discretizations of the distance functions
\beq\label{def RU}
{\mathcal R}_U(M):=\{ r_x^U :\ x\in M\}\subset C(U)
\eeq
defined on the open set $U$.

The {\it approximate interior distance function data}
 consist of the finite set $Y$, and a finite set  of functions on $Y$ given by
\begin{equation}
\widehat{{\mathcal R}}_Y := \{
\widehat{r}_{i}:Y\to \R \,|\,
i=1,2,\dots,I\}\subset \R^{J+1}.
\end{equation}
We assume that the family $\widehat{{\mathcal R}}_Y$ satisfies 
\begin {eqnarray}\label{Y hausdorff dist A}
d_{H} \big(\widehat{{\mathcal R}}_Y,\, {\mathcal R}_Y(M) \big)<\epsone
\end {eqnarray}
for some small parameter $\epsone>0$, that is,
 $\widehat{{\mathcal R}}_Y$ is an approximation of the set ${\mathcal R}_Y(M)$. Here $d_H$ stands for the
Hausdorff distance on $\R^{J+1}$, see \cite{Burago}. The following lemma motivates the conditions (a1)-(a2).

\begin{lemma} \label{lem: data equivalence}
Let $U\subset M$  be an open subset and $Y$ be a finite $\epszero$-net in $U$.
Then the  approximate interior distance function data $\{Y, \widehat{{\mathcal R}}_Y \}$ 
satisfy the condition \eqref{Y hausdorff dist A} if and only if 
 the {distance vector data} 
$
\vec R_{i}=[\hat R_{i,j}]
$, $i=1,2,\dots,I,$ defined by  
 $\hat R_{i,j}=\widehat{r}_{i}(y_j)$, 
satisfy the conditions (a1) and (a2).
\end{lemma}

The proof of Lemma \ref{lem: data equivalence} in given in Section \ref{prel}.


 %

\subsection{Main result}


Our  main result is a global result on the determination of a smooth manifold $(M,g)$ from 
the distance vector data consisting of the noisy distances of the points $x_i\in M$  to the points $y_j$ in an open set $U\subset M$ (see Figure 1).

\begin {theorem}\label{main 1 no boundary}
Let $n\in \Z_+$, $n\geq 2$, $\Lambda \geq 1$, $R>\Lambda^{-1}$.
Then there exist ${\hatepsone}>0$, $\Cmainthree,\Cmainfour,\Cmainfive>1$ explicitly depending only on $n, \Lambda$, such that the following holds for all $0<\epsone<\hatepsone$ and  $0<\epszero\leq \epsone$.

Let $(M,g)$ be a closed Riemannian manifold satisfying the bounds \eqref{basic 1} and $U\subset M$ be an open subset containing a ball $B(x_0,R)$.
Let $$Y=\{y_j: j=0,1,\dots,J\}\subset U$$ be a finite $\epszero$-net in $U$,
 $d(x_0,y_0)<\epszero$, and $\epsX=\Cmainthree \epsone^{1/2}$.

Assume that we are given  vectors $\vec R_i\in \R^{J+1}$, $i=1,2,\dots,I,$ such that conditions (a1) and (a2)
are valid with parameter $\epsone$.
Then the following statements hold.

\begin{itemize}
\item[(1)]  We can compute the numbers $\widehat{d}_{i,i'}$, $i,i'\in \{1,2,\dots,I\}$ directly from the given data 
 $\vec R_i\in \R^{J+1}$, $i=1,2,\dots,I,$ such that
 there exists an $\epsX$-net ${X}=\{x_1,\dots,x_I\}$ in $M$ for which 
\ba\label{final whole manifold estimate2}
\Big|\widehat{d}_{i,i'}-d(x_i,x_{i'}) \Big|\leq
\Cmainfour  \epsone^{\frac18},\quad \hbox{for all } i,i'\in \{1,2,\dots,I\}.
\ea


\item[(2)] 
The given data 
 $\vec R_i\in \R^{J+1}$, $i=1,2,\dots,I,$
determine a smooth Riemannian manifold $(\hat M,\hat g)$ that is
diffeomorphic to $M$. Moreover, there is a diffeomorphism
 $F: \hat M\to M$ such that 
 \ba
\frac 1L\leq \frac{d_{M}(F(x),F(x'))}{d_{\hat M}(x,x')}\leq L,\quad \hbox{for }x,x'\in \hat M,
\ea
where  $L=1+\Cmainfive \epsone^{1/12}$.
\end{itemize}
\end{theorem}

We will focus on proving the claim (1) in this paper, as  the claim (2)  essentially follows from the claim (1) and the part 1 of the paper, \cite{FIKLN}.

\begin{remark}\label{remark-main}
In the form of the inverse problem formulated in Figure \ref{pict1aaD}, Theorem \ref{main 1 no boundary} can be formulated in the following way. Suppose we are given an $\varepsilon_0$-net $Y=\{y_j\}$ in $U$ and an $\varepsilon_2$-net $X=\{x_i\}$ in $M$. Then the noisy distance data $\hat R_{i,j}=d(x_i,y_j)+\e_{i,j}$, where $|\e_{i,j}|<\varepsilon_1$, determine a smooth Riemannian manifold $(\hat M,\hat g)$ that is diffeomorphic to $M$. Moreover, there is a bi-Lipschitz diffeomorphism $F: \hat M\to M$ with Lipschitz constant $1+\Cmainfive \epsone^{1/12}$.
\end{remark}

\section{Inverse problem for the heat kernel with partial data and the local reconstruction of the manifold}

\subsection{Inverse problem for the heat kernel with noisy data}

%
Let $G(x,z,t)$ be the  heat kernel of a Riemannian manifold $(M,g)$, i.e., it satisfies
\begin{eqnarray*}
& &(\partial_t-\Delta_g)G(x,z,t)=0,\quad\hbox{for }(x,t)\in M\times \R_+,\\
& &G(x,z,t)|_{t=0}=\delta_z(x),
\end{eqnarray*}
where $\Delta_{g}$ is the Laplace-Beltrami operator on $(M,g)$  that operates in the $x$-variable
and $\delta_z$ is a point source at the point $z\in M$.
Let
\beq \label{kernel-measurements}
\tilde G(x,z,t)=\eta(x,z,t)\,G(x,z,t)
\eeq
be the values of the heat kernel with multiplicative noise $\eta(x,z,t)$ satisfying
\beq \label{eta-error}
 \Big| \log \eta(x,z,t)\Big | \leq \frac{\sigma}{t}, \quad\hbox{ for $0<t<1$,}
 \eeq
 where $\sigma \in (0,1)$ is small. We consider the stability of the following inverse problem on
 Riemannian manifolds with non-negative Ricci curvature.
  \medskip
 
 {\bf Inverse problem for heat kernel with separated sources and observations.} {\it Let $(M,g)$ be a compact Riemannian manifold and
 $U\subset M$ be a non-empty open subset. Assume that we are given the set $(U,g|_U)$
 as a Riemannian manifold  and the 
heat kernel $G(y,z,t)$ at observation points $y\in U$ at all times  $t\in (0,1)$ with the source points $z\in M\setminus \overline U$. Do these data uniquely determine, the topology,
 the differentiable structure and the metric of the manifold $(M,g)$?}
 \medskip

In  the case when
the heat kernel $G(y,z,t)$ is sampled on the set
$M\times M\times \R_{+}$, that is, sources and observations are on the whole manifold $M$, this problem in studied in
 the embedding of a manifold into an Euclidean space
using the heat kernel \cite{Besson,Portegies,WZ}, in diffusion distances in manifold learning \cite{CoifmanLafon,diffusion},
and in  manifold registration in shape theory \cite{Shape1,Shape2}. Also, an inverse problem where $G(y,z,t)$ is assumed to be known in the set
$U\times U\times \R_{+}$ for an open subset $U\subset M$ was studied in
\cite{Helin-etal,Helin-etal2,KrKuLa}, see also related studies \cite{AKKLT,BeKu,BKL,dCristo-Rondi,KKL,SHoop} for
manifolds with boundary.
In the methods used to study these problems it is essential that the set which contains the sources 
intersects the set where the solutions are observed. For partial data problems with separated sources and observations,
the inverse problem for the wave equation on a non-trapping manifold is considered in \cite{LO}.
Inverse problems for elliptic equations, where
the sources  and  the observations are on sets that are small but intersect, have been studied under geometric
convexity assumptions, see e.g. \cite{dsFKSU,KSU,KrU} and in the 2-dimensional case, see e.g. \cite{GT,LTU,LU}.

The following result proves that for manifolds with non-negative Ricci curvature, the Riemannian manifold structure depends in a stable way 
on the heat kernel with separated sources and observations when  the data have multiplicative errors.

\begin{theorem} \label{Thm-kernel}

Let $(M_1,g_1)$ and $(M_2,g_2)$ be two closed Riemannian manifolds of dimension $n$ satisfying the bounds \eqref{basic 1} with parameter $\Lambda$, and let $U_l=B_{g_l}(y_0^l,R) \subset M_l$ be open balls of radius $R\ge \Lambda^{-1}$. Suppose the Ricci curvatures of the manifolds $M_1$  and $M_2$ are non-negative. Then there exist constants $\hat{\sigma}, \Cmainfivadded>0$ explicitly depending only on $n,\Lambda$, such that the following holds for all $0<\sigma < \hat{\sigma}$ and $0<h\leq \sigma^{1/2}$.

 For $l=1,2$, let $\{z_i^l: i=1,2,\dots,I \}$ be an $h$-net in $M_l\setminus \overline U_l$, and $\{y_j^l: j=0,1,\dots,J \}$ be an $h$-net in the ball $U_l$. Suppose $\Big|d_{M_1}(y_j^1,y_{j'}^1)-d_{M_2}(y_j^2,y_{j'}^2) \Big|<h$
for all $j,j'=0,1,\dots,J$, and
the heat kernels of $M_1$ and $M_2$ satisfy
\begin{equation} \label{kernel-closeness}
e^{-\frac{\sigma}{t}} \leq \frac{G_2 (y_j^2,z_i^2,t)}{G_1 (y_j^1,z_i^1,t)}\leq e^{\frac{\sigma}{t}},
\end{equation}
for all $i=1,\dots,I$, $j=0,1,\dots,J$ and $0<t<1$.
Then $M_1$ and $M_2$ are diffeomorphic, and there is a diffeomorphism
 $F: M_1\to M_2$ such that 
 \ba
\frac 1L\leq \frac{d_{M_2}(F(x),F(x'))}{d_{M_1}(x,x')}\leq L,\quad \hbox{for }x,x'\in M_1,
\ea
where  $L=1+ \Cmainfivadded  \sigma^{1/24}$. 
\end{theorem}

As a corollary we obtain the unique solvability of  the inverse problem.

\begin{corollary} \label{cor: Thm-kernel}
Let $(M_1,g_1)$ and $(M_2,g_2)$ be two closed Riemannian manifolds of dimension $n$ with non-negative Ricci curvature. For $l=1,2$, let $U_l=B_{g_l}(y_0^l,R) \subset M_l$ be open balls of radius $R>0$.
If $\Phi:U_1\to U_2$ is an isometry and  $\Psi:M_1\setminus U_1\to M_2\setminus U_2$ is a bijection such that
\begin{equation*}
G_1 (y,z,t)=G_2 (\Phi(y),\Psi(z),t),\quad \hbox{for all }y\in U_1,\ z\in M_1\setminus U_1,\ 0<t<1,
\end{equation*}
then the Riemannian manifolds $(M_1,g_1)$ and $(M_2,g_2)$ are isometric.

\end{corollary}

We emphasize that in Corollary \ref{cor: Thm-kernel}, the map $\Psi:M_1\setminus U_1\to M_2\setminus U_2$ is  only assumed to be a bijection and thus we do not a priori assume that $M_1$ and $M_2$ are homeomorphic.

Theorem \ref{Thm-kernel} is proved by using the Cheeger-Yau asymptotics for heat kernel  \cite{CY} that generalize 
{Varadhan}'s classical formula   \cite{Varadhan}.


\subsection{Reconstruction of local coordinates and metric tensor from partial distance data}
 
In this subsection, we consider how the local coordinates and the distances near a point $x_{i_0}\in M$ can be approximately constructed, using the approximate distance functions $\widehat{{\mathcal R}}_Y\subset \R^{J+1}$
that are in a neighborhood of the function $\hat r_{i_0}$ corresponding to the point $x_{i_0}$.
The studied question is closely related to a manifold learning problem where  the local structure of the data set
needs to be constructed from the distances to the marker points, see subsection \ref {subsec: marker points}.
 \medskip

\begin{definition}\label{def-xp}
Consider the set $\widehat{{\mathcal R}}_Y$ satisfying \eqref{Y hausdorff dist A}.
We say that $x_i\in M$ is a point corresponding to $\hat r_i \in \widehat{{\mathcal R}}_Y$ if \eqref{r jp correspondence} holds.
For each element $\widehat{r}_i$ in $\widehat{{\mathcal R}}_Y$,  
we choose one corresponding point  $x_i \in M$ and denote the obtained set by ${X}$,
\beq\label{set X}
{X}:=\{x_i\in M\ :\ i=1,2,\dots,I\}.
\eeq 
\end{definition}

Note that the set ${X}$ of points are just known to exist, and they are not directly determined by the data $\widehat{{\mathcal R}}_Y$.

\smallskip
Let $U\subset M$ be an open subset and $Y=\{y_j: j=0,1,\dots, J\}$ be a finite $\epszero$-net in $U$.
In the space of real-valued functions on $Y$, we denote the
$\ell^\infty$-neighborhood of a function  $\hat r:Y\to \R$ by
\beq\label{Binfty}
\mathcal B_\infty(\hat r,\rho)=\{f:Y\to \R\ : \ \|f-\hat r\|_{\ell^\infty(Y)}<\rho\}\subset \R^{J+1},
\eeq
where $\rho>0$ is the radius of the neighborhood.
Suppose that we are given numbers $\hat d^{_Y}_{j,j'}$, $j,j'=0,1,\dots,J$ such that
\beq\label{Da-d estimate A}
\Big|\hat d^{_Y}_{j,j'}-d(y_j,y_{j'}) \Big|\leq 2\epsone,\quad \textrm{for all }j,j'=0,1,\dots,J.
\eeq
 
\smallskip
For $x_{{i_0}}\in M$, the map $\exp_{x_{{i_0}}}:T_{x_{{i_0}}} M\to M$ is the Riemannian exponential map at $x_{{i_0}}$. Let
 $\{v_k\}_{k=1}^n$  be unit vectors that form a basis in $T_{x_{{i_0}}}M$,
 and $x_{\ell}\in B(x_{{i_0}},r)$, $r<\hbox{inj}(M)$. 
 Then we say that 
 \beq\label{correct coordinates}
X(x_\ell)=\big(X_k(x_\ell) \big)_{k=1}^n \in \R^{n},\quad X_k(x_\ell)=\bra \exp_{x_{{i_0}}}^{-1}(x_\ell) ,v_k\cet_g,
\eeq
is the coordinate of the point $x_{\ell}$ in the Riemannian normal coordinates  centered at the point  $x_{{i_0}}$, associated to the (possibly non-orthogonal) basis $\{v_k\}_{k=1}^n$. 
Moreover, 
 \beq\label{correct coordinates metric}
g_{jk}(x_{{i_0}})=\bra v_j,v_k\cet_g, \quad j,k=1,\dots,n,
\eeq
are the components of the metric tensor in these Riemannian normal coordinates at the point $x_{{i_0}}$.

\begin{figure}[htbp]
\begin{center}
\includegraphics[width=7cm]{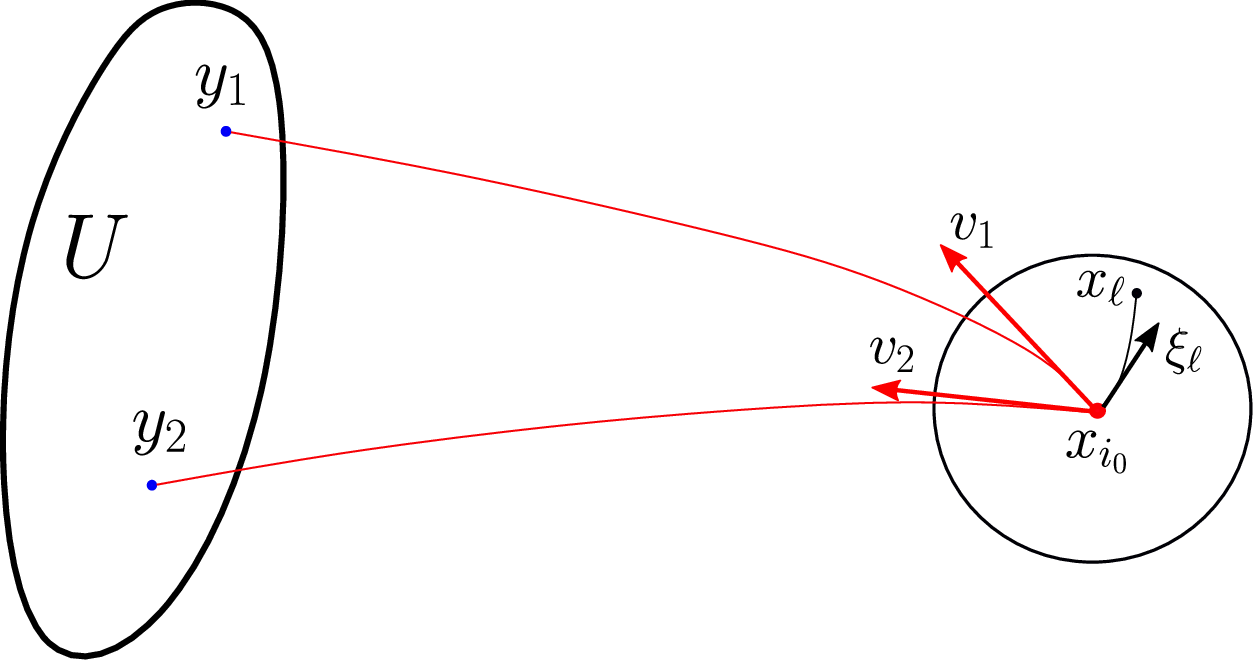} 
\label{figure_coordinate}
\end{center}
\caption{The Riemannian normal coordinates on the manifold $M$  when $\dim(M)=2$.
These coordinates are centered at the point $x_{{i_0}}$ and
 are associated to the basis $\{v_1,v_2\}$. We compute approximately
the coordinates of the points $x_\ell=\exp_{x_{i_0}}(\xi_{\ell})$, that is,  $X(x_\ell)=(X_1(x_\ell),X_2(x_\ell))\in \R^2$
and 
the metric
tensor $g_{jk}$ at the point $x_{{i_0}}$.}
\end{figure}
%
%

%

\begin {theorem}\label{main 1}
Let $n\in \Z_+$, $n\geq 2$, $\Lambda \geq 1$, $R>\Lambda^{-1}$.
Then there exist $\hatepsone>0$, $\crneighbor>0$,
$\cvdetlower>0$ and $\Cmainone>1$, explicitly depending only on
$n$  and $\Lambda$, such that
the following holds for all $0<\epsone<\hatepsone$, $0<\epszero\leq \epsone$.

Let $(M,g)$ be a closed Riemannian manifold satisfying the bounds \eqref{basic 1} with parameter $\Lambda$, and $U\subset M$ be an open subset containing a ball $B(x_0,R)$.
Assume that $Y$ is a finite $\epszero$-net in $U$, $\hat d^{_Y}_{j,j'}$ are numbers satisfying \eqref{Da-d estimate A},
and $\widehat{{\mathcal R}}_Y$ is an $\epsone$-approximation of ${\mathcal R}_Y (M)$ in the sense of \eqref{Y hausdorff dist A}. 
Let $x_{i}\in M$ be the points corresponding to $ \hat r_{i} \in \widehat{{\mathcal R}}_Y$ in the sense of Definition \ref{def-xp}.
Fix any $\hat r_{{i_0}}\in \widehat{{\mathcal R}}_Y$. Then in the tangent space $T_{x_{{i_0}}}M$, 
there 
exist unit vectors $\{v_k\}_{k=1}^n$ 
 satisfying $\det([\bra v_j,v_k\cet_g]_{j,k=1}^n)\geq \cvdetlower$,
 such that the following holds in the Riemannian normal coordinates centered at $x_{{i_0}}$ and associated to the basis $\{v_k\}_{k=1}^n$.

Assume that we are given $Y$, 
$\{\hat d^{_Y}_{j,j'} \}$, $\widehat{{\mathcal R}}_Y\cap
\mathcal B_\infty(\hat r_{{i_0}},\epsone^{1/4})$, and any element $\hat r_{\ell}\in \widehat{{\mathcal R}}_Y$ that satisfies
\begin{equation}\label{crneighbor-Thm1.2}
\hat r_\ell\in \widehat{{\mathcal R}}_Y\cap
\mathcal B_\infty(\hat r_{{i_0}},\rho_0),
\end{equation}
see \eqref{Binfty}. Then we can compute, directly from these data,  numbers $\hat X (x_\ell)\in \R^n$ 
and $\hat g_{jk}\in \R$  such that
%
%
%
 \beq\label{appr. coordinates}
\Big|{\hat X}(x_\ell)-X(x_\ell)\Big| \le \Cmainone \Big( 
d(x_\ell,x_{{i_0}})^{\frac43} + \epsone \Big),
 \eeq
 where $X(x_\ell)$ is the coordinate of the point $x_\ell$ in the Riemannian normal coordinates given in
 \eqref{correct coordinates}, and
 $\hat g_{jk}$ satisfies
\beq\label{metric error} 
\Big|\hat g_{jk}-g_{jk}(x_{{i_0}}) \Big|\leq \Cmainone {\btext \epsone^{\frac18}}, \quad j,k=1,\dots,n,
\eeq
where $g_{jk}(x_{{i_0}})$ are the components of the metric tensor in the Riemannian normal coordinates 
at $x_{{i_0}}$ given in  \eqref{correct coordinates metric}.
\end{theorem}

The  basis vectors $v_k\in T_{x_{{i_0}}}M$ in Theorem \ref{main 1} are 
the directions to nearby points $x_{i(k)}$ such that 
the geodesics $\gamma_{x_{{i_0}},v_{k}}$ can be continued as distance-minimizing geodesics
to points $y_{j(k)}\in U$, see Figure $2$.

\smallskip
Theorem \ref{main 1} has the following corollary.

\begin{corollary}\label{coro 1}
Let $n,\Lambda,R,
{\hatepsone},\crneighbor,\cvdetlower$ be as in Theorem \ref{main 1}. 
Then there  exists $ \Cmaintwo>1$ explicitly depending only on
$n, \Lambda$, such that the following holds for all $0<\epsone<\hatepsone$,
 $0<\epszero\leq \epsone$.

Let $(M,g), U,Y, \widehat{{\mathcal R}}_Y$, $\hat r_{{i_0}}\in \widehat{{\mathcal R}}_Y$ 
and $\hat r_{\ell}\in \widehat{{\mathcal R}}_Y \cap
\mathcal B_\infty(\hat r_{{i_0}},\rho_0)$
be as in Theorem \ref{main 1}.
Then we can compute a number $\widehat{d}_{\ell,{i_0}}$ directly from the given data $Y$, $\{\hat d^{_Y}_{j,j'} \}$, $\hat r_\ell$, and $\widehat{{\mathcal R}}_Y\cap
\mathcal B_\infty(\hat r_{{i_0}},\epsone^{1/4})$ such that
\beq\label{final coro estimate}
\Big|\widehat{d}_{\ell,{i_0}}-d
(x_{\ell},x_{{i_0}})\Big|
\leq 
\Cmaintwo \Big( 
d(x_\ell,x_{{i_0}})^{\frac43} +\epsone^{\frac12} \Big),
\eeq
where $x_{{i_0}},x_{\ell}\in M$ are the points in ${X}$ corresponding to $ \hat r_{{i_0}},\hat r_{\ell}$ in the sense of Definition \ref{def-xp}.
More explicitly, the number $\widehat{d}_{\ell,{i_0}}$ can be computed as
\beq\label{final coro def}
\widehat{d}_{\ell,{i_0}}:=
\bigg(\sum_{j,k=1}^n\hat g^{jk}\hat  X_j(x_{\ell})
\hat  X_k(x_{\ell})\bigg)^{\frac12},\hspace{-20mm}
\eeq
where $\hat X(x_{\ell})=\Big(\hat  X_k(x_{\ell}) \Big)_{k=1}^n$ and the inverse  $(\hat g^{jk})$  of the matrix  $(\hat g_{jk})$ are determined in Theorem \ref{main 1}.
%
%
%
\end{corollary}

Theorem \ref{main 1}  and Corollary \ref{coro 1} show that in a neighborhood of $x_i \in X$,
we can approximately
find the local coordinates of nearby points
 $x_{\ell} \in X$ and the distances $d(x_{\ell},x_i)$. 
 We use this result to reconstruct distances in a finite net of $M$ to prove Theorem \ref{main 1 no boundary}(1).
%
%
Similar results for manifolds with boundary have been studied in \cite{KatsudaKL}, assuming 
bounded derivative of the curvature tensor and using techniques that do not give explicit dependency on the geometric bounds. 
We emphasize that in Theorem \ref{main 1}  and Corollary \ref{coro 1} all the estimates can be made completely explicit in
terms of $\Lambda$ and $n$.

The method developed in the present paper replies on Toponogov's theorem. 
Suppose we do measurements in a ball $B$ having the center $y_0$ and radius $R$. Any point $x \in M\setminus B$ can be connected to $y_0$ by some distance-minimizing geodesic $[y_0 x]$. If $y_1\in [y_0 x]$ has distance $R/2<d(y_1,y_0)<R$, then the geodesic segment $[y_1 x]$ from $y_1$ to $x$ is minimizing and has no cut points. Furthermore, one can perturb $[y_1 x]$ slightly so that the perturbed geodesic has no cut points either, see Section \ref{sec-first-variation} for a quantitative formulation.
Technically, when
we construct approximate values of the metric tensor, we handle this construction related to long geodesics by using only one-sided comparison estimates for large triangles, such as Toponogov's theorem, or by using comparison estimates for small triangles near the point $x$.

A situation considered in Theorem \ref{main 1} is encountered in imaging applications,
where the wave speed (the metric tensor $g$) needs to be reconstructed near a point $x_{{i}}$
using the travel times  of waves from nearby points $x_\ell$ to the points $y_j\in Y$. For example, consider the case when a measurement device located at the point $y_j$  sends a wave at time $t=0$ which reflects from a small scatterer (e.g. a detail in the material) at the point $x_\ell$. If the reflected wave is observed
at the point $y_j$  at the time $t_{\ell,j}$, then the distance $d(x_\ell,y_j)$  is equal to $t_{\ell,j}/2.$


\subsection{Other applications and the stability of inverse problems under geometric a priori bounds}

The inverse problem of determining a Riemannian manifold $(M,g)$ from the distance functions $r_x:U\to \R$, 
$r_x(y)=d(x,y)$, defined on
an open subset $U\subset M$, is encountered in imaging problems that arise in geosciences,
medical imaging and non-destructive testing. In these applications, the speed of waves
defines a Riemannian metric on $M$ so that the travel time of the waves from a point $x$ to $y$ is equal to the Riemannian
 distance $d(x,y)$. For instance, in the seismic imaging of the Earth, the
inverse problem of finding the Riemannian metric in normal coordinates corresponds to
finding the physical material parameters of the Earth in the travel time coordinates.

 
Geometric inverse problems have been studied on closed manifolds with data measured on
an open subset of the manifold, see \cite{Helin-etal,Helin-etal2,KrKuLa},  as it is geometrically simpler to formulate the problems on a closed manifold than on a manifold with boundary. On the methods used to solve inverse problems for Riemannian manifold with boundary or a metric on it, see e.g.\ \cite{AKKLT,BeKu,KaKu,KV1,KV2, Ku3,Ku2, LTU,LU,UP,SHoop}.
 In many cases, it is also possible to reduce an inverse problem for a manifold with boundary to an inverse problem for a closed manifold, by means of extending the manifold with boundary to a closed manifold
 and extending measured boundary data to data on an open set.
%
%
%
%

%

\subsubsection{Inverse problems for linear equations}
Let us review a classical inverse problem that is related to the inverse problem studied in this paper.
\medskip


%
%
%
%
%
%
%
%

{\it 1. Inverse interior spectral problem}:
Let $(M,g)$ be a (unknown)
 Riemannian  manifold and $U\subset M$ be an open set.
Assume that we are given
the following data,
\begin {eqnarray}\label{data 2}
\big\{U,\, (\lambda_j)_{j=1}^\infty,\, (\phi_j|_{U})_{j=1}^\infty \big\}.
\end {eqnarray}
Here $\lambda_j$ are the  eigenvalues of the
Laplace-Beltrami operator $\Delta_g$ on $M$ and  $\phi_j$ are the corresponding
orthonormal eigenfunctions. Do data (\ref{data 2}) determine uniquely
  (up to an isometry)
the  Riemannian
manifold  $(M,g)$?

\medskip

The methods used to solve the inverse problem {\it 1}\  consist of
two steps. First,  the given data is used to construct the local distance function
representation ${\mathcal R}_U(M)$ (for a detailed exposition, see e.g.\ \cite{Helin-etal,KKL}). Second,
the manifold $(M,g)$ is reconstructed from
  ${\mathcal R}_U(M)$.

Analogous inverse problems  for the wave
 equation and also for the Maxwell and Dirac systems 
are studied in \cite{AKKLT,Be,KKL,Dirac,KLS}. The problem is closely related to inverse spectral problems 
where only eigenvalues are known, see \cite{Zelditch} 

\subsubsection{Inverse problems for non-linear equations}
The reconstruction of a manifold from partial distance measurements arises also in 
the study of the inverse problems for non-linear partial differential equations.
\medskip

{\it 2. Inverse problem for a non-linear wave equation}:
Let $(M,g)$ be a (unknown)
 Riemannian  manifold and $U\subset M$ be an open set.
Assume that we are given the source-to-solution map $L_U:C^\infty_0(U\times \R_+)\to C^\infty(U\times \R_+)$, $L_U(f)=v|_{U \times \R_+}$,
where $f \in C^\infty_0(U\times \R_+)$ is a source and $v$ is the solution of 
the following non-linear equation
\begin{eqnarray*}
&& (\p_t^2-\Delta_g)v(x,t)+a(x,t)v(x,t)^2=f(x,t),\quad\hbox{in }M\times \R_+, \\
&& v|_{t=0}=0,\quad \p_t v|_{t=0}=0,
\end{eqnarray*}
where $a(x,t)>0$.
Do the set $U$ and the map $L_U$ determine uniquely (up to an isometry)
the Riemannian manifold $(M,g)$ and the coefficient $a(x,t)$?

\medskip

In the study of this problem, the non-linear interaction of linearized waves produced by suitable sources in $U\times \R_+$ can be used to produce  ``artificial'' microlocal point sources at the points $y\in M$, including the unknown region $M\setminus U$ where the original source $f$ vanishes, see  \cite{HUZ21,KLOU2014,KLU18,Lassas, LUW18,UhWa18,WZ2019}. The wave fronts that are produced by
these point sources and are observed in $U$ determine the distances $d(x,y)$ for the points $x\in M$ and $y\in U$.
Thus the inverse problem {\it 2} for the non-linear wave equation is reduced to 
the reconstruction of a manifold from partial distance measurements.

\subsubsection{Manifold learning}\label{subsec: marker points}
In machine learning, an (invariant) manifold learning problem can be formulated as follows.
 \medskip
 
{\it 3. A manifold learning problem}:
Let $(M,g)$ be a
 Riemannian  manifold and $x_i\in M$, $i=1,2,\dots,I$, be an $\e$-dense set of sample points. We consider
 a small subset of these points,  $x_1,\dots,x_J$, $J<I$, as the marker points. Assume that we are given the distances
 $d(x_i,x_j)$, $i=1,2,\dots,I$, $j=1,2,\dots,J$, between the sample points and the marker points.
 Can we obtain an approximation of the manifold $(M,g)$ from these data?
 \medskip
 
 In the case when the marker points $x_1,\dots,x_J$ belong in an open set
 $U\subset M$ and form a $\delta$-dense set in $U$, this problem reduces to the problem studied in this paper. 
Machine learning problems analogous to the problem {\it 3} were studied e.g.\ in \cite{cheng,CoifmanLafon,diffusion,FILN,FMN,RS,TSL,Zha}.
%
%
%
%
%
%
%
%
%
%
%
%
%
%

\subsubsection{A priori bounds and conditional stability of inverse problems}

In the inverse problems above,
the problem of determining the metric from partial data measured on a subset is generally ill-posed  in the sense of Hadamard:
the map from the partial data
to the metric is not continuous so that small change in the data can lead to huge errors in the 
reconstructed metric.
One way out of this fundamental difficulty 
is to assume a priori bounds on the norms of the higher derivatives  of coefficients. Results under this type of conditions are called 
 conditional
stability results and were known mostly for conformally Euclidean metric tensors (see e.g.
\cite{Al,AlS,StU1,StU2}).
However, for inverse problems for general
metric, this approach bears significant difficulties.
The reason is that the usual $C^m$-norm bounds on coefficients
are not invariant, and thus this type of conditions does not suit
 the invariance of the problems under
 diffeomorphisms. Moreover, if 
the structure of the manifold
is not known a priori, this traditional approach cannot be used.

A natural way to overcome these difficulties
is to impose a priori 
{constraints} in an invariant form and consider
a class of manifolds that satisfy invariant a priori bounds,  for instance on curvature, second fundamental form, injectivity radius, etc. Under such type of conditions,  invariant stability results for various inverse problems  
{
have been proven 
}
in \cite{AKKLT,FIKLN,StU1,StU2}.
In particular, for the inverse  interior spectral problem on manifolds with non-trivial topology, stability 
results  
have been obtained in \cite{AKKLT,BKL}, see also \cite{BILL} for analogous results for manifolds with boundary. In the latter, it was shown that the convergence of the boundary spectral data implies
the convergence of the manifolds
with respect to the Gromov-Hausdorff distance. However, the stability for this problem is of $\log\log$  type.

 In addition to being
contaminated with errors, actual measurements typically provide
only a finite set of data. An example of this is 
the classical Whitney problem on
the extension of a function $f:X\to \R$, defined on $X\subset \R^n$, in an optimal way to a function $F \in C^m(\R^n)$, see \cite{W12}. This problem
has been answered in the works of E.\ Bierstone,
Y.\ Brudnyi, C.\ Fefferman, P.\ Milman, W.\ Pawluski, P.\ Shvartsman and others (see
\cite{Br3,Br4,F1,F2,F3,FK1}). 

\smallskip  

This paper is organized as follows. We review a few corollaries of Toponogov's theorem in Section \ref{prel}. In Section \ref{sec-first-variation}, we derive the first variation type of estimates for almost minimizing paths as our main tool. Sections \ref{sec-direction} and \ref{sec-frame} are technical preparations for the reconstruction of local coordinates. Section \ref{sec-mainproof} is devoted to proving the local result on the reconstruction of local coordinates and components of the metric tensor, that is, Theorem \ref{main 1}. We prove the global results, Theorem \ref{main 1 no boundary} and Theorem \ref{Thm-kernel}, in Sections \ref{sec-application} and \ref{Sec. heat kernel}.

\smallskip
\noindent{\bf Acknowledgement.} 
We thank the anonymous referee for thoroughly reading our paper and many valuable comments.
C. Fefferman was partially supported by AFOSR, grant   DMS-1265524, and NSF, grant FA9550-12-1-0425.
S. Ivanov was partially supported by RFBR, grants 14-01-00062 and 17-01-00128-A.
M. Lassas and J. Lu were partially supported by AF, grants 284715 and 312110.
H. Narayanan was partially supported by NSF grant DMS-1620102 and a Ramanujan Fellowship.

\section{Preliminary constructions}\label{prel}

\subsection{Notations}

Let $(M,g)$ be a closed (that is, compact without boundary) connected Riemannian manifold of dimension $n\geq 2$ satisfying the bounds \eqref{basic 1} with parameter $\Lambda$. Let $U\subset M$ be an open subset containing a ball $B(x_0,R)$ of radius $R>\Lambda^{-1}$ centered at $x_0\in U$.
By e.g. \cite[Thm IX.6.1]{Cha}, if
$r<\min\{\textrm{inj}(M)/2,\pi/(2K^{1/2})\}$ where $K$ is the upper bound for the absolute value of sectional curvatures,
the (open) metric balls $B(x,r)$ of $(M,g)$ having radius $r$  and center $x$ are convex.
Thus by making the ball $B(x_0,R)\subset U$ smaller, we can assume that the ball $B(x_0,R)$ is geodesically convex. We denote 
$$\Sigma_r:=\partial B(x_0,r).$$

\smallskip
Pairs $(x,v),\, (y,u)$, etc.\  stand for  points in the tangent bundle $TM$ with $v,u,$ etc.\ 
being tangent vectors.  We identify the vector space $T_v(T_xM)$ with $T_xM$ and
$T_{(x,v)}(TM)$ with $T_xM \times T_xM$, denoting by $(u,w) \in T_xM
\times T_xM$ a tangent vector in $T_{(x,v)}(TM)$. We denote by
$\gamma_{x,{  v}}(t) = \hbox{exp}_x(t{  v})$ the geodesic
emanating from $x$ in the direction ${  v}\in S_xM=\{v\in T_xM :\ |v|_g=1\}$. Geodesics as well
as other rectifiable curves are parametrized by the arclength.


\smallskip
Let $Y=\{y_j\}_{j=0}^J$ be a finite $\epszero$-net of $U$.
Suppose we are given finite number of data 
\begin{equation}
\widehat{{\mathcal R}}_Y := \{
\widehat{r}_{i}:Y\to \R \,|\,
i=1,\dots,I\},
\end{equation}
such that $\widehat{{\mathcal R}}_Y$ is an $\epsone$-approximation of ${\mathcal R}_Y(M)$ in the sense of \eqref{Y hausdorff dist A}.
The given data essentially consists of $I(J+1)$ numbers. To shorten notations, we sometimes denote
\begin{equation}
\widehat{r}_{ij}:=\widehat{r}_i (y_j),\quad y_j\in Y.
\end{equation}


Next, using the above notations, we prove Lemma \ref{lem: data equivalence}

\begin{proof} (of Lemma \ref{lem: data equivalence})
By the definition of the Hausdorff distance on $\R^{J+1}$, the condition \eqref{Y hausdorff dist A} holds if and only if the following two conditions are satisfied. 

\noindent (i) For any $\widehat{r}_i\in \widehat{{\mathcal R}}_Y$, there exists a point $x_i\in M$ such that
\beq\label{r jp correspondence}
\big|\hat r_{i}(y_j)-r^U_{x_i}(y_j) \big|<\epsone, \quad \hbox{for all }j=0,1,\dots,J.
\eeq
\noindent (ii) For any $x\in M$, there is $i\in \{1,2,\dots,I\}$ such that
$$\big|\hat r_{i}(y_j)-r^U_{x}(y_j) \big|<\epsone, \quad \hbox{for all }j=0,1,\dots,J.$$
Since we define $\hat R_{i,j}=\widehat{r}_{i}(y_j)$, the conditions (i),(ii) are the same as
the conditions   (a1),(a2).
\end{proof}

In the proof of Theorem \ref{main 1 no boundary}, we need to approximately determine 
the distances of the points in $Y$.
To do that, we define approximate distances for points in $Y$ as follows:
\beq \label{def-DY}
D^{a}_Y(y_j,y_k)=\inf_{i\in I} (\hat r_{ij}+\hat r_{ik}),\quad y_j,y_k\in Y.
\eeq
Then by the triangular inequality and (\ref{r jp correspondence}), we see that
\beq
D^{a}_Y(y_j,y_k)\geq d(y_j,y_k)-2\epsone.
\eeq
Let $x$ be a point on the shortest geodesic from $y_j$ to $y_k$. By
(\ref{Y hausdorff dist A}), there is $i$ such that $\|\hat r_{i}-r^U_x|_Y\|_{\ell^\infty(Y)}<\epsone$.
Then 
\ba
D^{a}_Y(y_j,y_k)\leq \hat r_{ij}+\hat r_{ik}
\leq d(y_j,x)+d(x,y_k)+2\epsone.\hspace{-1cm}
\ea
Thus we see that
\beq
D^{a}_Y(y_j,y_k)\leq d(y_j,y_k)+2\epsone.
\eeq
The above yields that
\beq\label{Da-d estimate}
|D^{a}_Y(y_j,y_k)-d(y_j,y_k)|\leq 2\epsone.
\eeq
Note that \eqref{Da-d estimate} is independent of $Y$ being a net of $U$. This shows we can find, up to an error of $2\epsone$, the distances between points in $Y$ using only the given data $\widehat{{\mathcal R}}_Y$. 

By \eqref{Da-d estimate}, the numbers $D^{a}_Y(y_j,y_k)$ satisfy the inequality  \eqref{Da-d estimate A}
that we required for the approximate distances
$\hat d^{_Y}_{j,k}$ in the set $Y$. Thus, instead of 
using the notation $\hat d^{_Y}_{j,k}$ 
in the proof of Theorem \ref{main 1}, we  identify these two notations and denote below 
\beq\label{Da-d estimate identification}
\hat d^{_Y}_{j,k}=D^{a}_Y(y_j,y_k).
\eeq

%
%
%
%
%


\medskip
In this paper, $C_1,C_2,\dots\in [1,\infty)$ and $c_1,c_2,\dots\in (0,1)$ denote uniform constants that explicitly depend only on $n$ and $\Lambda$, unless specified.
We also use a generic uniform constant $C>0$ that denotes a number
that explicitly depends only on $n,\Lambda$, unless specified, but its exact value can be different in each appearance even inside one single formula.

\subsection{Implications of Toponogov's theorem}
\label{subsection-Toponogov}

To begin with, we introduce some notations that we will frequently use.
We denote by $[ab]$ a minimizing geodesic (i.e.\ distance-minimizing
curve) connecting the points $a$ and $b$, and let $|ab|=d(a,b)$ denote
the distance between the points $a,b$. 
Let $\beta$ be the angle between the geodesics $[ab]$ and $[bc]$ at point $b$ and
$\theta=\pi-\beta$.

Let $H$ be the rescaled hyperbolic plane with the constant sectional curvature $-\Lambda^2$, and
$\overline d(\overline a,\overline b)$ denotes
the distance between the points $\overline a$ and $\overline b$ in $H$. Denote by
$[\overline a\,\overline  b]$ a minimizing geodesic connecting the points $\overline a$ and $\overline b$.
For our considerations, we usually take $\overline a,\overline b$ in the following way.
Let $\overline a,\overline b,$ and $\overline c$ be points
of $H$
such that $\overline{d}(\overline a,\overline b )=d(a,b)$
, $\overline{d}(\overline b,\overline c)=d(b,c)$  and 
 the angle between the geodesics $[\overline a\,\overline  b]$ and $[\overline b\,\overline  c]$ at 
 $\overline b$
 is  $\beta$. 
Then by Toponogov's
 theorem, the above triangle $abc$ on $M$ and the corresponding
 triangle $\overline a\overline b\overline c$ on $H$ satisfy  $d(a,c)\leq \overline d(\overline a,\overline c)$.

%
%
%
%
%

\smallskip
Now we present a corollary of Toponogov's theorem (e.g. \cite[Thm. 79]{Pe}).
The analogous results to the first variation inequality (\ref{Si far points new}) considered below are
well-known in Alexandrov geometry. Similar types of formulae are used in 
 Section 4.5 of \cite{Burago}
or Section 4 of \cite{Shiohama} or Section 3.6 of \cite{Plaut}. However, 
we present the results in the form needed later and give the proof for the convenience of the reader.


\begin{lemma}\label{lem: Topog estimates general}
There exist uniform constants $\CfirstToponogov,\CsecondToponogov>0$ such that the following holds.

Let $M$ be a closed Riemannian manifold with sectional curvature bounded below by ${\rm Sec}_M \geq -\Lambda^2$.
Let $a,b,c\in M$ and $\beta$ be the angle of the distance-minimizing
geodesics $[ab]$ and $[bc]$ at $b$. 
\medskip

\noindent
(i) Then
\beq
\label{Si far points new}
  |ac| \leq |ab|-|bc|\cos\beta + \CfirstToponogov |bc|^2/\min\{\Lambda^{-1},|ab|\}.
\eeq

\smallskip

\noindent
(ii) In addition to the assumptions above, assume that $|ab|= |bc|$, $|ab|\leq \Lambda$. 
Then
 \beq\label{eq. Topog.}
|ac|\leq 2|ab|(1- \CsecondToponogov \theta^2),\quad\hbox{where 
$\theta=\pi-\beta$.} 
\eeq

\end{lemma}

\noindent
{\bf Proof.} 
(i) See Lemma \ref{lem: Topog estimates general appen}.

(ii) Denote $
|ab|=|bc|=A\leq \Lambda$
so that also $\overline d(\overline a,\overline b)=d(\overline a,\overline b)=A$. 
Let $B=d(a,c)$ and $\overline B=d(\overline a,\overline c)\leq 2\Lambda$. Then by Toponogov's
 theorem $B\leq \overline B$. 
  Moreover, using  {the law of cosines} (\ref{the law of cosines}),
 we can estimate $B$ as follows: We have
 \ba
 \cosh(\Lambda^{} B )
&\leq &
\cosh(\Lambda^{}\overline B )=
\cosh^2(\Lambda^{} A)-\sinh^2(\Lambda^{}A)\cos(\beta)\\
&\leq &1+\sinh^2(\Lambda^{} A) (1-\cos(\beta))
%
\ea
and as $\cosh(2t)=1+2\sinh^2t$, or
\ba
 \cosh(\Lambda^{}  B)=1+2\sinh^2(\frac 12\Lambda^{} B),
\ea
we have,
\beq\label{eq. sinhs pre}
2\sinh^2(\frac 12\Lambda^{}B)&\leq &\sinh^2(\Lambda^{}A) (1-\cos(\beta)).
\eeq
Using the fact that there exists
a uniform constant $C>1$ so that 
\ba
\frac{ u}{w}-\frac{\sinh(u)}{\sinh (w)} \leq C \left (1-\frac{ u}{w}\right),\quad\hbox{for all }u,w\in (0,\Lambda^2],\ u\leq w,
\ea
we see  using (\ref{eq. sinhs pre}) that 
\beq\label{eq. Topog. pre}
\frac { B}{2 A} &\leq&  \frac 1{C+1}(C+\sqrt {\frac {(1-\cos(\beta))}2}).
\eeq

Let $\theta=\pi-\beta$. Using (\ref{eq. Topog. pre}), we see
that there exists   a uniform constant $\CsecondToponogov>0$ such that
\beq\label{eq. Topog. pre 2}
\frac { B}{2 A} &\leq& 1-\frac{1-\sin(\beta/2)}{C+1}=1-\frac{1-\cos(\theta/2)}{C+1} \leq  1- \CsecondToponogov \theta^2.
\eeq
This proves (ii).
\hfill $\square$\medskip

%


The following lemma is a variation of Lemma \ref{lem: Topog estimates general} when $|ab|$ is small.

\begin{lemma}\label{lem:Serg-02}
There exists
a uniform constant $\CthirdToponogov>0$ such that the following holds.

Let $M$ be a closed Riemannian manifold 
with sectional curvature bounded below by 
${\rm Sec}_M\geq -\Lambda^2$.
Let $a,b,c\in M$ be such that 
$
|ab|\geq |bc| .
$
Let $\beta$ be the angle at $b$ between (any pair of) shortest paths
$[ab]$ and $[bc]$. If
$|ab|\leq 1$ and
$|bc|\le \frac 12  |ab|$, then
\begin{equation}
\label{e:1stvarineq}
 |ac| \le |ab|-|bc|\cos\beta + \CthirdToponogov \beta {\color{black}\frac {|bc|^2}{|ab|} }.
\end{equation}
\end{lemma}

\noindent
{\bf Proof.}   Let us first prove the claim in the case when 
 $ |ab|=1$, $|bc|\leq 1/2$.

As above, let $H$ be a rescaled hyperbolic plane of curvature $-\Lambda^2$
and $\bar a,\bar b,\bar c\in H$ be such that $|\bar a\bar b|=|ab|$,
$|\bar b\bar c|=|bc|$ and $\angle \bar a\bar b\bar c=\beta$.
Here by $|xy|$, we denote the distance between points $x$ and $y$ in whatever space they belong.
Then by Toponogov's theorem, $|ac|\le |\bar a\bar c|$.

It remains to prove \eqref{e:1stvarineq} for $\bar a,\bar b,\bar c\in H$
in place of $a,b,c\in M$ and under the assumption that $|\bar a\bar b|=1$.
Let $e_1,e_2$ be an orthonormal basic
of $T_{\bar b}H$ such that $e_1$ is tangent to the geodesic segment
$[\bar b\bar a]$, i.e., $\bar a=\exp_{\bar b} (e_1)$. Let  $r_0=\frac 12$.
For every $r\in[-r_0,r_0]$ and $\eta\in[-\pi,\pi]$ define
a point $\xi(r,\eta)\in H$ by
$$
 \xi(r,\eta)=\exp_{\bar b} (r\cos(\eta) e_1+r\sin(\eta) e_2) .
$$
Note that $\bar c=\xi(|bc|,\beta)$.
Define
$$
 f(r,\eta) = \overline{d}(\bar a, \xi(r,\eta))-1+r\cos(\eta) .
$$
Clearly $\xi:[-r_0,r_0]\times [-\pi,\pi]\to H$ is a smooth map
and its image does not cover $\bar a$.
Therefore $f$ is a smooth function.
(It can be written explicitly using the cosine law
of the hyperbolic plane.)
Observe that $f(r,0)=0$ for all $r$,
hence
$$
\left|\frac {\pd^2}{\pd r^2}  f(r,\eta)\right| \le 2\CthirdToponogov |\eta|,\quad r\in [-r_0,r_0], 
$$
where
$$
 \CthirdToponogov = \frac12 \max_{r,\eta} \left|\frac{\pd^3}{\pd\eta\pd r^2} f(r,\eta)\right| .
$$
By the first variation formula we have
$$
 \frac{\pd}{\pd r} \overline{d}(\bar a,\xi(r,\eta))\bigg |_{r=0} = -\cos\eta ,
$$
and therefore
$$
\frac{\pd}{\pd r} f(r,\eta)\bigg |_{r=0} = 0 .
$$
This and the above estimate on $\pd^2f/\pd r^2$ imply that
$$
| f(r,\eta)| \le \CthirdToponogov|\eta| r^2
$$
for all $\eta\in[-\pi,\pi]$ and $r\in[-r_0,r_0]$.
Substituting $\eta=\beta$ and $r=|\bar b\bar c|$ yields
that
$$
 |\bar a\bar c| - 1 + |\bar b\bar c|\cos\beta \le \CthirdToponogov \beta |\bar b\bar c|^2
$$
or, equivalently,
$$
 |\bar a\bar c| \le |\bar a\bar b| - |\bar b\bar c|\cos\beta + \CthirdToponogov \beta |\bar b\bar c|^2 .
$$
As explained above, the the claim follows from this inequality, Toponogov's theorem
and the triangle inequality in~$M$. 

Thus we have proven the claim in the case when 
$ |ab|=1$. For the case when $|ab|<1$ we can
scale the metric $g$ with the constant  factor $|ab|^{-2}$.
Observe that after this scaling the curvature is still bounded from
below by $-\Lambda^2$. As
the inequality (\ref{e:1stvarineq}) is invariant  under metric scaling by a constant factor, we obtain the inequality (\ref{e:1stvarineq})
also on the case when $|ab|<1$.
\hfill $\square$




\section{First variation for almost minimizing paths}
\label{sec-first-variation}

In this section, let $M$ be a closed Riemannian manifold satisfying the bounds \eqref{basic 1}.
We consider the first variation type of estimates for geodesics, with the help of corollaries of Toponogov's theorem in Section \ref{subsection-Toponogov}.
To explain the idea of constructions we use, we first warm up by proving an improved version of the first variation formula
for geodesics that can be continued as a minimizing geodesic.

\begin{figure}[htbp]
\psfrag{1}{$p$}
\psfrag{2}{$q$}
\psfrag{3}{$x$}
\psfrag{4}{\hspace{-3mm}$z$}
\psfrag{5}{$\alpha$}
\psfrag{6}{$\beta$}
\psfrag{7}{}

\psfrag{9}{$q^\prime$}
\hspace{-8mm}\includegraphics[width=75mm]{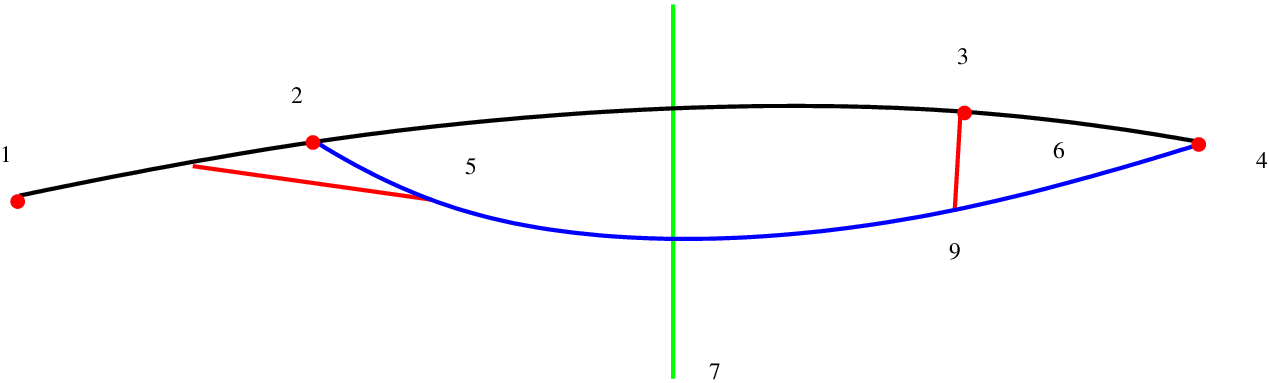} \label{pict11gg}
\psfrag{1}{$p$}
\psfrag{2}{$q$}
\psfrag{3}{$x$}
\psfrag{4}{\hspace{-3mm}$z$}
\psfrag{5}{$\alpha$}
\psfrag{6}{$\beta$}
\psfrag{7}{}
\psfrag{8}{$y$}
\psfrag{9}{$\dot\gamma_1(0)$}
\psfrag{0}{$\theta$}
\psfrag{A}{$\eta$}
\psfrag{B}{$z_2^\prime$}
\psfrag{C}{$z_1$}
\psfrag{D}{$z_2$}
\psfrag{F}{$q^\prime$}
\includegraphics[width=72mm]{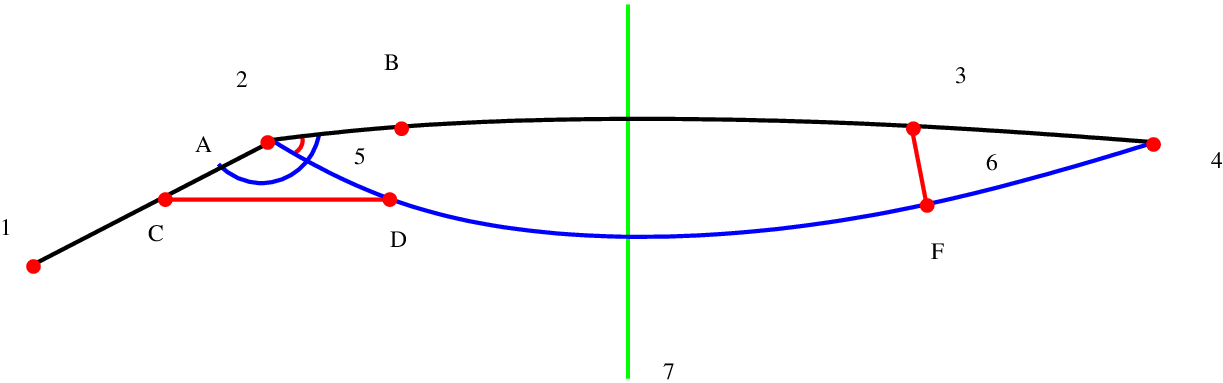} \hspace{-1cm}
%
%
\caption{\it {\bf Left:} Setting of Lemma \ref{lem: First variation lemma A}. The green vertical line corresponds to the boundary of the
set $U$ when the lemma is applied.\ {\bf Right:}  Setting of Lemma \ref{lem: First variation lemma delta}.}
\end{figure}

\begin{lemma}\label{lem: First variation lemma A} 
Let $M$ be a closed Riemannian manifold satisfying the bounds \eqref{basic 1}.
Let $\gamma_{q,v}([0,\ell])$, $\ell=\tilded(q,x)>\Lambda^{-1}$ be a distance-minimizing geodesic, parametrized by arc length, connecting $q,x\in M$. Assume that $p=\gamma_{q,v}(-\tau)$ and 
$\gamma_{q,v}([-\tau,\ell])$ is a distance-minimizing geodesic connecting $p$ and $x$,
where $\tau>\Lambda^{-1}/2$.
Let  $z=\gamma_{q,v}(\ell+r)$ be the point on the continuation of the geodesic $\gamma_{q,v}$,
and $\a$ be the angle of $\gamma_{q,v}$ and $[qz]$ at $q$. 
Then there are uniform constants $\ccalphafirst,\Calphafirst$ such 
that for all $0<r<\ccalphafirst$, we have
\beq\label{eq: alpha r} 
\alpha\leq \Calphafirst  r.
\eeq
\end{lemma}

\smallskip
\noindent {\bf Remark.} 
Lemma \ref{lem: First variation lemma A} can also be formulated in a scaling-invariant form. For example, assume that $\tau<\ell/2$ and $r<\ell/4$. Then we have
\beq
\alpha \leq C_9 \frac{r}{\tau}.
\eeq
This is also the case for the other two lemmas in this section. We note that an estimate of $\alpha\leq Cr^{1/2}$ can be obtained if one only assumes the lower sectional curvature bound.

\medskip
\noindent
{\bf Proof.} 
Without loss of generality, let us assume $\tau>1/2,\, \ell>1$, since this only changes the constants by a factor of $\Lambda$ due to scaling.
Let $\beta$ be the angle of $\gamma_{q,v}$ and $[zq]$ at $z$
and denote by $\a$ the angle of $[qx]$ and $[qz]$ at $q$, that is, $\a=\angle xqz$. 
Let $q^\prime$ be a point on $[zq]$ such that $|zq^\prime|=2r$.
Lemma \ref{lem:Serg-02}, applied for the triangle
$xzq^\prime$, implies that for $r<\min\{1/4,\textrm{inj}(M)/4\}$,
\ba
|xq^\prime|\leq |zq^\prime|-r\cos \beta+\CthirdToponogov \beta  {\color{black} r}.
\ea

Let $z_1\in [qp]$ and $z_2\in [qz]$ be points such
that $|z_1q|=1/4 $ and $|z_2q|=1/4$. 
Then using a shortcut argument near the point $q$,
see Fig.\ 3(Left) where the shortcut is the red segment,
we can compare the distances of $|z_1q|+|qz_2|$ and
$|z_1z_2|$, and   using (\ref{eq. Topog.}), we see that
\ba
|px|&\leq&|pz_1|+|z_1z_2|+|z_2q^\prime|+|q^\prime x|
\\  &\leq &
 |pz_1|+(|z_1q|+|qz_2|-\frac12 \CsecondToponogov \a^2)+|z_2q^\prime|+
(|q^\prime z|-r\cos \beta+\CthirdToponogov  {\color{black} r}\beta).
\ea
Here, $|pz_1|+|z_1q|=|pq|$ and $|qz_2|+|z_2q^\prime|+|q^\prime z|=|qz|$ due to $\tau>1/2$, $\ell>1$, $r<1/4$.
Thus
\beq\label{very nice estimate A}
|px|\leq |pq|+|qz|-\frac12 \CsecondToponogov \a^2-r\cos \beta+\CthirdToponogov  {\color{black} r}\beta.
\eeq
As $|px|=|pq|+|qx|$, this yields
\beq\label{very nice estimate A 2}
|qx|\leq |qz|-\frac12 \CsecondToponogov \a^2-r\cos \beta+\CthirdToponogov  {\color{black} r}\beta.
\eeq

Using (\ref{very nice estimate A 2}) and the fact that
$|qz|\leq |qx|+r$,
 we see that 
\begin{eqnarray}
|qx|&\leq& |qx| +(1-\cos \beta)r-\frac12 \CsecondToponogov \a^2 +\CthirdToponogov  {\color{black} r}\beta \nonumber
\\
&\leq& |qx| +\frac 12 \beta^2 r -\frac12 \CsecondToponogov \a^2+\CthirdToponogov  {\color{black} r}\beta. \label{alpha-initial}
\end{eqnarray}
Note that \eqref{alpha-initial} yields only $\alpha^2\leq Cr$. We need the following improvement to obtain the desired estimate.

\medskip
Pick the points $x_1\in [qx]$, $x_2\in [qz]$ such that $|zx_1|=|zx_2|=\Lambda^{-1}/2\leq {\rm inj}(M)/2$. 
Lemma \ref{holderbound} yields that
$$d(x_1,x_2) \leq \CexpLip(\alpha+r).$$
On the other hand, by the Rauch comparison theorem for $\textrm{Sec}_M\leq \Lambda^2$,
$$\CexpLip(\alpha+r)\geq d(x_1,x_2) \geq C(n,\Lambda)\sin \frac{\beta}{2} \geq \frac{1}{4} C(n,\Lambda)\beta.$$
This shows that for suitable $\Calphabeta>1$, we have
\begin{equation}\label{alphabetar}
\beta\leq \Calphabeta(\alpha+r).
\end{equation}

\smallskip
If $\alpha\leq r$, the claim is proven. Thus we may assume that $r\leq\alpha$. Then \eqref{alphabetar} becomes
\begin{equation}\label{alphabeta}
\beta\leq 2\Calphabeta\alpha.
\end{equation}
Hence by \eqref{alpha-initial}, we see that 
\ba
0&\leq&\frac 12 4\Calphabeta ^2 \a^2 r -\frac 12 \CsecondToponogov \a^2+\CthirdToponogov  {\color{black} r}\beta\\
&\leq&\frac 12 ( 4\Calphabeta ^2 r -\CsecondToponogov )\a^2+2\CthirdToponogov\Calphabeta   {\color{black} r}\alpha.
\ea
Assuming that $r<\ccalphafirst$ where $\ccalphafirst\leq \CsecondToponogov \Calphabeta ^{-2}/8$ we have 
$\CsecondToponogov - 4\Calphabeta^2r>\frac 12 \CsecondToponogov $, and hence
\ba
0&\leq&-\frac 14 \CsecondToponogov  \a^2+2\CthirdToponogov\Calphabeta   {\color{black} r}\a
\ea
or
\beq\label{eq: alpha r squared}
 \a\leq 8\CsecondToponogov ^{-1}\CthirdToponogov\Calphabeta   {\color{black} r}=:\Calphafirst  {\color{black} r}.
\eeq
\hfill $\square$
\medskip

We remark that if there are multiple distance-minimizing geodesics from $q$ to $z$, then \eqref{eq: alpha r} remains valid for each geodesic.

\smallskip
Next we modify the assumptions of Lemma \ref{lem: First variation lemma A} by
replacing the minimizing geodesic by an almost minimizing path.

%

\begin{lemma}\label{lem: First variation lemma delta}
Let $M$ be a closed Riemannian manifold satisfying the bounds \eqref{basic 1}.
Let $\gamma_{q,v}([0,\ell])$, $\ell=\tilded(q,x)>\Lambda^{-1}$ be a distance-minimizing geodesic, parametrized by arc length, connecting $q,x\in M$. 
Assume that there is a curve from $p$ to $x$ that
goes through $q$ that is almost distance-minimizing in the sense that
\beq\label{eq: almost minimizing}
|pq|+|qx| \leq |px|+\delta.
\eeq
Assume $|pq|>\Lambda^{-1}/2$.
Let $z=\gamma_{q,v}(\ell+r)$ be the point on the continuation of the geodesic $\gamma_{q,v}$,
and
 $\a$ be the angle of $\gamma_{q,v}$ and $[qz]$ at $q$. 
 Then there are uniform constants $\ccalphafirst,\Calphasecond,\Cbetadisfirst>1$ 
such that for all $0<r,\delta<\ccalphafirst$, we have 
\beq\label{first claim}
\alpha\leq \Calphasecond  ({\color{black} r^2}+\delta)^{1/2}.
\eeq
Moreover, 
\beq\label{second claim}
|qx|+|xz|\leq |qz|+
\Cbetadisfirst  {\color{black} r (r^2+\delta)^{1/2}}.
\eeq
\end{lemma}

\noindent
{\bf Proof.} 
Without loss of generality, let us assume $|pq|>1/2,\, |qx|>1$, since this only changes the constants by a factor of $\Lambda$ due to scaling.
Let $\beta$ be the angle of $\gamma_{q,v}$ and $[zq]$ at $z$
and let $\eta$ be the angle of the (distance-minimizing) geodesic segments $[qp]$ and $[qx]$ at $q$. Denote
$\a=\angle xqz$.
Let $q^\prime$ be a point on $[zq]$ such that $|zq^\prime|=2r$.
Then Lemma \ref{lem:Serg-02}, applied for the triangle
$xzq^\prime$ implies that for $r<\min\{1/4,\textrm{inj}(M)/4\}$,
\beq\label{apu1}
|xq^\prime|\leq |zq^\prime|-r\cos \beta+\CthirdToponogov \beta  {\color{black} r}.
\eeq

Let $z_1\in [qp]$, $z_2\in [qz]$, $z_2^\prime \in [qx]$ be the points such
that $|z_1q|=1/4$, $|z_2q|=1/4$, $|qz_2^\prime|=1/4$, see Figure $3$(Right).
Then \eqref{eq. Topog.} applied to the triangle $z_1 z z'_2$ gives
\beq
|z_1 z'_2|\leq |z_1q|+|q z'_2|-2|z_1 q| \CsecondToponogov (\pi-\eta)^2.
\eeq
Then we use a shortcut argument as follows.
\begin{eqnarray*}
|px| &\leq& |p z_1|+|z_1 z'_2|+|z'_2 x|  \\
&\leq& |p z_1|+\Big(|z_1q|+|q z'_2|-2|z_1 q| \CsecondToponogov (\pi-\eta)^2 \Big)+|z'_2 x| \\
&=& |pq|+|qx|-\frac{1}{2} \CsecondToponogov (\pi-\eta)^2.
\end{eqnarray*}
As $|px|\geq |pq|+|qx|-\delta $, this yields
\beq \label{estimate-eta}
(\pi-\eta)^2\leq 2 \CsecondToponogov^{-1} \delta.
\eeq

We use a shortcut argument near the point $q$,
see Fig.\ 3(Right) where the shortcut is the red segment $[z_1z_2]$. We compare the distances of $|z_1q|+|qz_2|$ and
$|z_1z_2|$ using (\ref{eq. Topog.}), and see that 
\beq\label{eq: sc 1}
|z_1z_2|\leq |z_1q|+|qz_2|-2|z_2q| \CsecondToponogov (\pi-\omega)^2,
\eeq
where $\omega$ denotes the angle between $[qp]$ and $[qz]$ at $q$.

Inequality (\ref{eq: sc 1}) implies then that
\ba
|px|&\leq&|pz_1|+|z_1z_2|+|z_2q^\prime|+|q^\prime x|
\\  &\leq &
 |pz_1|+\Big(|z_1q|+|qz_2|-\frac 12 \CsecondToponogov (\pi-\omega)^2 \Big)+|z_2q^\prime|+\\
 & &\quad +
(|q^\prime z|-r\cos \beta+\CthirdToponogov  {\color{black} r}\beta).
\ea
Here, $|pz_1|+|z_1q|=|pq|$ and $|qz_2|+|z_2q^\prime|+|q^\prime z|=|qz|$.
Thus
\beq\label{very nice estimate A BB}
|px|\leq |pq|+|qz|-\frac 12\CsecondToponogov (\pi-\omega)^2-r\cos \beta+\CthirdToponogov  {\color{black} r}\beta.
\eeq
As $|px|\geq |pq|+|qx|-\delta $, this yields
\beq\label{very nice estimate A 2 BB}
|qx|-\delta \leq |qz|-\frac 12\CsecondToponogov (\pi-\omega)^2-r\cos \beta+\CthirdToponogov  {\color{black} r}\beta.
\eeq
Using (\ref{very nice estimate A 2 BB}) and the fact that
$|qz|\leq |qx|+r$,
 we see that  
\beq \label{alphabeta-2-initial}
|qx|-\delta &\leq& |qx| +(1-\cos \beta)r-\frac 12\CsecondToponogov (\pi-\omega)^2 +\CthirdToponogov  {\color{black} r}\beta.
\eeq
Note that \eqref{alphabeta-2-initial} already implies $(\pi-\omega)^2\leq C(r+\delta)$, which combining with \eqref{estimate-eta} yields $\a \leq |\pi-\eta|+|\pi-\omega| \leq C(r+\delta)^{1/2}$. We still need the following improvement to obtain the desired estimate.


\medskip
Similar as in Lemma \ref{lem: First variation lemma A}, by using the Rauch comparison theorem in a ball of radius $\Lambda^{-1}/2$ centered at $z$, for sufficiently small $r,\delta$, we have
\begin{equation}\label{alphabetar-2}
\beta\leq \Calphabeta(\alpha+r).
\end{equation}
If $\alpha\leq r$, the claim is proven. Thus we may assume that $r\leq\alpha$. Then \eqref{alphabetar-2} becomes
\begin{equation}\label{alphabeta-2}
\beta\leq 2\Calphabeta\alpha.
\end{equation}
Hence using (\ref{alphabeta-2-initial}) and $\alpha \leq |\pi-\omega|+(2\CsecondToponogov^{-1}\delta)^{1/2}$ by \eqref{estimate-eta}, we see that 
\ba
-\delta&\leq&\frac 12 4\Calphabeta ^2 \a^2 r -\frac 12\CsecondToponogov (\pi-\omega)^2+2\CthirdToponogov \Calphabeta
 {\color{black} r}\alpha\\
&\leq&\frac 12 (8\Calphabeta ^2 r -\CsecondToponogov )(\pi-\omega)^2+2\CthirdToponogov \Calphabeta    {\color{black} r}|\pi-\omega| +\Caux r\delta+\Caux r\delta^{\frac12},
\ea
for some constant $\Caux>0$ depending on $\CsecondToponogov, \CthirdToponogov, \Calphabeta$.
Assuming that $r<\ccalphafirst$ where $\ccalphafirst\leq \CsecondToponogov \Calphabeta ^{-2}/16$, we have 
$\CsecondToponogov - 8\Calphabeta^2r>\frac 12 \CsecondToponogov $. Hence,
\ba
\frac 14 \CsecondToponogov  (\pi-\omega)^2-2 \CthirdToponogov \Calphabeta {\color{black} r}|\pi-\omega| -(1+\Caux r)\delta-\Caux r\delta^{\frac12} \leq0. \ea
Therefore,
\ba
\a &\leq& |\pi-\omega| +2\CsecondToponogov^{-\frac12} \delta^{\frac12} \\
&\leq& \frac {2\CthirdToponogov \Calphabeta   r+\sqrt{( 2\CthirdToponogov \Calphabeta  r)^2+ \CsecondToponogov (1+\Caux r) \delta +\CsecondToponogov \Caux r\delta^{1/2}}} {\CsecondToponogov /2}+2\CsecondToponogov^{-\frac12} \delta^{\frac12}
\\ &\leq &\Calphasecond  ( {\color{black} r^2}+\delta)^{1/2},
\ea
with some uniform constant $ \Calphasecond  >1$.
This proves the first claim (\ref{first claim}) of the lemma.

\medskip
For the second claim, let $q^\prime\in M$ be a point on $[zq]$ such that $|zq^\prime|=2r$.
Lemma \ref{lem:Serg-02}, 
applied for the triangle 
$xzq^{\prime}$ and its angle at $z$ implies that
\beq\label{apu2 A}
|q^{\prime}x|    \leq |q^{\prime}z|-|xz|\cos \beta +\CthirdToponogov  {\color{black} r} \beta.
\eeq
Adding $ |qq^{\prime}|$ on both sides of (\ref{apu2 A}) and using the facts that
$|qx|\leq |qq^{\prime}|+|q^{\prime}x|$ and  $|qz|= |qq^{\prime}|+|q^{\prime}z|$, we obtain
\beq\label{apu2}
|qx|&\leq& |qq^{\prime}|+|q^{\prime}x|\\ \nonumber
&\leq&( |qq^{\prime}|+|q^{\prime}z|)-|xz|\cos \beta + \CthirdToponogov  {\color{black} r}\beta \\
\nonumber
&=&|qz|-|xz|\cos \beta +  \CthirdToponogov  {\color{black} r}\beta.
\eeq
Using the triangle inequality, (\ref{alphabetar-2}), (\ref{first claim}),  and
the facts that $|1-\cos \beta|\leq {\frac 12}\beta^2$ and $|xz|=r$, we see that
\beq\label{apu3}
0 \leq |qx|+|xz|- |qz|&\leq&\CthirdToponogov \beta  {\color{black} r}+ {\frac 12}r\beta^2\\
\nonumber 
&\leq& \Cbetadisfirst {\color{black} r (r^2+\delta)^{1/2}}
\eeq
with some suitable uniform constant $\Cbetadisfirst>1$.
\hfill $\square$ \medskip
%


\begin{proposition}\label{lem: inner product lemma} 
Let $M$ be a closed Riemannian manifold satisfying the bounds \eqref{basic 1}.
Then there is a uniform constant $\hat \delta\in (0,\ccalphafirst^2)$ such that the following holds for all $0<\delta<\hat \delta$.

Let $p,q,x\in M$ such that $|pq|>\Lambda^{-1}/2$, $|qx|>\Lambda^{-1}$ and
\beq\label{eq: almost minimizing-2}
|pq|+|qx| \leq |px|+\delta.
\eeq
Denote by $\theta$ the angle of $\gamma_{q,v}$ and $[xy]$ at $x$.
Let $y\in M$ be such that $d(y,x)<\ccalphafirst^2$. 
Then there is a uniform constant $\betadissecond>1$ such that
\ba
\bigg| |xy|\cos\theta-(|xq|-|yq|)\bigg|\leq {\color{black}\betadissecond (  |xy|\delta^{1/4}+|xy|^{4/3})}.
\ea
\end{proposition}

Note that here $ |xy|\cos\theta$ is the inner product
of the vector $\xi\in T_xM$ for which $\exp_x\xi=y$ and
the unit vector $v\in S_xM$ for which $\gamma_{x,v}$ coincides with the 
geodesic $[xq]$.
\medskip


\begin{figure}[htbp]
\begin{center}
\psfrag{1}{$p$}
\psfrag{2}{$q$}
\psfrag{3}{$x$}
\psfrag{4}{$z$}
\psfrag{5}{$\alpha$}
\psfrag{6}{$\beta$}
\psfrag{7}{}
\psfrag{8}{$y$}
\psfrag{9}{$\dot\gamma_1(0)$}
\psfrag{0}{$\theta$}
\psfrag{A}{$\eta$}
\psfrag{G}{$q^\prime$}
\includegraphics[width=7.5cm]{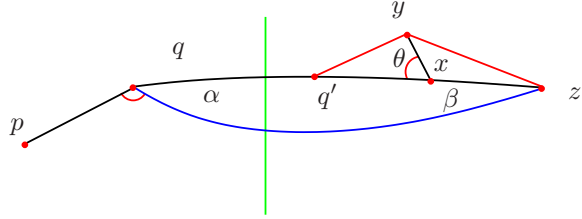} \label{pict13aad}
\end{center}
\caption{\it An auxiliary figure for Proposition \ref{lem: inner product lemma}.}
\end{figure}

\noindent
{\bf Proof.}
Let $v\in T_qM$ be the unit vector such that $[qx]$ is the path $\gamma_{q,v}([0,\ell])$, $\ell=|qx|$. Pick a number $r$ such that $|xy|\leq r<\ccalphafirst$. 
Let $z=\gamma_{q,v}(\ell+r)$ and $q^\prime=\gamma_{q,v}(\ell-{r})$,
so that $|xz|=r$ and $|xq^\prime|=r$. 
 Next we use Toponogov's theorem for the triangles $q^\prime xy$ and $zxy$.

Below in this proof, we use the notation $C$ denoting a generic  constant whose exact value can change even inside one formula.
The value value of each $C$ can be computed as an explicit function of $\Lambda$ and $n$.
Consider the triangle $q^\prime xy$, and we see from
 (\ref{Si far points new}) that for $r<\Lambda^{-1}$,
\ba
|q^\prime y|\leq |xq^\prime |-|xy|\cos\theta +C\frac {|xy|^2}r,
\ea
where $C$ is a uniform constant. 
Then, as $|qq^\prime |+ |xq^\prime |=|qx|$,
\beq\label{E3}
|qy|\leq |qq^\prime |+|q^\prime y|\leq |qx|-|xy|\cos\theta +C\frac {|xy|^2}r.
\eeq
Now considering triangle $xyz$, we see  from (\ref{Si far points new}) that
\beq\label{E1}\\ \nonumber
|yz|\leq |xz|-|xy|\cos(\pi-\theta) +C\frac {|xy|^2}r
= |xz|+|xy|\cos(\theta) +C\frac {|xy|^2}r.
\eeq
 By (\ref{second claim}),
\beq\label{E2}
|qx|+|xz|&\leq& |qz|+C{\color{black} r (r^2+\delta)^{1/2}}
\\ \nonumber 
&\leq& |qy|+|yz|+Cr(r^2+\delta)^{1/2}.
\eeq
Now (\ref{E1}) and (\ref{E2}) yield
\ba
|qx|+|xz|&\leq& |qy|+|yz|+C r(r^2+\delta^2)^{1/2}\\
&\leq& |qy|+(|xz|+|xy|\cos(\theta) +C\frac {|xy|^2}r)
+C {\color{black} r (r^2+\delta)^{1/2}},
\ea
yielding that, after cancelation of  $|xz|$,
\ba
|qx|\leq
 |qy|+|xy|\cos(\theta) +C\frac {|xy|^2}r
+C{\color{black} r (r^2+\delta)^{1/2}},
\ea
or
\beq\label{E5}
|qx|-|qy|\leq
 |xy|\cos(\theta) +C\frac {|xy|^2}r
+C {\color{black} r (r^2+\delta)^{1/2}}.
\eeq
Comparing this with (\ref{E3}), we obtain 
\ba
|xy|\cos\theta -C\frac {|xy|^2}r
\leq 
|qx|-|qy|\leq
 |xy|\cos(\theta) +C\frac {|xy|^2}r
+C{\color{black} r (r^2+\delta)^{1/2}}
\ea
and
hence 
\beq\label{Efinal}
\bigg|
(|qx|-|qy|)-|xy|\cos\theta \bigg| \leq
 C\frac {|xy|^2}r
+C{\color{black} r (r^2+\delta)^{1/2}}=:E.
\eeq

\smallskip
Now we optimize the value of $r$ so that ${|xy|^2}/r\approx {\color{black} r (r^2+\delta)^{1/2}}$
under the requirements that $|xy|\leq r$ and $r<\ccalphafirst$.
First , we  consider the case when $|xy|>\delta^{3/4}$.  Then a good choice is
 ${|xy|^2}= {\color{black} r^3}$, or $r={|xy|^{2/3}}$ so that 
 $|xy| < r$.
Then with some uniform constant $C$, we have
 \ba
E\leq C  \frac {|xy|^2}{|xy|^{2/3}}+C|xy|^{2/3}(2|xy|^{4/3})^{1/2}\leq 4C |xy|^{4/3}.
 \ea
As for the case when $|xy|\leq \delta^{3/4}$, a good choice is
 ${|xy|^2}/r= r\delta^{1/2}$,  or $r=|xy|\delta^{-1/4}$
 so that ${|xy|}= r\delta^{1/4}<r$ and  $r\leq \delta^{1/2}$.
 Then we see  that with some uniform constant $C$,
 \ba
 E\leq C  \frac {|xy|^2}{|xy|\delta^{-1/4} }+C(|xy|\delta^{-1/4}) \delta^{1/2} 
 =2C |xy|\delta^{1/4}.
 \ea
Note that the requirement $|xy|\leq r$ is valid in both cases above. The requirement that $r<\ccalphafirst$ is validated either by the condition $|xy|<\ccalphafirst^2$ in the former case, or by the choice $\hat{\delta}<\ccalphafirst^2$ in the latter case. 
Thus using
the above choices of $r$, we obtain the estimate
\ba
E\leq  \betadissecond (  |xy|\delta^{1/4}+|xy|^{4/3})= \betadissecond |xy|( \delta^{1/4}+|xy|^{1/3})
\ea
with some uniform constant  $\betadissecond>1$.
\hfill $\square$
\medskip

\section{Finding directions of minimizing paths}
\label{sec-direction}



Let $U\subset M$ be an open subset of a closed Riemannian manifold $M$ containing a ball $B(x_0,R)$, and $Y=\{y_j\}_{j=0}^J$ be an $\epszero$-net of $U$.
Suppose we are given an $\epsone$-approximation $\widehat{{\mathcal R}}_Y$ of ${\mathcal R}_Y (M)$ in the sense of \eqref{Y hausdorff dist A}. 
By the definition of Hausdorff distance, for any $\widehat{r}_i \in {\widehat{\mathcal R}}_Y$ 
there exists $x_i \in M$ such that
\begin {eqnarray}
\label{dH-rtox}
|\widehat{r}_i (y_j)-\tilded(x_i,y_j)|< \epsone,\quad\hbox{for all }y_j\in Y.
\end {eqnarray}
We choose for all $\hat r_i \in {\widehat{\mathcal R}}_Y$ 
some point $x_i$ satisfying (\ref{dH-rtox}), and call $x_i$ a corresponding point to the approximate distance function $\hat r_i$ (c.f.\ Definition \ref{def-xp}).
We denote by ${X}=\{x_i\}_{i=1}^I\subset M$ a set
of points chosen so that each point $x_i$ is a corresponding point 
to an approximate distance function  $\hat r_i\in \ell^\infty(Y)$.
We also denote
$\widehat {\mathcal R}_Y=\{\hat r_i\}_{i=1}^I$.


\begin{figure}[htbp]
 
\psfrag{1}{$U$}
\psfrag{2}{$M\setminus U$}
\psfrag{3}{$y_1$}
\psfrag{4}{$y_2$}
\psfrag{5}{$\hspace{-3mm}x_2$}
\psfrag{6}{$x_1$}
\psfrag{7}{$\mu_{\tau}$}
\psfrag{8}{$\dot\gamma_2(0)$}
\psfrag{9}{$\dot\gamma_1(0)$}
\hspace{-1cm}\includegraphics[width=5cm]{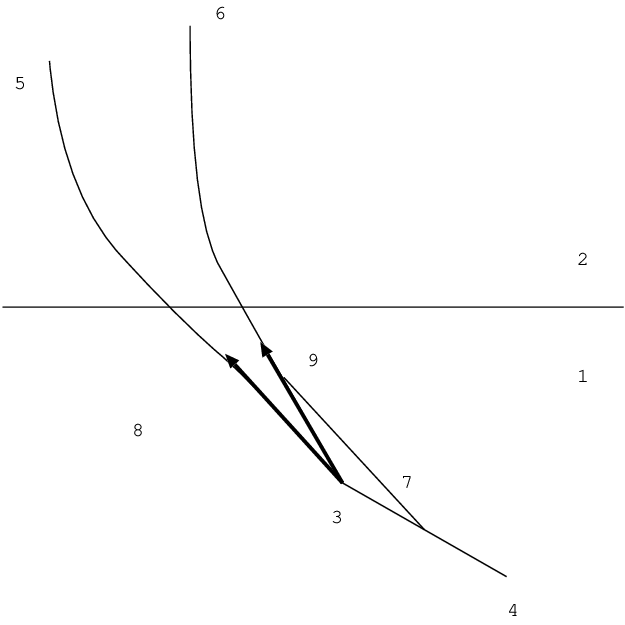} \label{pict4}\hspace{1cm}
\includegraphics[width=8cm]{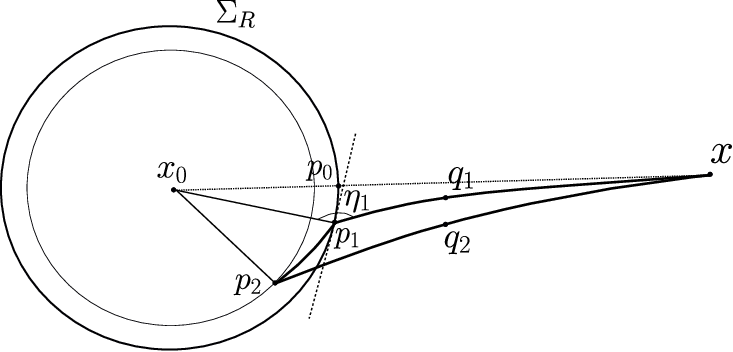} \label{figure-angle}\hspace{-1cm}

\caption{\it {\bf Left:}
Distance-minimizing paths $\gamma_i$ and a shortcut $\mu_{\tau}$ in Lemma \ref{observing same geodesic}; the
angles of $\mu$ and $\gamma_i$ at $y_1$ are close to $\pi$.
{\bf Right:} Setting of Lemma \ref{lemma-angle-separation}.
}
\end{figure}

First, we prove a result that improves \cite[Lemma 5.1]{KatsudaKL}.
Roughly speaking, the lemma states that we can identify
the directions of the distance-minimizing path from $y_1\in M$ to $x_1\in M$ up to a small error,
by considering the approximate distance function $\hat r_1$
corresponding to the point $x_1$.

\begin{lemma}\label{observing same geodesic}
Let $M$ be a closed Riemannian manifold with sectional curvature bounded below by ${\rm Sec}_M\geq -\Lambda^2$.
Let $Y\subset U\subset M$ and $y_1,y_2\in Y$.
Suppose we are given an $\epsone$-approximation $\widehat{{\mathcal R}}_Y$ of ${\mathcal R}_Y (M)$ in the sense of \eqref{Y hausdorff dist A}. 
Let $\widehat r_i \in \widehat{{\mathcal R}}_Y,$ $i=1,2$ and
let $x_i \in M$ be the points corresponding to $\widehat r_i$ (c.f.\ Definition \ref{def-xp}).
Denote by $\gamma_i(t)$ some distance-minimizing path
from $y_1$ to $x_i$, parametrized by arclength.
Assume that 
\beq
d(y_1,y_2)\geq \Lambda^{-1}, \quad d(y_1,x_i)\geq \Lambda^{-1},\quad i=1,2.
\eeq
Suppose the following is true for some $\delta\in (0,1]$ satisfying $\delta^2 \geq \epsone$:
\begin {eqnarray}\label{r1}
|\widehat r_i(y_2)-\widehat r_i(y_1)- \hatdN(y_1,y_2)|\leq \delta^{2},\quad\textrm{for } i=1,2.
\end {eqnarray}
Then there are uniform constants $\Cclosedirection,\,\Cclosedistance>1$ such that
\begin {eqnarray}\label{cl 1}
|\dot \gamma_1(0) -\dot \gamma_2(0)| \leq  \Cclosedirection \delta,
\end {eqnarray}
and
\begin {eqnarray}\label{cl 2}
\bigg|d(x_1,x_2)-|\widehat r_1(y_1)-\widehat r_2(y_1)|\bigg|\leq  \Cclosedistance
\delta.
\end {eqnarray}
\end{lemma}

\noindent {\bf Proof.}  
Due to \eqref{dH-rtox} and \eqref{Da-d estimate} and $\epsone\leq \delta^2$, the condition
(\ref{r1}) implies that
\begin {eqnarray} 
\label{r1-d}
|{\tilded}(x_i, y_2)-{\tilded}(x_i, y_1)-{\tilded}(y_1, y_2)| \leq 5 \delta^{2}, \quad i=1,2.
\end {eqnarray}
Let $\mu(t)$ be a distance-minimizing geodesic of $M$ from
 $y_1$ to  $y_2$, i.e. $\mu(0)=y_1,\,
 \mu({\tilded}(y_1, y_2))=y_2$.
 Denote by $\a_i>0$ the angle between
$\dot \gamma_i(0)$ and $\dot \mu(0)$, $i=1,2$. Now we show that $\alpha_i$ is close to $\pi$.

We restrict our attention to $i=1$; the case of $i=2$ follows in the same way. 
Let us use a shortcut argument.
Pick $z_1\in [y_1x_1]$, $z_2\in [y_1y_2]$ such that $d(y_1,z_1)=d(y_1,z_2)=\Lambda^{-1}/2$. Applying Toponogov's theorem \eqref{eq. Topog.}, we have
\beq \label{direction-shortcut}
d(z_1,z_2)\leq d(z_1,y_1)+d(y_1,z_2)-\Lambda^{-1}\CsecondToponogov (\pi-\alpha_1)^2.
\eeq
Then by the triangle inequality,
\begin{eqnarray*}
d(x_1,y_2) &\leq& d(x_1,z_1)+d(z_1,z_2)+d(z_2,y_2) \\
&\leq& d(x_1,z_1)+d(z_2,y_2)+d(z_1,y_1)+d(y_1,z_2)-\Lambda^{-1}\CsecondToponogov (\pi-\alpha_1)^2 \\
&=& d(x_1,y_1)+d(y_1,y_2)-\Lambda^{-1}\CsecondToponogov (\pi-\alpha_1)^2.
\end{eqnarray*}
On the other hand, \eqref{r1-d} gives
$$d(x_1,y_1)+d(y_1,y_2)\leq d(x_1,y_2)+5\delta^2.$$
Hence,
\begin{equation}\label{alpha_i-pi}
|\pi-\alpha_i|\leq (5\Lambda\CsecondToponogov^{-1})^{\frac12}\delta, \quad i=1,2.
\end{equation}
Denote by $\alpha$ the angle between $\dot \gamma_1(0)$ and $\dot \gamma_2(0)$. Then for suitable $\Cclosedirection>1$,
\begin{eqnarray} \label{alpha-small}
|\dot \gamma_1(0)-\dot \gamma_2(0)| \leq \alpha 
\leq |\pi-\alpha_1|+|\pi-\alpha_2| \leq \Cclosedirection \delta.
\end{eqnarray}

\smallskip
For the second claim, since $\gamma_i(d(x_i,y_1))=x_i$,
Lemma \ref{holderbound} and \eqref{alpha-small} yield that
\begin{eqnarray} \label{cl2 in d}
d(x_1,x_2) &\leq& |d(x_1,y_1)-d(x_2,y_1)|+\CexpLip \alpha \nonumber \\
&\leq& |d(x_1,y_1)-d(x_2,y_1)|+\CexpLip \Cclosedirection \delta.
\end{eqnarray}
Then the second claim follows from \eqref{dH-rtox} and the choice $\epsone\leq \delta^2$.
\hfill $\square$

\medskip

Observe that if there are several distance-minimizing paths from $y_1$ to $x_i$, the estimate (\ref{cl 1}) remains valid for each pair.

\medskip
Let $\hat r_{{i_0}}$ be given and let $x_{{i_0}}$ be such that $\|r_{x_{{i_0}}}-\hat r_{{i_0}}\|_{\ell^\infty(Y)}<\epsone$.
 We consider the elements
 $\hat r_\ell\in  \hat {\mathcal R}_Y,$
 $\ell=1,2,\dots,L$ for which
 \beq\label{e2 estimate}
 \|\hat r_{\ell}-\hat r_{{i_0}}\|_{\ell^\infty(Y)}\leq \crneighbor,
 \eeq
 where $\crneighbor$
  is a sufficiently small uniform constant to be determined later.

\smallskip
We will frequently use the following notation
\beq \label{def-Nr-section4}
N_{\epsilon}(x_0;r):=B(x_0,r+\epsilon)\setminus B(x_0,r-\epsilon),
\eeq
and we write $N_{\epsilon}(r)$ for short when the center is $x_0$.
 


%

\begin{proposition}\label{lem: properties} 
Let $M$ be a closed Riemannian manifold with sectional curvature bounded below by ${\rm Sec}_M\geq -\Lambda^2$, and $U\subset M$ be an open subset containing $B(x_0,R)$ with $R>\Lambda^{-1}$. Let $Y$ be an $\epszero$-net of $U$, and $\widehat{{\mathcal R}}_Y$
be an $\epsone$-approximation of ${\mathcal R}_Y (M)$ in the sense of \eqref{Y hausdorff dist A}.
Then the following statements hold for $0<\epszero\leq\epsone <\min\{1/16,\Lambda^{-1}/32\}$.

\begin{itemize}

\item[(1)] Let $x_{{i_0}},\, x_{\ell}$ be the points corresponding to $\hat r_{{i_0}}, \hat r_{\ell}$ (c.f.\ Definition \ref{def-xp}).
Suppose $\hat r_{{i_0}}, \hat r_{\ell}$ satisfies the following conditions:
\beq\label{bounds-section4}
\|\hat r_{\ell}-\hat r_{{i_0}}\|_{\ell^\infty(Y)}< \crneighbor<\min\{\frac{1}{4},\frac{\Lambda^{-1}}{16}\}.
\eeq
Assume we are given $y_0\in Y$ such that $d(x_0,y_0)<\epszero$.
Then there is a uniform constant $\Cclosedistance>1$ such that
\ba
d(x_\ell,x_{{i_0}})\leq  3\Cclosedistance \Big(\|\hat r_{\ell}-\hat r_{{i_0}}\|_{\ell^\infty(Y)} +3\epsone \Big)^{1/2}.
\ea

\item[(2)] The set ${X}$ is an $\epsX$-net of $M$, where $\epsX=\Cmainthree \epsone^{1/2}$ for some uniform constant $\Cmainthree>1$.
\end{itemize}
\end{proposition}

\noindent{\bf Proof.}
(1) Let us keep the parameter $R$ in the proof for clarity, and note that any dependency of $R$ in the constants can be replaced by $\Lambda$ using the condition $\Lambda^{-1}<R\leq \Lambda$. We divide into two cases depending on where $x_{{i_0}}$ lies.

\medskip
{\it $\bullet$ Case 1: $\hat r_{{i_0}}(y_0) > R/2$.}

\smallskip
Since $d(x_0,y_0)<\epszero\leq \epsone$, then $d(x_{{i_0}},x_0)>R/2-2\epsone$.
We take an arbitrary point $p\in N_{\epsone}(R/8)\cap Y$. The minimizing geodesic from $x_{{i_0}}$ to $p$ intersects with $\Sigma_{R/4}$ at some point $q'$, and we take a point $q\in Y$ such that $d(q,q')<\epszero\leq \epsone$. As a consequence, $q\in N_{\epsone}(R/4) \cap Y$. Thus \eqref{bounds-section4} and $R>\Lambda^{-1}$ yield
\beq \label{pqx-bounded-section4}
d(p,q)>\frac{R}{16},\quad d(q,x_{{i_0}})> \frac{R}{8}, \quad d(q,x_{\ell})> \frac{R}{16}.
\eeq
Moreover, by the triangle inequality,
\begin{equation} \label{satisfy-almost-minimizing-section4}
d(p,x_{{i_0}})=d(p,q')+d(q', x_{{i_0}})\geq d(p,q)+d(q,x_{{i_0}})-2\epsone.
\end{equation}

From \eqref{satisfy-almost-minimizing-section4}, \eqref{dH-rtox} and \eqref{Da-d estimate}, we see that
\beq\label{q-condition-section4}
|\hat r_{{i_0}}(p)-\hat r_{{i_0}}(q)-\hatd(p,q)|<6\epsone.
\eeq
Then pass to $\hat r_{\ell}$,
\beq\label{q-condition-l-section4}
|\hat r_{\ell}(p)-\hat r_{\ell}(q)-\hatd(p,q)|< 2\|\hat r_{\ell}-\hat r_{{i_0}}\|_{\ell^\infty(Y)}+6\epsone.
\eeq
Hence the assumptions of Lemma \ref{observing same geodesic} are satisfied with 
$$\delta^2=2\|\hat r_{\ell}-\hat r_{{i_0}}\|_{\ell^\infty(Y)}+6\epsone,$$
which satisfies $\delta^2>\epsone$ and $\delta<1$ when $\crneighbor<1/4$. Thus,
\begin{eqnarray*}
d(x_{{i_0}},x_{\ell}) &\leq&  |\hat r_{{i_0}}(q)-\hat r_{\ell}(q)|+\Cclosedistance \Big(2\|\hat r_{\ell}-\hat r_{{i_0}}\|_{\ell^\infty(Y)}+6\epsone \Big)^{1/2} \\
&\leq& 3\Cclosedistance (\|\hat r_{\ell}-\hat r_{{i_0}}\|_{\ell^\infty(Y)} +3\epsone)^{1/2}.
\end{eqnarray*}

\medskip
{\it $\bullet$ Case 2: $\hat r_{{i_0}}(y_0) \leq R/2$.}

\smallskip
In this case, $d(x_{{i_0}},x_0)\leq R/2+2\epsone$. To keep distances bounded away from zero, one can choose points $p,q$ from the outer layer $B(x_0,R)\setminus B(x_0,3R/4)$. More precisely,
we take an arbitrary point $p\in N_{\epsone}(R)\cap Y$. The minimizing geodesic from $x_{{i_0}}$ to $p$ intersects with $\Sigma_{3R/4}$ at some point $q'$, and we take a point $q\in Y$ such that $d(q,q')<\epszero\leq \epsone$. As a consequence, $q\in N_{\epsone}(3R/4) \cap Y$. Thus the bounds in \eqref{pqx-bounded-section4} still hold. Then the exact proof of Case 1 works in this case.
This concludes the proof of the first claim.

\medskip
(2) For any $x\in M$, there exists $\hat r \in \widehat{{\mathcal R}}_Y$ such that $\|r_{x}-\hat r\|_{\ell^\infty(Y)}\leq \epsone$ by \eqref{Y hausdorff dist A}. Let $x'\in M$ be a point corresponding to $\hat r$, i.e. satisfying $\|r_{x'}-\hat r\|_{\ell^\infty(Y)}\leq \epsone$. Then it follows that 
\begin{equation}\label{close-xx'}
\|r_{x}-r_{x'}\|_{\ell^\infty(Y)}\leq 2\epsone.
\end{equation}

If $x\in U$, then there exists $y\in Y$ such that $d(x,y)<\epszero\leq \epsone$, and hence $d(x',y)<3\epsone$ by \eqref{close-xx'}. Thus $d(x,x')<4\epsone$.

\smallskip
If $x\in M\setminus U$, we take an arbitrary point $p\in N_{\epsone}(R/4)\cap Y$. The minimizing geodesic from $x$ to $p$ intersects with $\Sigma_{R/2}$ at some point $q'$, and we take a point $q\in Y$ such that $d(q,q')<\epszero\leq \epsone$. As a consequence, $q\in N_{\epsone}(R/2)\cap Y$. Similarly as (1), by the triangle inequality,
\begin{equation} \label{satisfy-almost-minimizing-section4-2}
d(p,x) \geq d(p,q)+d(q,x)-2\epsone,
\end{equation}
and by \eqref{close-xx'},
\begin{eqnarray} \label{satisfy-almost-minimizing-section4-3}
d(p,x') \geq d(p,x)-2\epsone & \geq& d(p,q)+d(q,x)-4\epsone \nonumber \\
&\geq & d(p,q)+d(q,x')-6\epsone.
\end{eqnarray}
Thus the condition (\ref{r1-d}) is satisfied. Since $d(p,q),\ d(q,x), \ d(q,x')$ are all bounded below by $R/8$ by construction, one can repeat the proof of Lemma \ref{observing same geodesic} from \eqref{satisfy-almost-minimizing-section4-2} and \eqref{satisfy-almost-minimizing-section4-3}. In the end, we get the same conclusion as \eqref{cl2 in d}, namely
\begin{eqnarray}
d(x,x') \leq |d(x,q)-d(x',q)|+\Cclosedistance \epsone^{1/2},
\end{eqnarray}
which shows $d(x,x')\leq 2\epsone+\Cclosedistance \epsone^{1/2}$ by \eqref{close-xx'}.
\hfill$\square$\medskip

\smallskip
\noindent {\bf Remark.} One also has the other direction of Proposition \ref{lem: properties}(1):
\beq
\label{d-to-r} 
\|\hat r_{\ell}-\hat r_{{i_0}}\|_{\ell^\infty(Y)}\leq d(x_\ell,x_{{i_0}})+2 \epsone,
\eeq
which holds without any assumption on the parameters.
The proof is straightforward by the triangle inequality.

\section{Existence of a frame associated with measurement points}
\label{sec-frame}

We start with the following observation on geodesics connecting a point $x\in M$ to points near spheres $\Sigma_r:=\partial B(x_0,r)$. We denote 
\beq \label{def-Nr}
N_{\epsilon}(r):=B(x_0,r+\epsilon)\setminus B(x_0,r-\epsilon).
\eeq
For $x,y\in M$, we denote $|xy|=d(x,y)$, and denote by $[xy]$ a distance-minimizing geodesic from $x$ to $y$. 
The angle between $[xy]$ and $[xz]$ at $x$ is denoted by $\angle yxz$.

%

\begin{lemma} \label{lemma-angle-separation}
Let $M$ be a closed Riemannian manifold with diameter bound $\hbox{diam}(M) \leq \Lambda$ and sectional curvature bounded by $ \textrm{Sec}_M \geq -\Lambda^2$. Let $B(x_0,R)$ be an open ball with $R<{\rm inj}(M)/2$.
Then there exist uniform constants $\hat{\varepsilon},\cvdetneighbor,\Cangleseparation>0$ explicitly depending on $\Lambda,R$, such that the following holds for $0<\varepsilon<\hat{\varepsilon}$.

Given a point $x$ with $d(x,x_0)\geq 2R$, take $p_0$ to be a nearest point in $B(x_0,R)$ from $x$.
Let $p_1, p_2\in N_{\varepsilon}(R)$ such that 
\begin{equation}\label{condition-separation}
|p_1 p_2|\geq \Cnetseparation \varepsilon,\quad |p_0 p_1|<\cvdetneighbor,\quad |p_0 p_2|<\cvdetneighbor,
\end{equation}
for some $\Cnetseparation\geq 5$. Then 
\begin{itemize}
\item[(1)] $\angle p_1 x p_2 >\Cangleseparation \Cnetseparation \varepsilon$.

\item[(2)] Let $q_1\in [x p_1]$, $q_2\in [x p_2]$ such that $|x q_1| > R/2$ and $|x q_2| > R/2$. Then $|q_1 q_2|> R \Cangleseparation\Cnetseparation\varepsilon/4$.
\end{itemize}
\end{lemma}

\noindent
{\bf Proof.} (1)
Suppose $|x p_1| \leq |x p_2|$. By the triangle inequality,
\begin{eqnarray*}
|x x_0|=|x p_0|+ R \geq |x p_1|-\cvdetneighbor+ |x_0 p_1|-\varepsilon.
\end{eqnarray*}
Then by a shortcut argument using Lemma \ref{lem: Topog estimates general}(2), similar to \eqref{estimate-eta}, we have
\beq \label{estimate-eta-section5}
(\pi-\eta_1)^2\leq \CsecondToponogov^{-1} (\cvdetneighbor+\varepsilon),
\eeq
where $\eta_1 \in [0,\pi]$ is the angle between $[p_1 x]$ and $[p_1 x_0]$ at $p_1$. We choose $\cvdetneighbor,\varepsilon$ such that $\pi-\eta_1 \in [0,\pi/6]$.

Moreover, we can choose sufficiently small $\cvdetneighbor$ depending on $\Lambda,R$, such that the angle $\angle x_0 p_1 p_2$ is bounded below by $\pi/3$. This can be proved as follows. Applying Lemma \ref{lem: Topog estimates general}(1) to the triangle $x_0 p_1 p_2$, one has
\begin{eqnarray*}
2\varepsilon \geq |x_0 p_1|-|x_0 p_2| \geq |p_1 p_2| \cos(\angle x_0 p_1 p_2) - \CfirstToponogov \max\{\Lambda,R^{-1}\} |p_1 p_2|^2.
\end{eqnarray*} 
Then using the condition \eqref{condition-separation},
\begin{equation}
\cos(\angle x_0 p_1 p_2) \leq 2\varepsilon |p_1 p_2|^{-1}+ \CfirstToponogov \max\{\Lambda,R^{-1}\} |p_1 p_2|
< \frac{2}{5}+C(\Lambda,R) 2\cvdetneighbor.
\end{equation}
The above shows that we can choose sufficiently small $\cvdetneighbor,\varepsilon$, such that the angle $\angle x p_1 p_2$ satisfies
\begin{equation} \label{bound-eta}
\angle x p_1 p_2 \leq (\frac{\pi}{2}+\pi-\eta_1)+ (\frac{\pi}{2}-\frac{\pi}{3}) \leq \frac56 \pi.
\end{equation}

Next, we apply Lemma \ref{lem: Topog estimates general}(1) to the triangle $x p_1 p_2$ and use \eqref{bound-eta} as follows:
\begin{eqnarray*}
|x p_2| &\leq& |x p_1|-|p_1 p_2|\cos (\angle x p_1 p_2) + \CfirstToponogov \max\{\Lambda,R^{-1}\} |p_1 p_2|^2 \nonumber \\
&\leq& |x p_1|+\frac{\sqrt{3}}{2} |p_1 p_2|+ \CfirstToponogov \max\{\Lambda,R^{-1}\} |p_1 p_2|^2.
\end{eqnarray*}
Using $|p_1 p_2|< 2\cvdetneighbor$, we obtain for sufficiently small $\cvdetneighbor$,
\begin{equation} \label{p1p2-upper}
|x p_2| -|x p_1| \leq \frac{\sqrt{3}}{2}|p_1 p_2|+ \CfirstToponogov \max\{\Lambda,R^{-1}\} |p_1 p_2|^2 < \frac{9}{10}|p_1 p_2|.
\end{equation}
Thus by Lemma \ref{holderbound} and \eqref{p1p2-upper}, we obtain
\begin{eqnarray} \label{p1p2CA}
|p_1 p_2| &\leq& |x p_2|-|x p_1| + \CexpLip \angle p_1 x p_2 \nonumber \\
&<& \frac{9}{10}|p_1 p_2| + \CexpLip \angle p_1 x p_2,
\end{eqnarray}
which proves part (1) due to $|p_1 p_2| \geq \Cnetseparation\varepsilon$, with $\Cangleseparation=(10\CexpLip)^{-1}$.

\medskip
(2) Consider the triangle $xp_1 p_2$ and denote
its comparison triangle in the hyperbolic plane of constant sectional curvature $-\Lambda^2$ by $\overline{x} \overline{p_1} \overline{p_2}$, i.e. $|x p_1|=|\overline{x} \overline{p_1}|$, $|x p_2|=|\overline{x} \overline{p_2}|$, $|p_1 p_2|=|\overline{p_1} \overline{p_2}|$. The angle at $\overline{x}$ of this comparison triangle is denoted by $\overline{\angle} p_1 x p_2$.
Since the inequalities \eqref{p1p2-upper} and \eqref{p1p2CA} are also valid for the comparison triangle $\overline{x} \overline{p_1} \overline{p_2}$, it also holds that $\overline{\angle} p_1 x p_2 >\Cangleseparation \Cnetseparation \varepsilon$ in the comparison triangle.
By the angle-sidelength monotonicity version of Toponogov's theorem (e.g. \cite[Theorem 7.3.2]{AKP}), the function $(|xq_1|,|x q_2|)\mapsto \overline{\angle} q_1 x q_2$ is decreasing in both arguments. Hence 
\begin{equation}\label{comparison-angle-mono}
\overline{\angle} q_1 x q_2 \geq \overline{\angle} p_1 x p_2>\Cangleseparation \Cnetseparation \varepsilon.
\end{equation}
Then in the comparison triangle $\overline{x} \overline{q_1} \overline{q_2}$, by \eqref{the law of cosines} and \eqref{comparison-angle-mono},
\begin{eqnarray*}
\cosh({\Lambda}{  |\overline q_1\overline q_2|})&=& \cosh ({\Lambda}{ |\overline x\overline q_1|})\,\cosh ({\Lambda}{|\overline x\overline q_2|})-
 \sinh ({\Lambda}{ |\overline x\overline q_1|})\,\sinh({\Lambda}{ |\overline x\overline  q_2|})\,\cos(\overline{\angle} q_1 x q_2)\\
&=& \cosh \Big(\Lambda |x q_1|-\Lambda |x q_2| \Big) + \sinh ({\Lambda}{ |x q_1|})\,\sinh({\Lambda}{ |x q_2|})\,\Big(1-\cos(\overline{\angle} q_1 x q_2) \Big) \\
&>& 1+ \frac{\Lambda^2 R^2}{4} \Big(1-\cos(\Cangleseparation \Cnetseparation \varepsilon) \Big) \\
&\geq& 1+ \frac{\Lambda^2 R^2}{4} \frac{1}{4} (\Cangleseparation \Cnetseparation \varepsilon)^2.
\end{eqnarray*}
Since $\cosh({\Lambda}{  |\overline q_1\overline q_2|})=\cosh({\Lambda}{  |q_1 q_2|})\leq 1+\Lambda^2 |q_1 q_2|^2$ for small $|q_1 q_2|$, the second claim follows.
\hfill$\square$\medskip

Given a point $x\in M$ outside the ball $B(x_0,2R)$ and an $\varepsilon$-net of $N_{\varepsilon}(R)$, we consider all the unit initial vectors of minimizing geodesics from $x$ to points in the $\varepsilon$-net. The following lemma shows that there exist $n$ such unit vectors so that the corresponding determinant is bounded away from zero.

\begin{lemma}\label{lemma-frame}
Let $M$ be a closed Riemannian manifold with diameter bound $\hbox{diam}(M) \leq \Lambda$ and sectional curvature bounded by $|\textrm{Sec}_M|\leq \Lambda^2$. Let $B(x_0,R)$ be an open ball for $R<{\rm inj}(M)/2$, and $Y$ be an $\varepsilon$-net of $B(x_0,R)$.
Given a point $x$ with $d(x,x_0)\geq R$, take $p_0$ to be a nearest point in $B(x_0,R/2)$ from $x$.
Then for every $\Cnetseparation\geq 5$, there are uniform constants $\hat{\varepsilon},\cvdetlower,\cvdetneighbor>0$ explicitly depending on $n,\Lambda,R,\Cnetseparation$, such that the following holds for $0<\varepsilon<\hat{\varepsilon}$.

\smallskip
We can find a separated set $\{p_{j(l)}:\ l=1,2,\dots, L\}\subset  Y\cap N_{\varepsilon}(R/2)$ satisfying 
\beq\label{sparse condition}
|p_{j(l)} p_{j(m)}|\geq \Cnetseparation \varepsilon, \quad\hbox{for $l\neq m$, and}\quad  |p_0 p_{j(l)}|<\cvdetneighbor,
\eeq
such that the following statement is true: there exist $n$ unit vectors $w_{l_1},\dots,w_{l_n}$ such that the volume
of the simplex in $T_{x}M$ with the vertices $0,w_{l_1},\dots,w_{l_n}$ is 
larger than $\cvdetlower$,
where $w_{l} \in S_{x}M$, $l=1,2,\dots,L$ is the unit initial vector of a minimizing geodesic in $M$ from ${x}$ to $p_{j(l)}$.
\end{lemma}

\noindent
{\bf Proof.}
Let us take $\Cnetseparation=10$ for the sake of argument.
First, we choose a suitable separated set.
We can choose a separated set of $Y\cap N_{\varepsilon}(R/2)$ satisfying the condition \eqref{sparse condition} such that it is also a $11\varepsilon$-net. Namely, first pick any point, say $p_1$, from $Y\cap N_{\varepsilon}(R/2)$, and then take all points $p_2,\dots,p_{l}$ in $Y\cap N_{\varepsilon}(R/2)$ with distance at least $10\varepsilon$ away from $p_1$. For the second step, pick all points with distance at least $10\varepsilon$ from all of $p_1,p_2,\dots,p_l$. Repeat this procedure and the procedure stops in finite steps. The chosen points form a $11\varepsilon$-net of $N_{\varepsilon}(R/2)$. Indeed, for any $z\in N_{\varepsilon}(R/2)$, there exists a point $y\in Y\cap N_{\varepsilon}(R/2)$ such that $d(y,z)<\varepsilon$ by the $\varepsilon$-net condition for $Y$. If $y$ is also chosen by our procedure above, then we are done. If $y$ is not chosen, it means that there must be at least one chosen point $p_m$ such that $d(y,p_m)< 10\varepsilon$ otherwise $y$ would have been chosen. In such case, we have $d(z,p_m)\leq d(z,y)+d(y,p_m)<11\varepsilon$.

\smallskip
The separated set of our choice above satisfying the condition \eqref{sparse condition} has cardinality at least 
\begin{equation} \label{L-bound}
L\geq C(n,\Lambda)\cvdetneighbor \varepsilon^{-(n-1)}.
\end{equation}
Indeed, the separated set that we chose above is a $11\varepsilon$-net of the set $N_{\varepsilon}(R/2)\cap B(p_0,\cvdetneighbor)$, and the latter set has volume bounded below. More precisely,
\begin{equation*}
C(n) (11 \varepsilon)^{n} L \geq \textrm{vol}_n \Big(N_{\varepsilon}(R/2)\cap B(p_0,\cvdetneighbor)\Big) \geq C(n,\Lambda)\cvdetneighbor \varepsilon,
\end{equation*}
which yields the lower bound \eqref{L-bound}.

\smallskip
On the other hand, by Lemma \ref{lemma-angle-separation}(1), the set $\{w_l\}_{l=1}^L \subset S_{x} M$ of unit vectors is $\Cangleseparation \varepsilon$-separated. We claim that we can choose sufficiently small $h$ explicitly depending on $n,\Lambda,\cvdetneighbor$, such that for any unit vector $\xi\in S_x M\setminus \{0\}$,
\begin{equation} \label{vl-notsubset}
\{w_l\}_{l=1}^L \not\subset \{v \in S_{x}M: \bra v,\xi\cet_g\in (-h,h)\}.
\end{equation}
The claim \eqref{vl-notsubset} can be proved as follows. Suppose \eqref{vl-notsubset} is not true: there exists $\xi$ such that $\{w_l\}_{l=1}^L \subset A(\xi,h)$, where $A(\xi,h)$ denotes the set on the right-hand side of \eqref{vl-notsubset}. Since $\{w_l\}_{l=1}^L \subset S_{x} M$ is $\Cangleseparation\varepsilon$-separated, this implies
$$C(n) (\Cangleseparation \varepsilon)^{n-1} L \leq \textrm{vol}_{n-1}(A(\xi,h)).$$
However, the $(n-1)$-dimensional volume of $A(\xi,h)$ is bounded above by $C(n)h$. Hence \eqref{L-bound} yields
$$C(n,\Lambda)\cvdetneighbor \leq C(n) (\Cangleseparation \varepsilon)^{n-1} L \leq \textrm{vol}_{n-1}(A(\xi,h)) \leq C(n)h,$$
which cannot be true if $h<C(n,\Lambda)\cvdetneighbor$. The claim \eqref{vl-notsubset} is proved.

\smallskip
Thus we can construct the desired $w_{l_1},\dots,w_{l_n}$ as follows.
Let $w_{l_1}$ be chosen arbitrarily. Then when $w_{l_1},\dots,w_{l_k}$ are chosen,
let $\xi$ be a unit vector orthogonal to  $w_{l_1},\dots,w_{l_k}$. From the claim \eqref{vl-notsubset}, there exists
$w_{l_{k+1}}$ such that $w_{l_{k+1}}\not \in A(\xi,h)$ with $h$ properly chosen as above.
In other words, the angle between $w_{l_{k+1}}$ and the linear subspace spanned by $w_{l_1},\dots,w_{l_k}$ is bounded below by $h$. Hence the simplex with vertices $0,w_{l_1},\dots,w_{l_k},w_{l_{k+1}}$ has volume bounded below by $C(n,h)$.
\hfill$\square$\medskip

In other words, there exist $n$ unit vectors in $\{w_l\}_{l=1}^L$ such that the corresponding determinant  $\det([\langle w_k,w_m \rangle_g]_{k,m=1}^n])>\cvdetlower$.

\smallskip
The next lemma is a technical modification of the previous lemma.

\begin{lemma}\label{lemma-frame-q}
Under the setting of Lemma \ref{lemma-frame}, there are uniform constants $\hat{\varepsilon},\cvdetlower,\cvdetneighbor,\cvdetwidth>0$ explicitly depending on $n,\Lambda,R$, such that the following holds for $0<\varepsilon<\hat{\varepsilon}$.

Let $\{p_{j(l)}:\ l=1,2,\dots, L\}\subset  Y\cap N_{\varepsilon}(R/2)$ be a choice of points satisfying the condition \eqref{sparse condition} for sufficiently large $\Cnetseparation$. For each $l$, suppose a minimizing geodesic $[x p_{j(l)}]$ intersects with $\partial B(R/2+\cvdetwidth)$ at $q'_{j(l)}$, and we take $q_{j(l)}\in Y\cap N_{\varepsilon}(R/2+\cvdetwidth)$ to be a point in $Y$ such that $|q'_{j(l)} q_{j(l)}|< \varepsilon$.

\smallskip
Then we have
\beq\label{sparse condition-q}
|q_{j(l)} q_{j(m)}|\geq 10 \varepsilon, \quad\hbox{for $l\neq m$, and}\quad  |p_0 q_{j(l)}|<2\cvdetneighbor.
\eeq
\noindent As a consequence, there exist $n$ unit vectors $v_{l_1},\dots,v_{l_n}$ such that the volume
of the simplex in $T_{x}M$ with the vertices $0,v_{l_1},\dots,v_{l_n}$ is 
larger than $\cvdetlower$,
where $v_{l} \in S_{x}M$, $l=1,2,\dots,L$ is the unit initial vector of a minimizing geodesic in $M$ from ${x}$ to $q_{j(l)}$.
\end{lemma}

\noindent
{\bf Proof.}
Fix the parameters $\hat{\varepsilon},\cvdetneighbor$ as chosen in Lemma \ref{lemma-frame}.
Let us consider two points $p_1,p_2$ in the maximal set constructed in Lemma \ref{lemma-frame} with $\Cnetseparation$ to be determined later. The intersection of $[xp_1]$, $[x p_2]$ with $\partial B(R/2+\cvdetwidth)$ is $q'_1,q'_2$. By Lemma \ref{lemma-angle-separation}(2), $|q'_1 q'_2|>R\Cangleseparation \Cnetseparation \varepsilon/4$. When we take the points $q_1,q_2$ in the $\varepsilon$-net, we have 
\begin{equation}
|q_1 q_2|>\Big(\frac{R\Cangleseparation \Cnetseparation}{4}-2 \Big) \varepsilon.
\end{equation}
Then we can choose $\Cnetseparation=48R^{-1}\Cangleseparation^{-1}$ so that $|q_1 q_2|>10 \varepsilon$.
Furthermore, since the incident angle is bounded by \eqref{estimate-eta-section5}, we can choose sufficiently small $\cvdetwidth>0$ such that $|p_1 q'_1|$ is bounded by $\cvdetneighbor/2$. Hence $|p_0 q_1|<2\cvdetneighbor$. This proves the claim \eqref{sparse condition-q}.

Let $q_0$ be the intersection of $[xp_0]$ with $\partial B(R/2+\cvdetwidth)$, and hence $q_0$ is a nearest point in $B(R/2+\cvdetwidth)$ from $x$. Moreover, $|q_0 q_1|<3\cvdetneighbor$ if we choose $\cvdetwidth<\cvdetneighbor$. Thus the condition \eqref{condition-separation} in Lemma \ref{lemma-angle-separation} is satisfied by the triangle $x q_1 q_2$, which gives $\angle q_1 x q_2>\Cangleseparation \varepsilon$. Observe that the total number of points $\{q_{j(l)}\}$ is equal to the total number $L$ of points $\{p_{j(l)}\}$, because of $\{q_{j(l)}\}$ being $10\varepsilon$-separated.
The total number $L$ is bounded below by \eqref{L-bound}. Then the same argument yields \eqref{vl-notsubset} for $\{v_l\}_{l=1}^L$, and the second claim follows from the last part of the proof of Lemma \ref{lemma-frame}.
\hfill$\square$\medskip

\begin{remark} \label{remark-inner}
Lemma \ref{lemma-angle-separation}(1) is also true if $x$ is inside the ball, say $x\in B(x_0,R/2)$, assuming two-side bounds on the sectional curvature $|\textrm{Sec}_M|\leq \Lambda^2$. The proof is similar and can be found in Lemma \ref{lemma-angle-separation-inner}. As a consequence, Lemma \ref{lemma-frame} is still valid by the same argument if $x\in B(x_0,R/2)$, in which case we can find the desired separated set in $Y\cap N_{\varepsilon}(R)$.
We will use this observation in the next section.
\end{remark}

\section{Local reconstructions from partial distance data}
\label{sec-mainproof}

This section is the proof of Theorem \ref{main 1} and consequently Corollary \ref{coro 1}.
Let $x_i\in M$ be the points corresponding to $\hat r_i\in \widehat{\mathcal R}_Y$, $i=1,2,\dots,I$, i.e. satisfying \eqref{dH-rtox}. Let us fix one element $\hat r_{{i_0}}\in \widehat{\mathcal R}_Y$ and the corresponding point $x_{{i_0}}\in M$. The basic idea of the proof is to find appropriate points in $B(x_0,R)\subset U$, and apply geometric lemmas in previous sections to approximate the inner product. One important point is to keep distances of points bounded away from zero, as required by previous lemmas. This is possible because we assumed the knowledge of a point $y_0\in Y$ such that $d(x_0,y_0)<\epszero$. This assumption enables us to determine where $x_{{i_0}}$ lies in reference to $B(x_0,R)$ up to a small error. Furthermore, it is possible to use only part of all measurement points in the ball. This allows us to simply take, for example $R=(4\Lambda)^{-1}$, and only consider the measurement points in $Y$ within this smaller ball. In the proof, we keep the parameter $R$ for clarity, and note that any dependency of $R$ in the constants can be replaced by $\Lambda$.

\smallskip
We divide the proof of Theorem \ref{main 1} into two cases depending on where $x_{{i_0}}$ lies, in a similar way as we considered in Proposition \ref{lem: properties}(1). We will focus on the first case, as the second case is a simple modification from the first case.

\bigskip
\textbf{Case 1:} $\hat r_{{i_0}}(y_0) > R/2$.

\medskip
Let us set 
\begin{equation}\label{epsilon-bound}
0<\epszero\leq\epsone <\min \big\{\frac{1}{16},\frac{R}{32} \big\}.
\end{equation}
We consider the elements $\hat r_{\ell}\in \widehat{\mathcal R}_Y$ in the neighborhood of $\hat r_{{i_0}}$,
\beq \label{c2neighbor}
\|\hat r_{\ell}-\hat r_{{i_0}}\|_{\ell^\infty(Y)}<\rho_0,\quad \crneighbor< \min\{\frac{1}{4},\frac{R}{16}\},
\eeq
with the parameter $\crneighbor$ to be determined later.
Let $x_{\ell}\in M$ be a point corresponding to $\hat r_{\ell}$. 


\medskip
{\it $\bullet$ Step 1: Applying first variation formula.}

\smallskip

The condition $\hat r_{{i_0}}(y_0) > R/2$ implies that $d(x_{{i_0}},x_0)>R/2-2\epsone$. 
For an arbitrary point $p\in N_{\epsone}(R/8)\cap Y$, we pick a point $q\in N_{\epsone}(R/4)\cap Y$ such that
\beq\label{q-condition-step1}
\Big|\hat r_{{i_0}}(p)-\hat r_{{i_0}}(q)-\hatd(p,q) \Big|<6\epsone.
\eeq

Recall $\hatd$ defined in \eqref{def-DY}.
It is clear that condition \eqref{q-condition-step1} can be tested using the given data $\widehat{\mathcal R}_Y$ only.
Observe that the set of points $q$ (for each $p$) satisfying \eqref{q-condition-step1} is nonempty. This is because the point in $Y$ within $\epszero$-distance from the intersection of a minimizing geodesic $[x_{{i_0}}p]$ with $\Sigma_{R/4}$ satisfies the condition \eqref{q-condition-step1} due to \eqref{satisfy-almost-minimizing-section4} and \eqref{q-condition-section4}.
In particular, the following bounds are valid:
\beq \label{pqx-bounded}
d(p,q)>\frac{R}{16},\quad d(q,x_{{i_0}})> \frac{R}{8}, \quad d(q,x_{\ell})> \frac{R}{16}.
\eeq
Moreover, by \eqref{q-condition-step1}, \eqref{dH-rtox}, \eqref{Da-d estimate}, and \eqref{Da-d estimate identification} we see that
\beq \label{satisfy-almost-minimizing}
\Big|d(x_{{i_0}},p)-d(x_{{i_0}},q)-d(p,q) \Big| \leq 10\epsone.
\eeq
Hence the assumptions of Proposition \ref{lem: inner product lemma} are satisfied when $\epsone<\hat{\delta}/10$.

By Proposition \ref{lem: properties}(1), we know
\beq
d(x_{\ell},x_{{i_0}}) \leq 3\Cclosedistance \big(\|\hat r_{\ell}-\hat r_{{i_0}}\|_{\ell^\infty(Y)}+3\epsone \big)^{\frac12}.
\eeq
Hence if we choose $\crneighbor$ in \eqref{c2neighbor} such that 
\beq
3\Cclosedistance(\crneighbor+3\epsone)^{1/2}<\ccalphafirst^2\, ,
\eeq
we can apply Proposition \ref{lem: inner product lemma} (with $x=x_{{i_0}},\, y=x_{\ell}$) and obtain
\beq
\bigg| \langle \xi_{\ell}, v \rangle_g - \Big(|x_{{i_0}}q|-|x_{\ell} q| \Big)\bigg|\leq \betadissecond |x_{{i_0} } x_{\ell}| \Big(\epsone^{1/4}+|x_{{i_0}} x_{\ell}|^{1/3} \Big),
\eeq
where $\xi_{\ell}=\exp_{x_{{i_0}}}^{-1}(x_{\ell})$ (of length $|x_{{i_0}}x_l|$), and $v$ is the unit initial vector of $[x_{{i_0}} q]$.
Then using \eqref{dH-rtox}, we have
\beq \label{eq-xil}
\hspace{10mm} \bigg| \langle \xi_{\ell}, v \rangle_g - \Big(\hat r_{{i_0}}(q)-\hat r_{\ell}(q) \Big)\bigg|\leq \betadissecond |x_{{i_0} } x_{\ell}| \Big(\epsone^{1/4}+|x_{{i_0}} x_{\ell}|^{1/3} \Big) +2\epsone.
\eeq

\medskip
{\it $\bullet$ Step 2: Finding the length $|x_{{i_0}} x_{\ell}|$.}

\smallskip
We aim to construct an approximate inner product to the actual one, i.e. the first term in \eqref{eq-xil}. 
However, the first term $\langle \xi_{\ell}, v \rangle_g$ involves the length $|\xi_{\ell}|_g=|x_{{i_0}} x_{\ell}|$, which cannot be exactly computed from the data $\widehat{{\mathcal R}}_Y$. What we can do is to use Lemma \ref{observing same geodesic} to approximate it.

Let us take a small parameter $s$ whose value is determined later:
\beq \label{def-s}
s\in (\epsone^{1/2},\frac{\crneighbor}{2}).
\eeq
Let $p\in N_{\epsone}(R/8)\cap Y$ be arbitrary, and let $q\in N_{\epsone}(R/4)\cap Y$ be chosen according to the condition \eqref{q-condition-step1}.
Now we choose an element $\hat r_{\ell}\in \widehat{\mathcal R}_Y$ satisfying \eqref{c2neighbor} such that the following conditions \eqref{criteria1} and \eqref{criteria2} hold:
\begin{eqnarray}
&&\Big|\hat r_{\ell}(p)-\hat r_{\ell}(q)-\hatd(p,q) \Big|\leq 9\epsone, \label{criteria1} \\
&&\Big|\hat r_{{i_0}}(q)-(\hat r_{\ell}(q)+s) \Big|\leq 9\epsone. \label{criteria2}
\end{eqnarray}
In fact, as we will later specify our choice $s=\epsone^{3/8}$, 
we can actually choose this $\hat r_{\ell}$ from $\widehat{{\mathcal R}}_Y\cap\mathcal B_\infty(\hat r_{{i_0}},\epsone^{1/4})$.
Indeed,
considering \eqref{d-to-r}, 
it is straightforward to check that the element in $\widehat{{\mathcal R}}_Y\cap\mathcal B_\infty(\hat r_{{i_0}},\epsone^{1/4})$ corresponding to $\gamma_{x_{{i_0}},w}(s)$ satisfies these two conditions above, where $w$ is the unit initial vector of $[x_{{i_0}}p]$. Essentially, these conditions can be understood as a test to search for $\gamma_{x_{{i_0}},w}(s)$ up to a small error, see Figure $6$.

\smallskip
Next, let us discuss the properties of such $\hat r_{\ell}$ satisfying the criteria \eqref{criteria1} and \eqref{criteria2}.
Due to Lemma \ref{observing same geodesic}, \eqref{q-condition-step1} and \eqref{criteria1} imply that
\beq \label{ds-close}
\bigg| |x_{{i_0}} x_{\ell}|- |\hat r_{{i_0}}(q)-\hat r_{\ell}(q)| \bigg| \leq 3\Cclosedistance \epsone^{1/2}.
\eeq
Hence \eqref{ds-close}, \eqref{criteria2} and \eqref{dH-rtox} yield for some suitable $\Cs>1$,
\beq
&& \Big| |x_{{i_0}} x_{\ell}|-s \Big| \leq \Cs \epsone^{1/2}, \label{ds-close-2} \\
&& \Big| |x_{{i_0}} x_{\ell}|- (|x_{{i_0}} q|-|x_{\ell} q|) \Big| \leq \Cs \epsone^{1/2}. \label{ds-close-3}
\eeq

Furthermore, let $\gamma_{x_{{i_0}},v}(\cdot)$ be a minimizing geodesic from $x_{{i_0}}$ to $q$, and $\xi_{\ell}=\exp_{x_{{i_0}}}^{-1}(x_{\ell})$. We claim that for some uniform constant $\Csl>1$,
\beq \label{d-gammas-xl}
d(\gamma_{x_{{i_0}},v}(s), x_{\ell}) \leq \Csl \epsone^{1/2}.
\eeq
This {claim can be proved as follows. Denote $z'=\gamma_{x_{{i_0}},v}(s)$
and observe that  $z'$ is on a minimizing geodesic from $x_{{i_0}}$ to $q$. Due to the fact that
$$d(x_{{i_0}},p)-d(x_{{i_0}},q) \leq d(x_{{i_0}},z')+ d(z',p)-d(x_{{i_0}},z')-d(z',q) = d(z',p)-d(z',q),$$
the inequality \eqref{satisfy-almost-minimizing} 
is still valid after replacing $x_{{i_0}}$ with $z'=\gamma_{x_{{i_0}},v}(s)$,
that is,
\beq \label{satisfy-almost-minimizing  x replaced by z'}
\Big|d(z',p)-d(z',q)-d(p,q) \Big| \leq 10\epsone.
\eeq
Thus Lemma \ref{observing same geodesic} (with $x_1=z',\ x_2=x_{\ell}$) and \eqref{criteria1} yield that
$$\Big| d(z', x_{\ell}) - |d(z', q)-d(x_{\ell},q)| \Big|\leq \Cclosedistance \epsone^{1/2}.$$
Since 
\begin{eqnarray*}
d(z', q)-d(x_{\ell},q) &=& |z' q| +|x_{{i_0}} z'| -(|x_{\ell} q|+|x_{{i_0}}x_{\ell}|) +|x_{{i_0}}x_{\ell}|-|x_{{i_0}} z'| \\
&=& |x_{{i_0}}q|- (|x_{\ell} q|+|x_{{i_0}}x_{\ell}|) +(|x_{{i_0}}x_{\ell}|-s),
\end{eqnarray*}
then the claim \eqref{d-gammas-xl} follows from  \eqref{ds-close-2} and \eqref{ds-close-3}.}

As a consequence, by applying the Rauch comparison theorem (see e.g. \cite{Pe}) for ${\rm Sec}_M\leq \Lambda^2$ (in a neighborhood of $x_{{i_0}}$), \eqref{d-gammas-xl} yields that for some $\Cxiv>1$,
\beq \label{close-xiv}
|\xi_{\ell}-sv| \leq \Cxiv \epsone^{1/2}.
\eeq

\medskip
{\it $\bullet$ Step 3: Approximating the inner product.}

\smallskip
Denote by $s_0,s_1$ the lower and upper bounds for $|x_{{i_0}} x_{\ell}|$. From \eqref{ds-close-2}, we set
\beq
s_0:=s- \Cs \epsone^{1/2}, \quad s_1:=s+ \Cs \epsone^{1/2}.
\eeq
We require that $s>2\Cs \epsone^{1/2}$ so that $s_0>0$ and $s_1<2s$. 
Observe that \eqref{ds-close-2} yields
\begin{equation}\label{ds-close-4}
\bigg| \langle \xi_{\ell},v\rangle_g - s\langle \frac{\xi_{\ell}}{|\xi_{\ell}|_g},v\rangle_g \bigg| =\bigg| |x_{{i_0}}x_{\ell}|\langle \frac{\xi_{\ell}}{|\xi_{\ell}|_g},v\rangle_g - s\langle \frac{\xi_{\ell}}{|\xi_{\ell}|_g},v\rangle_g \bigg| \leq \Cs \epsone^{1/2}.
\end{equation}
Then \eqref{eq-xil} gives
\begin{equation} 
\bigg| s\langle \frac{\xi_{\ell}}{|\xi_{\ell}|_g}, v \rangle_g - \Big(\hat r_{{i_0}}(q)-\hat r_{\ell}(q) \Big)\bigg|\leq \betadissecond s_1(\epsone^{1/4}+s_1^{1/3}) +(\Cs+2) \epsone^{1/2}.
\end{equation}
Hence dividing by $s$ gives
\begin{eqnarray*} 
\bigg| \langle \frac{\xi_{\ell}}{|\xi_{\ell}|_g}, v \rangle_g - \frac{1}{s}\Big(\hat r_{{i_0}}(q)-\hat r_{\ell}(q) \Big)\bigg| &\leq& \betadissecond s_1(\epsone^{1/4}+s_1^{1/3})s^{-1} +(\Cs+2)\epsone^{1/2} s^{-1} \nonumber \\
&\leq& 4(\betadissecond+\Cs+2) s^{-1}(\epsone^{1/2} + \epsone^{1/4}s +s^{4/3}) \label{def-E}.
\end{eqnarray*}
We obtain a good estimate when we choose 
\begin{equation} \label{choice-s}
s=\epsone^{3/8}.
\end{equation}
Note that when $\epsone$ is sufficiently small, the requirement that $s>2\Cs \epsone^{1/2}$ is satisfied. In such case, we have the estimate
\begin{equation} \label{estimate-E}
\bigg| \langle \frac{\xi_{\ell}}{|\xi_{\ell}|_g}, v \rangle_g - \frac{1}{s}\Big(\hat r_{{i_0}}(q)-\hat r_{\ell}(q) \Big)\bigg| \leq \CE \epsone^{1/8}.
\end{equation}

\begin{figure}[htbp]
\begin{center}
\psfrag{a}{$\Sigma_{R/2}$}
\psfrag{c}{$\Sigma_{R/4}$}
\psfrag{b}{$\Sigma_{R/8}$}
\psfrag{1}{\vspace{-4mm}$p_{j(k)}$}
\psfrag{2}{$q_{j(k)}$}
\psfrag{3}{$x_{{i_0}}$}
\psfrag{4}{$y'$}
\psfrag{5}{$x_{i(m)}$}
\psfrag{6}{$x_\ell$}
\includegraphics[width=7cm]{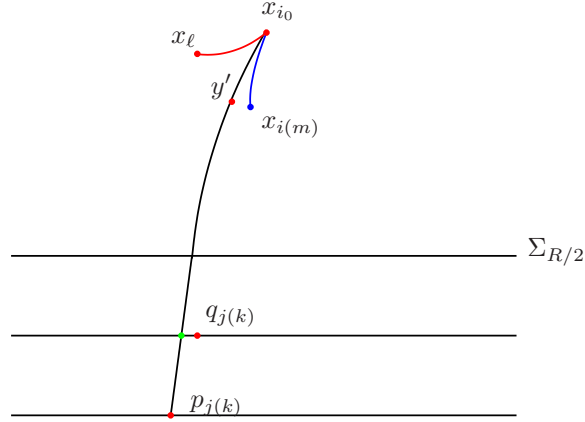} \label{pict13ac}
\end{center}
\caption{\it The geodesic  $[x_{{i_0}}p_{j(k)}]$ intersects
$\Sigma_{R/4}=\partial B(x_0,R/4)$ at a green point that is close
to the red point $q_{j(k)}$.
The other red point is $y'=\gamma_{x_{{i_0}},w}(s)$,
where $w\in S_{x_{{i_0}}}M$ is the direction vector of the geodesic $[x_{{i_0}}p_{j(k)}]$
and $s=d(x_{{i_0}},y')$.
When the criteria (\ref{criteria1}) and  (\ref{criteria2}) are satisfied for $p_{j(k)},q_{j(k)}$,
the point $x_{i(m)}$ corresponding to $\hat r_{i(m)}$ is close to $y'$.
The blue curve is the geodesic $[x_{{i_0}} x_{i(m)}]$ and
$\xi_{i(m)}=\exp^{-1}_{x_{{i_0}}}(x_{i(m)})$. Similarly, the red curve is the geodesic $[x_{{i_0}}x_{\ell}]$ and
$\xi_{\ell}=\exp^{-1}_{x_{{i_0}}}(x_{\ell})$. In the proof  we show that
we can find indexes $i(m)$, $m=1,2,\dots,n$, such that the unit vectors
$\overline{\xi}_{i(m)}:=|\xi_{i(m)}|_g^{-1}\xi_{i(m)}$ form a good basis of $T_{x_{{i_0}}}M$.
Moreover, we show that we can 
approximately compute
the inner product of the vectors $\overline{\xi}_{i(m)}$ and $\xi_\ell$,
that is, approximately find the coordinates of the points $x_\ell$ in the normal coordinates at $x_{{i_0}}$.}
\end{figure}

\medskip
{\it $\bullet$ Step 4: Approximating the metric.}

\smallskip
Let us find a proper frame in $T_{x_{{i_0}}} M$ and apply the estimate \eqref{estimate-E} to each vector in the frame to approximate the metric. First, we search for a point $p_0\in N_{\epsone}(R/8)\cap Y$ such that
\begin{equation}\label{search-p0}
\hat{r}_{{i_0}}(p_0)=\min_{y\in N_{\epsone}(R/8)\cap Y} \hat{r}_{{i_0}}(y).
\end{equation}
Let $p_{j(k)} \in N_{\epsone}(R/8)\cap Y$, $j(k)\in \{0,\dots,J\}$, $k=1,\dots,n$, be arbitrary $n$ points in $Y$ satisfying the condition \eqref{sparse condition} (with $\varepsilon=\epsone$), which we will vary later.
For each $p_{j(k)}$, we search for $q_{j(k)}\in N_{\epsone}(R/4)\cap Y$ such that \eqref{q-condition-step1} hold, i.e.
\beq \label{q-condition-i(k)}
\Big|\hat r_{{i_0}}(p_{j(k)})-\hat r_{{i_0}}(q_{j(k)})-\hatd(p_{j(k)},q_{j(k)}) \Big|<6\epsone.
\eeq
Then we choose an element $\hat r_{i(k)}\in \widehat{{\mathcal R}}_Y\cap\mathcal B_\infty(\hat r_{{i_0}},\epsone^{1/4})$, $i(k)\in \{1,\dots,I\}$ such that \eqref{c2neighbor}, \eqref{criteria1} and \eqref{criteria2} hold for $p_{j(k)},q_{j(k)}$, namely
\begin{eqnarray}
&&\Big|\hat r_{i(k)}(p_{j(k)})-\hat r_{i(k)}(q_{j(k)})-\hatd(p_{j(k)},q_{j(k)}) \Big|\leq 9\epsone, \label{criteria1-jm} \\
&&\Big|\hat r_{{i_0}}(q_{j(k)})-(\hat r_{i(k)}(q_{j(k)})+s) \Big|\leq 9\epsone. \label{criteria2-jm}
\end{eqnarray}
By the discussion following \eqref{criteria2}, the set of $\hat r_{i(k)}$ satisfying the conditions above is non-empty, as we set $s=\epsone^{3/8}$.

\smallskip
Now we apply the estimate \eqref{estimate-E} for each index $j(k)$.  
Denote by $v_{j(k)}$ the unit initial vectors of minimizing geodesics $[x_{{i_0}} q_{j(k)}]$. For every $k,m\in \{1,\dots,n\}$, we have
\begin{equation} \label{estimate-E-i(k)}
\bigg| \langle \frac{\xi_{i(m)}}{|\xi_{i(m)}|_g}, v_{j(k)} \rangle_g - \frac{1}{s}\Big(\hat r_{{i_0}}(q_{j(k)})-\hat r_{i(m)}(q_{j(k)}) \Big)\bigg| \leq \CE \epsone^{1/8},
\end{equation}
where $\xi_{i(m)}=\exp_{x_{{i_0}}}^{-1}(x_{i(m)})$.
By \eqref{close-xiv},
\begin{equation}
\Big| \langle \xi_{i(m)},v_{j(k)}\rangle_g - s\langle v_{j(m)},v_{j(k)}\rangle_g \Big|  \leq \Cxiv \epsone^{1/2}, \quad \forall \, k,m=1,\dots,n.
\end{equation}
Hence \eqref{ds-close-4}, \eqref{choice-s} and \eqref{estimate-E-i(k)} yield
\begin{equation} \label{estimate-i(k)}
\bigg| \langle v_{j(m)}, v_{j(k)} \rangle_g - \frac{1}{s}\Big(\hat r_{{i_0}}(q_{j(k)})-\hat r_{i(m)}(q_{j(k)}) \Big)\bigg| \leq \CEv \epsone^{1/8},\; \forall \, k,m=1,\dots,n.
\end{equation}

The formula (\ref{estimate-i(k)}) shows that we can compute the numbers 
\begin{equation} \label{def-Gkm}
G_{k,m}:=\frac 1{s} \Big(\hat r_{{i_0}}(q_{j(k)})-\hat r_{i(m)}(q_{j(k)}) \Big),\quad \textrm{where }s=\epsone^{3/8},
\end{equation}
such that
\begin{equation}\label{g products of v:s}
\Big|\bra v_{j(k)},v_{j(m)}\cet_g-G_{k,m} \Big|\leq  \CEv \epsone^{1/8},  \quad \forall \, k,m=1,\dots,n.
\end{equation}

\medskip
As above, we have considered any $n$  indices ${j(k)}$, the points
$p_{j(k)}, q_{i(k)}$ and the unit initial vectors $v_{j(k)}\in S_{x_{{i_0}}} M$ for geodesics $[x_{{i_0}} q_{j(k)}]$, $k=1,2,\dots,n$.
In view of Lemma \ref{lemma-frame-q}\footnote{Note that Lemma \ref{lemma-frame-q} is applicable here if $R/8\leq \cvdetwidth$.
However if $R/8>\cvdetwidth$, one can replace the radius $R/4$ with $R/8+\cvdetwidth$, and instead find $q_{j(k)}\in N_{\epsone}(R/8+\cvdetwidth)\cap Y$. All relevant distances would be bounded below depending on $\cvdetwidth$, which again depends on $n, \Lambda, R$.}, the formula (\ref{g products of v:s}) shows that when $\epsone$ is smaller than some uniform constant, there exist some indices $j(k)$ and points $p_{j(k)},q_{j(k)}$ such that
\begin{equation} \label{final-test-G}
\det([G_{k,m}]_{k,m=1}^n])>\frac{3}{4}\cvdetlower.
\end{equation}
Hence we can search for such indices so that \eqref{final-test-G} is satisfied by computing the determinant $\det([G_{k,m}]_{k,m=1}^n])$ for each choice of indices.


\bigskip
\noindent {\bf Reconstruction of metric.} Now let us summarize our procedure for the reconstruction of metric using the given data $\widehat{\mathcal R}_Y$ only. Let $\hat r_{{i_0}}$ be given satisfying $\hat r_{{i_0}}(y_0) > R/2$. Fix sufficiently small $\epszero\leq \epsone$ explicitly depending only on $n,\Lambda$. First, we search for a point $p_0 \in N_{\epsone}(R/8)\cap \YY$ by \eqref{search-p0}.
Due to \eqref{Da-d estimate}, see also \eqref{Da-d estimate identification},
 and the assumption that $d(x_0,y_0)<\epszero$, the point $p_0$ can be chosen using the given data (up to an error of $3\epsone$).

For $k=1,\dots,n$, we arbitrarily choose $n$ points $p_{j(k)}\in N_{\epsone}(R/8)\cap \YY$ satisfying the condition \eqref{sparse condition}.
The points $q_{j(k)}\in N_{\epsone}(R/4)\cap \YY$ are chosen according to \eqref{q-condition-i(k)}. For each $k$, choose one element $\hat r_{i(k)}\in \widehat{{\mathcal R}}_Y\cap\mathcal B_\infty(\hat r_{{i_0}},\epsone^{1/4})$ such that \eqref{criteria1-jm} and \eqref{criteria2-jm} are satisfied.
Thus we can compute the numbers $G_{k,m}$ defined by \eqref{def-Gkm}.
We test all possible choices of $n$ points $p_{j(k)}$ satisfying \eqref{sparse condition}, and find one choice such that \eqref{final-test-G} is satisfied. 
Note that by \eqref{g products of v:s}, we know that the metric corresponding to the basis that $G_{k,m}$ approximates must satisfy $\det([\langle v_{j(k)},v_{j(m)}\rangle_g]_{k,m=1}^n])>\cvdetlower/2$.

The discussion above shows that in the Riemannian normal coordinates, denoted below
by $X:B(x_{{i_0}},r)\to \R^n$, at $x_{{i_0}}$ associated with the specific basis $\{v_{j(k)}\}_{k=1}^n$ that $G_{k,m}$ approximates,
we can find the metric tensor $g_{km}(x_{{i_0}})=\bra v_{j(k)},v_{j(m)}\cet_g$ up to a uniformly bounded error, that is,
we can find numbers ${\hat g}_{km}:=G_{k,m}$ such that
\begin{equation} \label{close-ghat}
\Big|g_{km}(x_{{i_0}})-{\hat g}_{km}\Big|\leq  \CEv \epsone^{1/8},  \quad \forall \, k,m=1,2,\dots,n.
\end{equation}
In particular, since the basis $\{v_{j(k)}\}_{k=1}^n$ are unit, we have $|g_{km}|\leq 1, \ |\hat{g}_{km}|\leq 2$ for all $k,m$.

\medskip
\noindent {\bf Reconstruction of normal coordinates.} 
Suppose we have already picked indices $j(k)$ and points $p_{j(k)},q_{j(k)}$ such that the metric is approximated as above.
We can also find the coordinates of the points corresponding to elements $\hat r_\ell$ in a neighborhood of $\hat r_{{i_0}}$ 
in the a normal coordinate $X:B(x_{{i_0}},r)\to \R^n$, up to a uniformly bounded error.
Indeed, for any $\hat r_\ell$ satisfying \eqref{c2neighbor}, we can compute  the numbers
\begin{equation} \label{def-Xk}
\hat X_k(x_\ell) := \hat r_{{i_0}}(q_{j(k)})-\hat r_\ell(q_{j(k)}),\quad k=1,2,\dots,n.
\end{equation}
Note that due to Proposition \ref{lem: properties}, we can choose sufficiently small $\crneighbor$ such that $d(x_{\ell},x_{{i_0}})<\textrm{inj}(M)$.
Then by (\ref{eq-xil}), we have
\begin{eqnarray*} \label{eq-xil-jk}
\Big| X_k(x_{\ell}) - \hat X_k (x_{\ell}) \Big| &\leq& \betadissecond |x_{{i_0} } x_{\ell}| \Big(\epsone^{1/4}+|x_{{i_0}} x_{\ell}|^{1/3} \Big) +2\epsone \nonumber \\
&\leq& 3\betadissecond \Big( |x_{{i_0}} x_{\ell}|^{4/3}+\epsone \Big),
\end{eqnarray*}
where 
\begin{equation} \label{def-Xk-true}
X_k(x_{\ell}):=\langle \xi_{\ell}, v_{j(k)} \rangle_g=\langle \exp_{x_{{i_0}}}^{-1}(x_{\ell}), v_{j(k)} \rangle_g
\end{equation}
are the true values of the coordinates of the point $x_\ell$ in the $X$-coordinates. Here we have used Young's inequality $ab\leq a^{4/3}+b^4$ in the last inequality. This concludes the proof for Case 1.

\begin{figure}[htbp]
\begin{center}
\includegraphics[width=7cm]{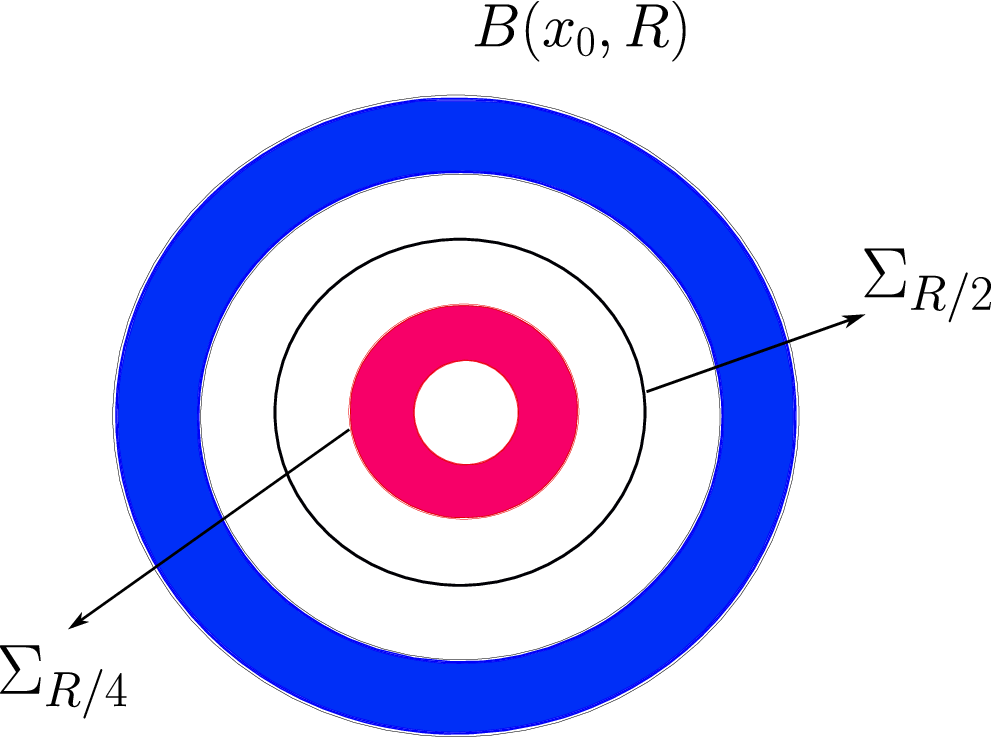} \label{figure-circles}
\end{center}
\caption{\it When $x_{{i_0}}\in B(x_0,R/2)$, we choose measurement points $p,q$ from the outer (blue) band such that the distances between $x_{{i_0}},p,q$ are bounded away from zero. When $x_{{i_0}}\notin B(x_0,R/2)$, we choose measurement points in the inner (red) band.}
\end{figure}

\bigskip
$\bullet$ \textbf{Case 2:} $\hat r_{{i_0}}(y_0) \leq R/2$.

\medskip
In this case, $d(x_{{i_0}},x_0)\leq R/2+2\epsone$. To keep distances bounded away from zero, one can choose points $q_{j(k)},p_{j(k)}$ from the outer layer $B(x_0,R)\setminus B(x_0,3R/4)$, see Figure $7$.
More precisely, we first choose arbitrary $n$ points $\{q_{j(k)}\}$ in $N_{\epsone}(3R/4)\cap \YY$ satisfying \eqref{sparse condition}. 
Note that Lemma \ref{lemma-frame} still holds in this case due to Lemma \ref{lemma-angle-separation-inner}, see Remark \ref{remark-inner}. 
For each ${q_{j(k)}}$, we can choose one point $p_{j(k)}\in N_{\epsone}(R)\cap \YY$ such that \eqref{q-condition-i(k)} is valid.
Observe that the set of points satisfying \eqref{q-condition-i(k)} (for each $q_{j(k)}$) is nonempty. This is because one can extend the minimizing geodesic $[x_{{i_0}} q_{j(k)}]$ further until it intersects with $\Sigma_R$, since we are within the injectivity radius. Thus the point in $N_{\epsone}(R)\cap \YY$ within $\epszero$-distance from the intersection point satisfies \eqref{q-condition-i(k)}.
In particular, the bounds in \eqref{pqx-bounded} still hold. From this point, the exact proof of Case 1 works in this case.

\medskip
Now we prove Corollary \ref{coro 1}.
\medskip

\noindent {\bf Proof of Corollary \ref{coro 1}.} 
Let us fix the basis $\{v_k\}_{k=1}^n$ for which the metric has been approximated in Theorem \ref{main 1}.
Observe that
$$d(x_{\ell},x_{{i_0}})=\Big| \exp_{x_{{i_0}}}^{-1}(x_{\ell}) \Big|_g= \bigg(\sum_{j,k=1}^n g^{jk} X_j X_k \bigg)^{1/2} ,$$ 
where $X_k$ is defined in \eqref{correct coordinates} and $(g^{jk})$ is the matrix inverse of $(g_{jk})$.
Hence we define an approximate distance by
\begin{equation}
\widehat{d}_{\ell,{i_0}}=\bigg(\sum_{j,k=1}^n {\hat g}^{jk} {\hat X}_j {\hat X}_k \bigg)^{1/2}.
\end{equation}
We use the following notation for convenience:
\begin{equation}
|\xi|_g^2=\sum_{j,k=1}^n g^{jk} X_j X_k,\quad |\xi|_{\hat g}^2=\sum_{j,k=1}^n {\hat g}^{jk} X_j X_k, \quad \xi=(X_1,\dots,X_n).
\end{equation}
Denote ${\hat\xi}=({\hat X}_1,\dots, {\hat X}_n)$.
Then
\begin{eqnarray*}
\Big| \widehat{d}_{\ell,{i_0}}- d(x_{\ell},x_{{i_0}}) \Big| &\leq& \Big| |\hat{\xi}|_{\hat g}-|\xi|_{\hat g}\Big|+\Big| |\xi|_{\hat g}-|\xi|_g\Big| \\
&\leq& \Big| \hat{\xi}-\xi \Big|_{\hat g} +|\xi|_g^{-1} \Big| |\xi|_{\hat g}^2-|\xi|_g^2 \Big|.
\end{eqnarray*}

Due to \eqref{metric error} and $\det([\bra v_j,v_k\cet_g]_{j,k=1}^n)\geq \cvdetlower$, all components $g^{jk},\hat{g}^{jk}$ are bounded above by $C(n,\cvdetlower)$. Hence by \eqref{appr. coordinates},
\begin{equation} \label{close-d-1}
\Big| \hat{\xi}-\xi \Big|_{\hat g} \leq C(n,\cvdetlower) \Cmainone \Big( d(x_{{i_0}}, x_{\ell})^{4/3}+\epsone \Big).
\end{equation}
Moreover, the largest eigenvalue of the matrix $(g_{jk})$ is bounded above by $C(n)$, and thus the eigenvalues of $(g^{jk})$ are bounded below by $C(n)^{-1}$. Hence by \eqref{metric error},
\begin{eqnarray} \label{close-d-2}
|\xi|_g^{-1} \Big| |\xi|_{\hat g}^2-|\xi|_g^2 \Big| &\leq& C(n) |X|_{\R^n}^{-1} \sum_{j,k=1}^n \Big(\hat{g}_{jk}-g_{jk} \Big) X_j X_k 
\leq n C(n) \Cmainone \epsone^{\frac18} |\xi|_g.
\end{eqnarray}
Thus by \eqref{close-d-1} and \eqref{close-d-2}, we obtain
\begin{eqnarray*}
\Big| \widehat{d}_{\ell,{i_0}}- d(x_{\ell},x_{{i_0}}) \Big| &\leq& C(n,\cvdetlower) \Cmainone \Big( d(x_{{i_0}}, x_{\ell})^{4/3}+\epsone \Big) + n C(n) \Cmainone \epsone^{\frac18} d(x_{{i_0}}, x_{\ell}) \\
&\leq& \Cmaintwo \Big(d(x_{{i_0}}, x_{\ell})^{4/3}+\epsone^{1/2} \Big).
\end{eqnarray*}

\vspace{-8mm}\hfill$\square$



\section{Global constructions}
\label{sec-application}

In this section we prove Theorem \ref{main 1 no boundary}.
\medskip

\noindent {\bf Proof of Theorem \ref{main 1 no boundary}.} 
(1) The first claim has been proved in Proposition \ref{lem: properties}(2), and we prove the second claim here. 
Let $\rho>0$ be a parameter which is determined later.
Given two indices $i,j\in \{1,\dots, I\}$, we consider the following minimization problem
\begin{equation} \label{def-minimization}
\hat{d}_{i,j}:=\min \bigg( \hat{d}_{i,\sigma(1)} +\sum_{k=1}^{N-1} \hat{d}_{\sigma(k), \sigma(k+1)} +\hat{d}_{\sigma(N),j} \bigg)
\end{equation}
over all chain of $N$ indices $\sigma(k)\in \{1,\dots,I\}$ with $N\leq 1+\Lambda/\rho$, and with 
\begin{equation}\label{def-minimization-condition}
\|\hat r_{\sigma(k)}-\hat r_{\sigma(k+1)}\|_{\ell^\infty(Y)}\leq\rho+\Cmainthree\epsone^{1/2}+2\epsone<\crneighbor,\quad\textrm{for }k=0,\dots,N,
\end{equation}
where we denote $\sigma(0):=i$ and $\sigma(N+1):=j$. In this minimization problem, the numbers $\hat{d}_{\sigma(k), \sigma(k+1)}$ are determined by the data $\widehat{{\mathcal R}}_Y\cap\mathcal B_\infty(\hat r_{\sigma(k)},\epsone^{1/4})$ in Corollary \ref{coro 1} since the condition \eqref{def-minimization-condition} is valid for all $k$. Thus the solution of the minimization problem can be found using the given data $\widehat{{\mathcal R}}_Y$ only. Note that the numbers $\hat d^{_Y}_{j,j'}$ assumed in Corollary \ref{coro 1} can be determined by $\widehat{{\mathcal R}}_Y$ due to \eqref{Da-d estimate}.

\smallskip
Now let us analyze the property of the solution of the minimization problem \eqref{def-minimization}.
Given any pair of points in ${X}$, consider a shortest path $\gamma$ connecting these two points. 
On this shortest path $\gamma$, one can choose a chain of points with at most $N\leq 1+\Lambda/\rho$ points such that each pair of adjacent points has distance at most $\rho$. 
Since the set ${X}$ is an $\epsX$-net of $M$, we can replace this chain of points on $\gamma$ by points in ${X}$, and thus each pair of adjacent points has distance at most $\rho+\epsX=\rho+\Cmainthree\epsone^{1/2}$.
We require $\rho+\Cmainthree\epsone^{1/2}<\crneighbor/2$ such that the condition \eqref{crneighbor-Thm1.2} is satisfied by the elements in $\widehat{{\mathcal R}}_Y$ corresponding to this chain of points in ${X}$, by virtue of \eqref{d-to-r}.
Let us relabel this chain of points by $x_1,\dots,x_{N}$ with endpoints $x_1,x_N$ for convenience.
Then the estimate (\ref{final coro estimate}) shows that for this particular chain of points, we have
\ba 
E:=\bigg| d(x_1,x_N) -\sum_{i=k}^{N-1} \hat{d}_{k,k+1} \bigg|&\leq& N\Cmaintwo \Big(\epsone^{1/2}+(\rho+\epsone^{1/2})^{4/3} \Big)
\\
\\
&< &8\Lambda\Cmaintwo \Big(\epsone^{1/2}\rho^{-1}+\rho^{1/3} \Big).
\ea
We obtain a good estimate when we choose $\rho=\epsone^{3/8}$, namely
\beq\label{D-distance estimate}
E < 16 \Lambda \Cmaintwo \epsone^{1/8}.
\eeq

This particular chain of points above satisfies the conditions of the minimization problem \eqref{def-minimization}, and thus the solution of the minimization problem also satisfies \eqref{D-distance estimate} when we choose $\rho=\epsone^{3/8}$.
Therefore, the solution of the minimization problem \eqref{def-minimization} gives us an approximate distance $\hat{d}(x_i,x_j):=\hat{d}_{i,j}$ that satisfies
\beq\label{D-distance estimate 0}
 \Big|\hat{d}(x_i,x_j)-d(x_i,x_j) \Big|<16 \Lambda \Cmaintwo \epsone^{1/8}, \quad \forall \, i,j\in \{1,\dots,I\}.
 \eeq
This proves part (1) of Theorem \ref{main 1 no boundary}.

\smallskip
Part (2) is a direct consequence of the estimate \eqref{D-distance estimate 0} and \cite[Corollary 1.10]{FIKLN}.
 \hfill $\square$\medskip

\section{Reconstruction of the manifold from the noisy heat kernel}\label{Sec. heat kernel}

In this section, we consider the reconstruction of a manifold from noisy heat kernel measurements \eqref{kernel-measurements} satisfying \eqref{eta-error}.

\begin{theorem} \label{Thm-kernel-reconstruction}
Let $M$ be a closed Riemannian manifold of dimension $n$ satisfying the bounds \eqref{basic 1} with parameter $\Lambda$, and let $U=B(y_0,R)$ be a ball of radius $R>\Lambda^{-1}$.
 Suppose the Ricci curvature of $M$ is non-negative. Then there exist constants $\hat{\sigma}, \Cmainfivelast>0$ explicitly depending only on $n,\Lambda$, such that the following holds for $0<\sigma < \hat{\sigma}$ and $0<h\leq \sigma^{1/2}$. 

Let $Y=\{y_j:\  j=0,1,\dots,J\}$ be an $h$-net in the ball $U$. 
Assume that either
\smallskip

\noindent
(i)   $\{z_i:\ i=1,\dots, I\}$ is an $h$-net in $M\setminus \overline U$,
and we are given $\hat d^{_Y}_{j,j'}$, $j,j'=0,1,\dots,J$ such that $|\hat d^{_Y}_{j,j'}-d(y_j,y_{j'})|<h$,

\noindent
 or
\smallskip
 
\noindent
(ii) $\{z_i:\ i=1,\dots, I\}$ is an $h$-net in $M$. 
\smallskip

\noindent Moreover, assume that we are given the data 
\begin{equation} \label{kernel-data}
\Big \{\tilde{G}(y_j,z_i,t):\ i=1,\dots, I, \; j=0,1,\dots,J, \; 0<t<1 \Big\}
\end{equation}
which satisfy
\begin{equation} \label{heat-kernel-recon}
e^{-\frac{\sigma}{t}} \leq \frac{\tilde{G}(y_j,z_i,t)}{G(y_j,z_i,t)}\leq e^{\frac{\sigma}{t}}, \;\textrm{ for all }i=1,\dots,I,
\ j=0,\dots, J,\ 0<t<1.
\end{equation}
Then the given data \eqref{kernel-data} determine a smooth Riemannian manifold $(\hat M,\hat g)$ that is
diffeomorphic to $M$. Moreover, there is a diffeomorphism
 $F: \hat M\to M$ such that 
 \beq
\label{Lip-conditionB}
\frac 1L\leq \frac{d_{M}(F(x),F(x'))}{d_{\hat M}(x,x')}\leq L,\quad \hbox{for }x,x'\in \hat M,
\eeq
where  $L=1+\Cmainfivelast \sigma^{1/24}$.
\end{theorem}


\smallskip 
\noindent {\bf Proof.} 
By Corollary 3.1 and Theorem 4.1 in \cite{LY} (see also \cite{SY}), for every $\epsilon\in (0,1)$, $t>0$,
\begin{equation} \label{upperbound-kernel}
G(y,z,t) \leq C_{\epsilon} v_{z,t} \exp\bigg( -\frac{d^2(y,z)}{(4+\epsilon)t} \bigg),
\end{equation}
and 
\begin{equation} \label{lowerbound-kernel}
G(y,z,t) \geq C_{\epsilon}^{-1} v_{z,t} \exp\bigg( -\frac{d^2(y,z)}{(4-\epsilon)t} \bigg),
\end{equation}
where $v_{z,t}={\rm vol}^{-1}(B(z,\sqrt{t}))$, and $C_{\epsilon}\to \infty$ as $\epsilon\to 0$.

Taking $\log$ on both sides of \eqref{upperbound-kernel}, we have
$$\frac{d^2}{4}+t\log G \leq t\log (C_{\epsilon} v_{z,t})+\frac{d^2}{4}-\frac{d^2}{4+\epsilon} \leq t\log (C_{\epsilon} v_{z,t}) + C(\Lambda) \epsilon.$$
Similarly, from \eqref{lowerbound-kernel},
$$\frac{d^2}{4}+t\log G \geq t\log (C_{\epsilon}^{-1} v_{z,t})+\frac{d^2}{4}-\frac{d^2}{4-\epsilon} \geq t\log (C_{\epsilon}^{-1} v_{z,t}) - C(\Lambda) \epsilon.$$
Combining the two inequalities above, we obtain
\begin{equation} \label{kernel-d}
\Big|d^2(y,z)+4t \log G(y,z,t) \Big| \leq 4t |\log C_{\epsilon} | + 4t |\log v_{z,t} | + C(\Lambda) \epsilon.
\end{equation}
From the noisy measurements $\tilde G(y,z,t)=\eta(y,z,t) G(y,z,t)$, \eqref{kernel-d} and \eqref{eta-error} yield that
\begin{eqnarray}
\Big|d^2(y,z)+4t \log \tilde{G}(y,z,t) \Big| &\leq& \Big|d^2(y,z)+4t \log G(y,z,t) \Big| +4t |\log \eta | \nonumber \\
&\leq& 4\sigma+4t|\log C_{\epsilon} | + 4t |\log v_{z,t} | + C(\Lambda) \epsilon.
\end{eqnarray}

Now we pick sufficiently small $\epsilon>0$ such that $C(\Lambda) \epsilon<\sigma$. For small $t>0$, we know
$$t|\log v_{z,t}| = t \big|\log {\rm vol}(B(z,\sqrt{t})) \big| \leq C(n)t|\log t|.$$
Thus one can pick sufficiently small $t>0$ such that $4t|\log C_{\epsilon} | + 4t |\log v_{z,t} |<\sigma$. Hence,
$$\Big|d^2(y,z)+4t \log \tilde{G}(y,z,t) \Big| < 6\sigma,$$
which yields that 
\begin{equation} \label{kernel-d-appro}
\bigg|d(y,z)-\sqrt{4t \big|\log \tilde{G}(y,z,t) \big|} \bigg| < 6\sigma^{\frac12}.
\end{equation}

First, consider the case when the condition (i) is valid.
Since $\{z_i\}$ is an $h$-net in $M\setminus \overline U$ for $h\leq \sigma^{1/2}$, \eqref{kernel-d-appro} shows that the data $\tilde{G}(y_j,z_i,t)$ for some suitable choice of $t$ (depending on $\sigma$) give a $7\sigma^{1/2}$-approximation of the distances of the pairs  $(y_j,z_i)$. Moreover, we are already given 
an $h$-approximation $\hat d^{_Y}_{j,{j'}}$ of the distances of the pairs  $(y_j,y_{j'})$ in $Y$. 
Since the set
$X=\{z_i:\ i=1,\dots,I\}\cup Y$ is a $(2h)$-net in $M$, 
thus the sets $X,Y$ and the given data
satisfy the conditions (a1) and (a2) with parameter $\epsone=7\sigma^{1/2}$.

Second,  consider the case when the condition (ii) is valid.
Since $\{z_i\}$ is an $h$-net in $M$ for $h\leq \sigma^{1/2}$, \eqref{kernel-d-appro} shows that the data $\tilde{G}(y_j,z_i,t)$ for some suitable choice of $t$ (depending on $\sigma$) give a $7\sigma^{1/2}$-approximation of the interior distance functions on $Y$. 
Thus the sets $\{z_i\},Y$ and the given data
satisfy the conditions (a1) and (a2) with parameter $\epsone=7\sigma^{1/2}$.

These considerations imply that in both cases (i) and (ii) the claim follows by applying Theorem \ref{main 1 no boundary}(2).
 \hfill $\square$\medskip

\smallskip
Finally, we obtain the uniqueness and the stability for the inverse problem for heat kernel.

\medskip 
\noindent {\bf Proof of Theorem \ref{Thm-kernel}.} 
Let $X_1=\{z^1_1,\dots,z^1_I\}\cup \{y^1_0,\dots,y^1_J\}\subset M_1$,
$Y_1= \{y^1_0,\dots,y^1_J\}\subset M_1$
 and 
$X_2=\{z^2_1,\dots,z^2_I\}\cup \{y^2_0,\dots,y^2_J\}\subset M_2$, $Y_2= \{y^2_0,\dots,y^2_J\}\subset M_2$.

Due to the condition \eqref{kernel-closeness}, the heat kernel data $G_2(y_j^2,z_i^2,t)$ of $M_2$ can be used as the noisy observations of the heat kernel of $M_1$ at $(y_j^1,z_i^1,t)$.
Then by the proof of Theorem \ref{Thm-kernel-reconstruction}, the heat kernel data $G_2(y_j^2,z_i^2,t)$ of $M_2$ and $d_{M_2}(y_j^2,y_{j'}^2)$  determine 
the distances
$d_{M_1}(x,y)$ of $M_1$  for $(x,y)\in X_1\times Y_1$  up to an error $7\sigma^{1/2}$.
Similarly, the heat kernel data $G_2(y_j^2,z_i^2,t)$ of $M_2$ and $d_{M_2}(y_j^2,y_{j'}^2)$
also determine the 
distances
$d_{M_2}(x,y)$ of $M_2$ for $(x,y)\in X_2\times Y_2$  up to an error $7\sigma^{1/2}$.
Thus, by enumerating the points in $X_l$  as $\{x^l_i:\  i=0,1,\dots,I'\}$ and the points in $Y_l$  as $\{y^l_j:\ j=0,1,\dots ,J\}$,
we see that
\ba
|d_{M_1}(x^1_i,y^1_j)-d_{M_2}(x^2_i,y^2_j)|<14\sigma^{\frac12},\quad i=0,1,\dots,I',\ j=0,1,\dots J.
\ea
Thus, the numbers $d_{i,j}=d_{M_1}(x^1_i,y^1_j)$ 
can be used
as the noisy distance data both for the manifold $M_1$ and for the manifold $M_2$ with an error $14\sigma^{\frac12}$.
Hence, by Theorem \ref{main 1 no boundary} we have 
$$
|d_{M_1}(x^1_i,x^1_{i'})-d_{M_2}(x^2_i,x^2_{i'})|\leq \Cmainfour  (14\sigma^{\frac12})^{\frac18}<2 \Cmainfour \sigma^{\frac{1}{16}}, \quad i,i'=0,1,\dots,I',
$$
and there is a manifold $\hat M$   such that there are $L_1$-bi-Lipschitz diffeomorphisms  $F_1:M_1\to \hat M$
and $F_2:M_2\to \hat M$ with $L_1$ given in Theorem\ \ref{main 1 no boundary}. Hence there is an $(L_1^2)$-bi-Lipschitz diffeomorphism  $F=F_2^{-1}\circ F_1:M_1\to M_2$, which proves the claim.
 \hfill $\square$\medskip

\medskip 
\noindent {\bf Proof of Corollary
\ref{cor: Thm-kernel}.} Let $\Lambda>R^{-1}$ be such that both $M_1$ and $M_2$ satisfy the
geometric bounds  bounds \eqref{basic 1}.
Consider an arbitrary $\sigma>0$, $h=\sigma^{1/2}$. Let $\{y^1_j: j=0,1,\dots, J\}$ be an $h$-net in $U_1$ and $\{y^2_j=\Phi(y^1_j): j=0,1,\dots, J\}$ be an $h$-net in $U_2$. 

We recall that the map $\Psi:M_1\setminus U_1\to M_2\setminus U_2$  is assumed only to be 
a bijection, and thus we need to do some additional considerations to obtain suitable
$h$-nets on sets $M_1\setminus U_1$ and $M_2\setminus U_2$ .
To that end, let  $\{\tilde  z^1_i,$ $i=1,2,\dots,I_1\}$ be an $h$-net in $M_1\setminus U_1$ and
  $\{\tilde  z^2_i$, $i=1,2,\dots,I_2\}$ be an $h$-net in $M_2\setminus U_2$.
  Then we define   
  \ba
  z^1_i=\begin{cases} \tilde  z^1_i, &\hbox{for $i=1,2,\dots,I_1,$}\\ 
   \Psi^{-1}(\tilde  z^2_{i-I_1}), &\hbox{for $i=I_1+1,I_1+2,\dots,I_1+I_2,$}
   \end{cases}
   \ea
and  
  \ba
  z^2_i=\begin{cases} \Psi(\tilde  z^1_i), &\hbox{for $i=1,2,\dots,I_1,$}\\ 
  \tilde  z^2_{i-I_1}, &\hbox{for $i=I_1+1,I_1+2,\dots,I_1+I_2$}.
   \end{cases}
   \ea
  Then $\{  z^1_i:\ i=1,\dots,I_1+I_2\}$  is an $h$-net in $M_1\setminus U_1$ and 
$\{  z^2_i:\ i=1,\dots,I_1+I_2\}$  is an $h$-net in $M_2 \setminus U_2$, and we have
$$
G_1(y_j^1,z_i^1,t)=G_2(y_j^2,z_i^2,t),\quad  i=1,\dots,I_1+I_2,\ j=0,\dots,J,\ 0<t<1.
$$
%
%
 By applying Theorem \ref{Thm-kernel}, it follows that the Gromov-Hausdorff
distance of the metric spaces $(M_1,d_{g_1})$ and $(M_2,d_{g_2})$ is smaller than
$\Clast\sigma^{1/24}$, where $\Clast>0$ depends only on $n$ and $\Lambda$,
see e.g. Corollary 7.3.28 in \cite{Burago}. 
Letting $\sigma\to 0$, we see that Gromov-Hausdorff
distance of $(M_1,d_{g_1})$ and $(M_2,d_{g_2})$ is zero, which
 implies that $(M_1,d_{g_1})$ and $(M_2,d_{g_2})$
are isometric as (compact) metric spaces, see e.g. \cite{Burago,Pe}.
 By the Myers-Steenrod theorem,  
there is a  diffeomorphism $F:M_1\to M_2$ between Riemannian manifolds such that $g_2=F_*g_1$. This proves the claim.
 \hfill $\square$\medskip

\bigskip
\appendix

\section{Auxiliary Lemmas}

\begin{lemma}\label{holderbound} 
Let $M$ be a closed Riemannian manifold with sectional curvature bounded below by ${\rm Sec}_M\geq -\Lambda^2$. Suppose $\gamma_{x,v_1}(t),\ \gamma_{x,v_2}(t)$ are two distance-minimizing geodesics emanating from $x\in M$ with unit initial vectors $v_1,v_2\in S_x M$. Denote by $\alpha$ the angle between $v_1$ and $v_2$.
Then there is a uniform constant $\CexpLip>1$, explicitly depending only on $\Lambda$, such that
\begin{equation} 
d\big(\gamma_{x,v_1}(t_1),\ \gamma_{x,v_2}(t_2)\big)\leq |t_1-t_2|+\CexpLip \alpha ,\quad \forall \, t_1,t_2 \in[0,\Lambda].
\end{equation}
\end{lemma}

\noindent {\it Proof.} 
Assume that $t_1\leq t_2$. Let us denote $a=\gamma_{x,v_1}(t_1)$ and $b=\gamma_{x,v_2}(t_1)$.
We can compare the triangle $axb$ with a triangle $\overline{a}\overline{x}\overline{b}$ in the rescaled hyperbolic plane $H$ with constant sectional curvature $-\Lambda^2$, satisfying that $d(x,a)=\overline{d}(\overline{x}, \overline{a})=t_1$, $d(x,b)=\overline{d}(\overline{x} ,\overline{b})=t_1$,
$\alpha=\angle \overline{a}\overline{x}\overline{b}$.
Then Toponogov's theorem yields $d(a,b)\leq \overline{d}(\overline{a},\overline{b})$.

On the hyperbolic plane $H$, the exponential map is smooth everywhere, and its differential is uniformly bounded. Hence,
\begin{equation} \label{eq-exp-Lip}
 \overline{d}(\overline{a},\overline{b}) \leq
C(\Lambda) |t_1v_1-t_1v_2| \leq C(\Lambda) t_1 \alpha.
\end{equation}
Then,
\begin{eqnarray*}
d\big(\gamma_{x,v_1}(t_1),\ \gamma_{x,v_2}(t_2)\big) &\leq&
d(a,b)+d\big(\gamma_{x,v_2}(t_1),\ \gamma_{x,v_2}(t_2)\big) \\
&\leq& \overline{d}(\overline{a},\overline{b}) +|t_2-t_1| \leq |t_2-t_1|+C(\Lambda)\Lambda \alpha.
\end{eqnarray*}
\hfill $\square$\medskip

\begin{lemma}\label{lem: Topog estimates general appen}
There exists a
uniform constant $\CfirstToponogov>1$ such that the following holds. Let $N$ be a compact Riemannian manifold with boundary $\partial N$ with sectional curvature bounded below by $Sec_N\geq -\Lambda^2$.
Let $a,b,c\in N$ and $\beta$ be the angle of the length minimizing
curves $[ab]$ and $[bc]$ at $b$. Then we have
\beq
\label{Si far points new appen}
  |ac| \leq |ab|-|bc|\cos\beta + \CfirstToponogov |bc|^2/\min\{\Lambda^{-1},|ab|,d(b,\p N)\}.
\eeq
\end{lemma}

\noindent
{\bf Proof.} 
To prove the statement, we apply Toponogov's Theorem (e.g. \cite[Thm. 79]{Pe}) to the triangle $abc$.
Below, let $H$ be the rescaled hyperbolic plane of constant sectional curvature $-\Lambda^2$,
and $\bar a,\bar b,\bar c\in H$ be such that $|\bar a\bar b|=|ab|$,
$|\bar b\bar c|=|bc|$ and $\angle \bar a\bar b\bar c=\angle abc=\beta$. Here by $|xy|$, we denote the distance between points $x$ and $y$ in whatever space they belong.

(i) Let us first consider the case when 
\beq\label{assumption 1}
|ab|=\frac 14,\quad 
|bc|<\frac 14,\quad  \hbox{and}\quad d(b,\p N)\geq 1.
\eeq
 Then 
$|ac|\leq |ab|+|bc|\leq \frac 12$. 
Applying Toponogov's theorem to triangle $abc$ implies that $|ac|\leq |\bar a\bar c|$.

By \cite[Prop.\ 48]{Pe}, the law of cosines on $H$ gives
\beq\label{the law of cosines}\\
\nonumber
\cosh({\Lambda}{  |\overline a\overline c|})&=& \cosh ({\Lambda}{ |\overline a\overline b|})\,\cosh ({\Lambda}{|\overline b\overline c|})-
 \sinh ({\Lambda}{ |\overline a\overline b|})\,\sinh({\Lambda}{ |\overline b\overline  c|})\,\cos\beta
\eeq
and thus 
\ba
\cosh({\Lambda}{|a c|})&\leq & \cosh ({\Lambda}{ | a b|})\,\cosh ({\Lambda}{| b c|})-
 \sinh ({\Lambda}{ | a b|})\,\sinh( {{\Lambda} | b  c|})\,\cos\beta.
\ea
Using Taylor series, the above yields
\ba& &
\cosh( {\Lambda}{|ac|})-\cosh ( {\Lambda}{|ab|})
 \\
 &\leq &
\cosh ( {\Lambda}{ |ab|})\,(\cosh ( {\Lambda}{ |bc|})-1)- \sinh ( {\Lambda}{ |ab|})
\,\sinh( {\Lambda}{ |bc|})\,\cos\beta
\\
&=&
- V |bc| \,\cos\beta 
+E_1,
\ea
where $V=\Lambda \sinh({\Lambda}{|ab|})>0$ and  $E_1$ satisfies $|E_1|\leq C|bc|^2$,
where $C$ is a uniform constant.

By triangle inequality, $- |bc|\leq |ac|-|ab|\leq |bc|$.
Then one can show that
 \ba
 \cosh( {\Lambda}{|ac|})-\cosh ( {\Lambda}{|ab|})
  \geq V ( |ac|-|ab|)+E_2,
\ea
where $E_2$ satisfies $|E_2|\leq C(|ac|-|ab|)^2\leq C|bc|^2$, where $C$ is a uniform constant.
Combining these we see that
\ba
V ( |ac|-|ab|)+E_2\leq 
- V |bc| \,\cos\beta
+E_1,
\ea
or
\ba
 |ac|-|ab|\leq 
-  |bc| \,\cos\beta
+V^{-1}(E_1-E_2),
\ea
which yields the  inequality 
\beq
\label{Si far points A}
|ac| \le |ab|-|bc|\cos\beta + \CfirstToponogov |bc|^2,
\eeq
where $\CfirstToponogov $ is a uniform constant. We
can assume that $\CfirstToponogov>8$.

Consider next the case when 
\beq\label{assumption 2}
|ab|=\frac 14,\quad 
  \hbox{and}\quad d(b,\p N)\geq 1.
\eeq
If it holds that $|bc|\geq \frac 1 4$, then $\CfirstToponogov>8$ yields
that $\CfirstToponogov |bc|^2>2|bc|$. Considering $|ac|\leq |ab|+|bc|$, we see that  
(\ref{Si far points A}) automatically holds. Since we have already proven 
(\ref{Si far points A}) when $|bc|<\frac 1 4$, we
can conclude that (\ref{Si far points A}) holds under 
assumptions (\ref{assumption 2}).

Next, consider the case when
\beq\label{assumption 3}
|ab|\geq \frac 14,\quad  d(b,\p N)\geq 1.
\eeq

Let $a^{\prime}$ be the point on $[ab]$ with $|a^{\prime}b|=\frac 14$. 
By the triangle inequality we have $|ac|\le|aa^{\prime}|+|a^{\prime}c|$.
Moreover, we have $|aa^{\prime}|=|ab|-|a^{\prime}b|$   and thus (\ref{Si far points A})
for the triangle $a^{\prime}bc$ implies
$$
 |ac|-|ab| \le |aa^{\prime}|+|a^{\prime}c|-|ab|\leq |a^{\prime}c|-|a^{\prime}b| \le -|bc|\cos\beta +  \CfirstToponogov |bc|^2.
$$
Hence, (\ref{Si far points A}) holds under 
assumption (\ref{assumption 3}).

The
 inequality (\ref{Si far points A}) yields that
 the
 inequality   (\ref{Si far points new appen}) holds.
%
%
%
%
%
%
%
%
%
%
 \hfill $\square$\medskip

\begin{lemma} \label{lemma-angle-separation-inner}
Let $M$ be a closed Riemannian manifold with sectional curvature bounded by $ |{\rm Sec}_M| \leq \Lambda^2$. Let $B(x_0,R)$ be an open ball for $R\leq \min\{ {\rm inj}(M)/2, \pi/(4\Lambda)\}$.
Then there exist uniform constants $\hat{\varepsilon},\cvdetneighbor,\Cangleseparation>0$ explicitly depending on $\Lambda,R$, such that the following holds for $0<\varepsilon<\hat{\varepsilon}$.

Given a point $x\in B(x_0,R/2)$, take $z_0$ to be the nearest point on $\partial B(x_0,R)$ from $x$.
Let $z_1, z_2\in N_{\varepsilon}(R)$ such that 
\begin{equation}\label{condition-separation-inner}
|z_1 z_2|\geq \Cnetseparation \varepsilon,\quad |z_0 z_1|<\cvdetneighbor,\quad |z_0 z_2|<\cvdetneighbor,
\end{equation}
for some $\Cnetseparation \geq 32(\Lambda R)^{-1}$.
Then $\angle z_1 x z_2 > \Cangleseparation \Cnetseparation \varepsilon$.
\end{lemma}

\noindent {\it Proof.} 
The proof is similar to Lemma \ref{lemma-angle-separation}(1): we use the upper bound for the angle $\angle x z_1 z_2$ to derive an lower bound for $\angle z_1 x z_2$. First, we show that the incident angle of $[x z_1]$, i.e. the angle of $[x z_1]$ with the tangent space $T_{z_1}\Sigma_{|x_0 z_1|}$ is bounded from below. 
Suppose $|x_0 z_1|\leq |x_0 z_2|$. 

We take the point $z'_1\in \Sigma_{3R/2}$ such that $|x_0 z'_1|-|x_0 z_1|=|z_1 z'_1|$. Let $z'_0$ be the nearest point in $\Sigma_{3R/2}$ from $x$. Then for sufficiently small $\cvdetneighbor,\varepsilon$, we have
$$|x z'_1|\geq |x z'_0|-|z'_0 z'_1| \geq |x z_0|+\frac{R}{2}-2\cvdetneighbor \geq |x z_1|+|z_1 z'_1|-3\cvdetneighbor.$$
Hence the same argument as \eqref{estimate-eta-section5} gives
$(\pi-\angle x z_1 z'_1)^2\leq \CsecondToponogov^{-1} 3\cvdetneighbor$,
which is
\begin{equation}\label{estimate-eta-inner}
\angle x z_1 x_0 \leq C(\Lambda) \cvdetneighbor^{1/2}.
\end{equation}

For the upper bound for the angle $\angle x_0 z_1 z_2$, one can apply Rauch comparison theorem to compare with the sphere of constant sectional curvature $\Lambda^2$. Namely, we take the triangle $\overline{x_0}\overline{z_1}\overline{z_2}$ on the sphere such that $|x_0 z_1|=|\overline{x_0}\overline{z_1}|$, $|z_1 z_2|=|\overline{z_1}\overline{z_2}|$, $\angle x_0 z_1 z_2=\angle \overline{x_0}\overline{z_1}\overline{z_2}$.
Due to Rauch comparison theorem, $|x_0 z_2|\geq |\overline{x_0} \overline{z_2}|$. Hence,
\begin{eqnarray*}
\cos({\Lambda}{ |x_0 z_2|})\leq \cos ({\Lambda}{ |x_0 z_1|})\,\cos ({\Lambda}{| z_1 z_2|})+
 \sin ({\Lambda}{ |x_0 z_1|})\,\sin({\Lambda}{ | z_1 z_2|})\,\cos(\angle x_0 z_1 z_2).
\end{eqnarray*}
Assume $\angle x_0 z_1 z_2>\pi/2$. Since $|x_0 z_1|\leq R+\varepsilon <\pi/(2\Lambda)$, then
\begin{eqnarray*}
\cos({\Lambda}{ |x_0 z_2|})-\cos ({\Lambda}{ |x_0 z_1|}) &\leq& \cos ({\Lambda}{ |x_0 z_1|}) \Big(\cos ({\Lambda}{| z_1 z_2|})-1 \Big) \\
&&+ \sin ({\Lambda}{ |x_0 z_1|})\,\sin({\Lambda}{ | z_1 z_2|})\,\cos(\angle x_0 z_1 z_2) \\
&\leq& \frac{\Lambda^2 R}{8} |z_1 z_2| \angle x_0 z_1 z_2.
\end{eqnarray*}
On the other hand,
$$\cos({\Lambda}{ |x_0 z_2|})-\cos ({\Lambda}{ |x_0 z_1|}) \geq \Lambda( |x_0 z_1|-|x_0 z_2|).$$
Hence,
$$2\varepsilon\geq |x_0 z_2|-|x_0 z_1| \geq -\frac{\Lambda R}{8} |z_1 z_2| \angle x_0 z_1 z_2,$$
which yields $\angle x_0 z_1 z_2 \leq 2\pi/3$ due to the condition \eqref{condition-separation-inner}.
Thus combining with \eqref{estimate-eta-inner}, we can choose sufficiently small $\cvdetneighbor$ such that
\begin{equation}
\angle x z_1 z_2 \leq \angle x z_1 x_0 +\angle x_0 z_1 z_2 \leq \frac{5}{6} \pi.
\end{equation}
Then the lemma follows from the exact same argument in the last part of Lemma \eqref{lemma-angle-separation}(1).
\hfill $\square$

\vfill

\medskip



\bigskip
\begin{thebibliography}{99}







\bibitem{Al} G. Alessandrini, \emph{Stable determination of conductivity by boundary measurements}, Appl. Anal. {\bf 27} (1988), 153--172.


\bibitem {AlS} G. Alessandrini, J. Sylvester, \emph{Stability for a multidimensional inverse spectral theorem}, Comm. PDE. {\bf 15} (1990), 711--736.




\bibitem{AA} R. Alexander, S. Alexander, \emph{Geodesics in Riemannian
manifolds-with-boundary}, Indiana Univ. Math. J. {\bf 30} (1981), 481--488.

\bibitem{AKP} S. Alexander, V. Kapovitch, A. Petrunin, 
\emph{Alexandrov geometry: preliminary version no. 1}, arXiv:1903.08539.





\bibitem{AKKLT}
M. Anderson, A.  Katsuda, Y. Kurylev, M. Lassas, M. Taylor,
\emph{Geometric Convergence, and Gel'fand's Inverse Boundary Problem},
Invent. Math. {\bf 158} (2004), 261--321.


%
%
%

%

%
%
%


\bibitem{Be}
M. Belishev,  \emph{An approach to multidimensional inverse problems
for the wave equation}, (Russian) Dokl. Akad. Nauk SSSR {\bf 297} (1987), 524--527



\bibitem{BeKu} 
M. Belishev, Y. Kurylev,
\emph{To the reconstruction of a Riemannian manifold via its spectral
data (BC-method)}, Comm. PDE {\bf 17} (1992), 767--804.



\bibitem{Besson}
P. B\'erard, G. Besson, S. Gallot, \emph{Embedding Riemannian manifolds
by their heat kernel},
Geom. Funct. Anal. {\bf 4} (1994), 373--398.

%
%
%
%

\bibitem{BKL}
R. Bosi, Y. Kurylev, M. Lassas, \emph{Reconstruction and stability in Gel'fand's inverse interior spectral problem}, To appear in Analysis and PDE.


%



\bibitem{Br3}
Y. Brudnyi, P. Shvartsman, {\it A linear extension operator for a space of smooth functions
defined on closed subsets of $\R^n$}, Dokl. Akad. Nauk SSSR {\bf 280} (1985), 268--270.
English transl. in Soviet Math. Dokl. {\bf 31}, No. 1 (1985), 48--51.

\bibitem{Br4}
Y. Brudnyi, P. Shvartsman, {\it Generalizations of Whitney's extension theorem}, Int.
Math. Research Notices {\bf 3} (1994), 129--139.



\bibitem{Burago}
D. Burago, Y. Burago, S. Ivanov,  {\it A course in metric geometry}, Graduate Studies in Mathematics \textbf{33}, AMS, 2001. xiv+415 pp.

\bibitem{BIK}
D. Burago, S. Ivanov, Y. Kurylev, \emph{A graph discretization of the Laplace-Beltrami operator}, J. Spectr. Theory \textbf{4} (2014), no. 4, 675--714. 

\bibitem{BILL} D. Burago, S. Ivanov, M. Lassas, J. Lu,
\emph{Quantitative stability of Gel'fand's inverse boundary problem}, Anal. PDE \textbf{18} (2025), no. 4, 963--1035.


\bibitem{Cha}
I. Chavel,  {\it  Riemannian geometry--a modern introduction.}
Cambridge U.\ Press, 1993.



\bibitem{Ch}
J. Cheeger, \emph{Finiteness theorems for Riemannian manifolds}, Amer. J. Math. {\bf 92} (1970), 61--75.

\bibitem{CY}
J. Cheeger, S. T. Yau, \emph{A lower bound for the heat kernel}, Comm. Pure Appl. Math. \textbf{34} (1981), 465--480.


\bibitem{Chen2019} 
X. Chen, M. Lassas, L. Oksanen, G. Paternain, \emph{Detection of Hermitian connections in wave equations with cubic non-linearity}, To appear in JEMS. arXiv:1902.05711.

\bibitem{cheng}
S. Cheng, T. Dey, E. Ramos,
\emph{Manifold reconstruction from point samples},
SODA (2005), 1018--1027.


\bibitem{CoifmanLafon}
R. Coifman, S. Lafon, \emph{Diffusion maps}, Appl.  Comp. Harm. Anal. \textbf{21} (2006), 5--30.

 \bibitem{diffusion}
R. Coifman, S. Lafon, A. Lee, M. Maggioni, B. Nadler, F. Warner, S. Zucker,
\emph{Geometric diffusions as a tool for harmonic analysis and structure
  definition of data: Multiscale methods},
Proc. of Nat. Acad. Sci. \textbf{102} (2005), 7432--7438.


\bibitem{dH2019}
M. de Hoop, G. Uhlmann, Y. Wang, \emph{Nonlinear responses from the interaction of two progressing waves at an interface}, Ann. de l'Inst. Henri Poincar\'e C, Anal. non lin. \textbf{36} (2019), 347--363.

\bibitem{dH2020} 
M. de Hoop, G. Uhlmann, Y. Wang, \emph{Nonlinear  interaction  of  waves  in  elastodynamics  and  an  inverse problem}, Math. Ann. \textbf{376} (2020), 765--795. 


\bibitem{dCristo-Rondi} 
M. Di Cristo, L. Rondi, S. Vessella, 
\emph{Stability properties of an inverse parabolic problem with unknown boundaries.} Mat. Pura Appl. (4) {\bf 185} (2006), 223--255.


  \bibitem{dsFKSU}
 D. Dos Santos Ferreira, C. Kenig; J. Sj\"strand, G. Uhlmann,
 \emph{Determining a magnetic Schr\"odinger operator from partial Cauchy data.} Comm. Math. Phys. \textbf{271} (2007),  467--488.
     

\bibitem{F1} 
C. Fefferman, \emph{A sharp form of Whitney's extension theorem}, Ann. of Math. {\bf 161} (2005), 509--577.

\bibitem{F2} 
C. Fefferman, \emph{Whitney's extension problem for $C^m$}, Ann. of Math. {\bf 164} (2006), 313--359.

\bibitem{F3} 
C. Fefferman, \emph{$C^m$-extension by linear operators}, Ann. of Math. {\bf 166} (2007), 779--835.

\bibitem{FIKLN}
C. Fefferman, S. Ivanov, Y. Kurylev, M. Lassas, H. Narayanan, \emph{Reconstruction and interpolation of manifolds I: The geometric Whitney problem}, Found. Comp. Math. \textbf{20} (2020), 1035--1133.

\bibitem{FILN}
C. Fefferman, S. Ivanov, M. Lassas, H. Narayanan: \emph{Reconstruction of a Riemannian manifold from noisy intrinsic distances}, SIAM J. Math. Data Science \textbf{2} (2020), No. 3, 770--808.


\bibitem{FK1} 
C. Fefferman, B. Klartag, \emph{Fitting $C^m$-smooth function to data I}, Ann. of Math. {\bf 169} (2009), 315--346.


\bibitem{FK2} 
C. Fefferman, B. Klartag, \emph{Fitting $C^m$-smooth function to data II}, Rev. Mat. Iberoam. {\bf 25} (2009), 49--273.


\bibitem{FMN} 
C. Fefferman, S. Mitter, H. Narayanan, \emph{Testing the manifold
hypothesis}, JAMS \textbf{29} (2016), 983--1049.


\bibitem{Ge} I. Gel'fand, \emph{Some aspects of functional analysis
and algebra}, Proc. Intern. Cong. Math. {\bf 1} (1954),
253--277.


%
%
%

\bibitem{Gr2} M. Gromov, \emph{Filling Riemannian manifolds}, J. Diff. Geom. {\bf 18} (1983),  1--147.



%
%
%
%
%
%
%
%
%
%
%
%
%
%
%
%
%
%
%
%
%
%
%
%
%
%
%
%
%
%
%
%
%
%
%
%

      \bibitem{GT}
   C. Guillarmou, L. Tzou, \emph{Calderon inverse problem with partial data on Riemann surfaces.} Duke Math. J. \textbf{158} (2011), 83--120.
    

\bibitem{Helin-etal}
T. Helin, M. Lassas, L. Oksanen, T. Saksala, \emph{Correlation based passive imaging with a white noise source}, J.  Math. Pures et Appl. {\bf 116} (2018), 132--160.


\bibitem{Helin-etal2}
T. Helin, M. Lassas, L. Ylinen, Z. Zhang, \emph{Inverse problems for heat equation and space-time fractional diffusion equation with one measurement.} J. Diff. Eq. \textbf{269} (2020), no. 9, 7498--7528.

\bibitem{HUZ20} 
P. Hintz, G. Uhlmann, J. Zhai, \emph{An inverse boundary value problem for a semilinear wave equation on Lorentzian manifolds}, 
 Int. Math. Res. Not., rnab088, 2021.



\bibitem{HUZ21} 
P. Hintz, G. Uhlmann, J. Zhai, \emph{The Dirichlet-to-Neumann map for a semilinear wave equation on Lorentzian manifolds}, arXiv:2103.08110.

\bibitem{Ivanov}
S. Ivanov, \emph{Distance difference representations of Riemannian manifolds}, Geom. Dedicata \textbf{207} (2020), 167--192.



\bibitem{Karcher}
H. Karcher, \emph{Riemannian Comparison Constructions}. In
Global differential geometry (Ed. S. S. Chern),
Studies in Mathematics, Vol. 27 (1989) MAA., pp. 170--222.


\bibitem{Kasue}
A. Kasue, \emph{Convergence of Riemannian manifolds and Laplace
operators. I}, Ann. Inst. Fourier {\bf 52}  (2002), 1219--1257.



\bibitem{KaKu} A. Katchalov, Y. Kurylev,
\emph{Multidimensional inverse problem with incomplete boundary spectral
data}, Comm. Part. Diff. Eq. {\bf 23} (1998),  55--95.




\bibitem{KKL} A. Katchalov, Y. Kurylev, M. Lassas,
{\it Inverse Boundary Spectral Problems}. Chapman Hall/CRC, Pure
and Applied Mathematics \textbf{123}. (2001), 290pp.





%

\bibitem{KatsudaKL}
 A. Katsuda, Y. Kurylev, M. Lassas, \emph{Stability of boundary distance representation and reconstruction of Riemannian manifolds}, Inverse Probl. Imag. {\bf 1} (2007), 135--157. 


%
%
%
%
%
%
 
 \bibitem{KSU}
C. Kenig; J. Sj\"strand, G. Uhlmann,
\emph{The Calder\'on problem with partial data}. 
Ann. of Math. {\bf 165} (2007), 567-591.

 \bibitem{KV1}
 R. Kohn, M. Vogelius, \emph{Determining conductivity by boundary measurements.} Comm. Pure Appl. Math.{\bf  37} (1984),289-298.

 \bibitem{KV2}
 R. Kohn, M. Vogelius, \emph{Determining conductivity by boundary measurements. II. Interior results}. Comm. Pure Appl. Math. {\bf 38} (1985),  643-667. 


\bibitem{KrKuLa}
K. Krupchyk, Y. Kurylev, M. Lassas, \emph{Inverse spectral problems on a closed manifold}, J. Math. Pures  Appl. \textbf{90} (2008), 42--59. 

 \bibitem{KrU}
 K. Krupchyk, G. Uhlmann, \emph{Inverse problems for advection diffusion equations in admissible geometries.} Comm. Part. Diff. Eq. \textbf{43} (2018), 585--615.

\bibitem{Ku3} Y. Kurylev,
\emph{Inverse boundary problems on Riemannian manifolds}, Contemp.
Math. {\bf 173} (1994), 181--192.


\bibitem{Dirac} Y. Kurylev, M. Lassas,
  \emph{Inverse problem
  for a Dirac-type equation on a vector bundle}, Adv. Math. {\bf 221} (2009), 170-216.


\bibitem{KLOU2014} 
Y. Kurylev, M. Lassas, L. Oksanen, G. Uhlmann, \emph{Inverse  problem  for  Einstein-scalar  field  equations}.
To appear in Duke Math. J.
%
%
%



\bibitem{KLS} Y. Kurylev, M. Lassas, E. Somersalo,
\emph{Maxwell's equations with a polarization independent wave velocity: direct and inverse problems},
J. Math. Pures Appl. \textbf{86} (2006), 237--270.

\bibitem{Ku2} 
Y. Kurylev, M. Lassas, G. Uhlmann, \emph{Rigidity of broken geodesic flow and inverse problems}, 
Amer. J. Math. {\bf 132} (2010), 529--562.


\bibitem{KLU18} 
Y. Kurylev, M. Lassas, G. Uhlmann,  \emph{Inverse problems for Lorentzian manifolds and non-linear hyperbolic equations}, Invent. Math. \textbf{212} (2018), 781--857. 


\bibitem{KOP}
Y. Kurylev, L. Oksanen, G. Paternain, \emph{Inverse problems for the connection Laplacian}, J. Diff. Geom. \textbf{110} (2018), 457--494. 


\bibitem{Lassas}
M. Lassas, \emph{Inverse problems for linear and non-linear hyperbolic equations.} Proc. Int. Congress of Math. ICM 2018, Rio de Janeiro, Brazil, Vol III, 3739--3760, 2018.


  \bibitem{LO}   
M.    Lassas, L.  Oksanen,  \emph{Inverse problem for the Riemannian wave equation with Dirichlet data and Neumann data on disjoint sets.} Duke Math. J. \textbf{163} (2014), 1071--1103. 

\bibitem{SakL}
M. Lassas, T. Saksala, \emph{Determination of a Riemannian manifold from the distance difference functions}. Asian J. Math. \textbf{23} (2019), 173--200.


\bibitem{LTU}
M. Lassas, M. Taylor, G. Uhlmann,
\emph{The Dirichlet-to-Neumann map for complete Riemannian manifolds with
boundary}, Comm. Anal. Geom. {\bf  11}  (2003), 207--221.



\bibitem{LU}
M.  Lassas, G. Uhlmann,  \emph{On determining a Riemannian manifold
from the Dirichlet-to-Neumann map}, Ann. Sci. Ecole Norm.
Sup. {\bf  34}  (2001), 771--787.

\bibitem{LUW18} 
M. Lassas, G. Uhlmann, Y. Wang, \emph{Inverse Problems for Semilinear Wave Equations on Lorentzian Manifolds}, Comm. in Math. Phys. \textbf{360} (2018), 555--609.

\bibitem{LY}  
P. Li, S. T. Yau, \emph{On the parabolic kernel of the Schr\"odinger operator}, Acta Math. \textbf{156} (1986), 153--201.

%
%
%
%
%
%
%
%
%

\bibitem{Shape1} 
F. Memoli. \emph{Spectral Gromov-Wasserstein distances for
shape matching.} In Proc. NORDIA, 2009.

\bibitem{Shape2} 
M. Ovsjanikov et al., \emph{One Point Isometric Matching with the Heat Kernel}.
  Computer Graphics Forum 29 (2010):1555-1564.

\bibitem{UP}
L. Pestov, G. Uhlmann, \emph{Two dimensional compact simple Riemannian
manifolds are boundary distance rigid}, Ann. of Math. {\bf 161} (2005), 1093--1110.







\bibitem{Pe}
P. Petersen, {\it Riemannian geometry}, 1st Ed. Springer, 1998. xvi+432pp.

\bibitem{Plaut} C. Plaut, 
\emph{Metric spaces of curvature $\geq k$}, Handbook of geometric topology, 819--898, North-Holland, Amsterdam, 2002. 


\bibitem{Portegies} J. Portegies, \emph{Embeddings of Riemannian manifolds with heat kernels and eigenfunctions}, Comm. Pure Appl. Math. \textbf{69} (2016), 478--518.





\bibitem{RS} 
S. Roweis, L. Saul, \emph{Nonlinear dimensionality reduction by locally
linear embedding}, Science \textbf{290} (2000), 2323--2326.



\bibitem{SY}
R. Schoen, S. T. Yau, \emph{Lectures on differential geometry}, International Press, 1994.


\bibitem{Shiohama}
K. Shiohama, \emph{An introduction to the geometry of Alexandrov spaces}, Lecture Notes Series, 8. Seoul National University, Seoul, 1993. ii+78 pp.


%
%



\bibitem{StU1}
P. Stefanov, G. Uhlmann, \emph{Stability estimates for the hyperbolic
Dirichlet to Neumann map in anisotropic media}, J. Funct.
Anal. {\bf  154}  (1998), 330--358.







\bibitem{StU2}
P. Stefanov, G. Uhlmann, \emph{Boundary rigidity and stability for
generic simple metrics}, J. Amer. Math. Soc. {\bf 18}
(2005), 975--1003.







\bibitem{SHoop}
C. Stolk, M. de Hoop,  \emph{Microlocal analysis of seismic inverse
scattering in anisotropic elastic media}, Comm. Pure Appl.
Math. {\bf  55}  (2002), 261--301.




\bibitem{TSL} 
J. Tenenbaum, V. de Silva,  J. Langford, \emph{A global geometric framework for nonlinear dimensionality reduction}, Science \textbf{290} (2000),  2319--2323.

\bibitem{U}
G. Uhlmann, \emph{Inverse boundary value problems for partial
differential equations}, Proceedings of the International
Congress of Mathematicians, Vol. III (Berlin, 1998). Doc. Math., 77--86.




%
%


\bibitem{UhWa18} 
G. Uhlmann, Y. Wang. \emph{Determination of space-time structures from gravitational perturbations}, To appear in Comm. Pure App. Math.
%

\bibitem{Varadhan}
S. Varadhan,  \emph{On the behaviour of the fundamental solution of the heat equation with variable coefficients.} Comm. Pure Appl. Math. {\bf 20} (1967), 431--455.

\bibitem{Varopoulos}
N. Varopoulos, \emph{The Poisson kernel on positively curved manifolds}, J. Funct. Anal. \textbf{44} (1981), 359--380.




\bibitem{WZ2019} 
Y. Wang, T. Zhou, \emph{Inverse problems for quadratic derivative nonlinear wave equations}, Comm. PDE \textbf{44} (2019), 1140--1158.

\bibitem{WZ} 
X. Wang, K. Zhu, \emph{Isometric embeddings via heat kernel}. J. Diff. Geom. {\bf 99} (2015), 497--538.



%

\bibitem{W12} H. Whitney, {\it Functions differentiable on the boundaries of regions,} Ann. of Math. {\bf 35}
(1934), 482--485.



%

\bibitem{Zelditch}
S. Zelditch, \emph{Survey on the inverse spectral problem.} ICCM Not. 2 (2014), no. 2, 1-20. 

\bibitem{Zha}
H. Zha, Z. Zhang,
\emph{Continuum Isomap for manifold learnings},
Comp. Stat. Data Anal. \textbf{52} (2007), 184--200.


\end{thebibliography}
\end {document}